\documentclass[11pt]{article}
\newcommand{\mylabel}[1]{\label{#1}}
\usepackage{amsthm}
\usepackage{amsfonts}
\usepackage{amscd}
\usepackage{amsmath}
\usepackage{amssymb}
\usepackage{color}
\usepackage{array}
\usepackage{epsfig}
\usepackage{verbatim}

\def\cbar{\overline{\C}}

\def\al{\alpha}
\def\be{\beta}

\def\g{\gamma}
\def\G{\Gamma}

\def\ep{\varepsilon}
\def\la{\lambda}
\def\La{\Lambda}
\def\De{\Delta}
\def\de{\delta}

\def\wtq{\widehat S}

\def\R{\mbox{$\mathbb R$}}

\def\C{\mbox{$\mathbb C$}}

\def\D{\mbox{$\mathbb D$}}
\def\ds{\displaystyle}
\def\P{\mbox{$\mathcal P$}}
\def\Q{\mbox{$\mathcal Q$}}
\def\S{\mbox{$\mathbb S$}}
\def\N{\mbox{$\mathbb N$}}

\def\AAA{{\mathcal A}}
\def\A{{\mathcal A}}
\def\BBB{{\mathcal B}}
\def\CCC{{\mathcal C}}

\def\E{{\mathcal E}}
\def\FFF{{\mathcal F}}
\def\GGG{{\mathcal G}}

\def\JJJ{{\mathcal J}}
\def\LLL{{\mathcal L}}
\def\MMM{{\mathcal M}}

\def\OOO{{\mathcal O}}

\def\RRR{{\mathcal R}}

\def\XXX{{\mathcal X}}

\newtheorem{newthm}{Theorem}
\newtheorem{theorem}{Theorem}[section]
\newtheorem{lemma}[theorem]{Lemma}
\newtheorem{proposition}[theorem]{Proposition}
\newtheorem{corollary}[theorem]{Corollary}
\newtheorem{defthm}[theorem]{Definition et \th}
\newcommand{\REFEQN}[1] { \begin{equation}\mylabel{#1} }
\newcommand{\ENDEQN}{\end{equation}}
\newcommand{\REFTHM}[1] { \begin{theorem}\mylabel{#1} }
\newcommand{\ENDTHM}{\end{theorem}}
\newcommand{\REFNTH}[1] { \begin{newthm}\mylabel{#1} }
\newcommand{\ENDNTH}{\end{newthm}}
\newcommand{\REFPROP}[1]{\begin{proposition}\mylabel{#1} }
\newcommand{\ENDPROP}{\end{proposition} }
\newcommand{\REFLEM}[1]{\begin{lemma}\mylabel{#1} }
\newcommand{\ENDLEM}{\end{lemma} }
\newcommand{\REFCOR}[1]{\begin{corollary}\mylabel{#1} }
\newcommand{\ENDCOR}{\end{corollary} }

\newcommand{\REFDEFTHM}[1] { \begin{defthm}\mylabel{#1} }
\newcommand{\ENDDEFTHM}{\end{defthm}}
\newcommand{\Ref}[1]{ (\ref{#1}) }
\def\beginp{ {\noindent\em Proof.} }
\def\smm{ {\smallsetminus }}
\def\qr{quasi-regular}
\def\qc{quasi-conformal}
\def\rc{bordered Riemann surface}
\def\ps{puzzle surface}
\def\ce{c-equivalent}
\def\rs{repelling system}
\def\cm{of constant complexity}
\def\gf{geometrically finite}
\def\pf{postcritically finite}
\def\sr{semi-rational map}

\begin{document}

\begin{center}{\Large\bf A characterization of hyperbolic rational maps \vspace{0.1cm} } \\
CUI Guizhen and TAN Lei,\ \today\end{center}

\begin{abstract}
In the early 1980's Thurston gave a topological  characterization of  rational maps  whose critical points
have finite iterated orbits (\cite{Th,DH1}): given a topological branched covering $F$ of the two sphere with finite critical orbits, if $F$ has no Thurston obstructions then $F$ possesses an invariant complex structure (up to isotopy), and is combinatorially equivalent to a
rational map.

We extend this theory to the
setting of rational maps  with infinite critical orbits, assuming a certain kind of hyperbolicity.
Our study includes also  holomorphic dynamical systems that arise as coverings over disconnected Riemann surfaces of finite type. The obstructions we encounter are similar to those of Thurston. We give concrete criteria
for verifying whether or not such obstructions exist.

Among many possible applications, these results can be used for example to
construct holomorphic maps with prescribed dynamical properties; or to give a
parameter description, both local and global, of bifurcations of complex dynamical systems.
\end{abstract}

Subject class [2000]: {Primary 37F; Secondary 32G}.
\section{Introduction}

Thurston's characterization of \pf\ rational maps is one of the
major tools in complex dynamics. It enables us to produce various
kind of rational maps with prescribed dynamical properties, as well as to produce combinatorial models for parameter spaces. There are many applications of
Thurston's theory. Just to mention a few, we may cite for example
Douady's proof of monotonicity of entropy for unimodel maps (\cite{Do}),
Mary Rees' descriptions of parameter spaces (see \cite{Re}),
McMullen's rational quotients (\cite{Mc2}), Kiwi's
characterization of polynomial laminations (using previous work of
Bielefield-Fisher-Hubbard (\cite{BFH}) and Poirier (\cite{Po})),
etc.

One drawback of Thurston's theorem is that it can only be applied
to \pf\ rational maps. On one hand, these maps all have a connected Julia set;
on the other hand, they form a totally
disconnected subset in the parameter space (except the Latt\`es examples). Therefore the theorem can not
characterize the combinatorics of  disconnected Julia sets,
nor the bifurcations through continuous parameter perturbations.

Over the years, there has been several attempts to extend
Thurston's theory beyond \pf\ maps. For example, David Brown (see
[Br]), supported by previous work of Hubbard and Schleicher
([HS]), has succeeded in extending the theory to the uni-critical
polynomials with an infinite postcritical set (but always with a
connected Julia set), and pushed it even further to the infinite
degree case, namely the exponential maps. See also Jiang and Zhang
[JZ].

We mention also a recent work of Hubbard-Schleicher-Shishikura \cite{HSS}
extending Thurston's theorem to \pf\ exponential maps.

In this paper, supported by previous works of Cui, Jiang and
Sullivan ([CJS]), as well as unpublished manuscripts of Cui, we
extend Thurston Theorem to the full setting of arbitrary
non-postcritically finite hyperbolic or sub-hyperbolic rational
maps. Our analysis leads naturally to the concept of \rs s over disconnected
Riemann surfaces of finitely type, and allows us to establish an analog of Thurston Theorem for these dynamical systems.

This work consists of the first step of a long program,
as exposed in [C2], to the study of deformations and bifurcations
of rational maps. In a forthcoming paper (\cite{CT}), we will
extend our characterization to the setting of \gf\ rational maps
(i.e.  maps with parabolic periodic points), and then give a
detailed study of their relations with hyperbolic rational maps. A
\gf\ map $g$ often sits on the boundary of several hyperbolic
components, and does so in a quite subtle way: if you approach it
algebraically, you may or may not get an different geometric
limit, depending very much how you approach it. This subtlety
makes the study for the deformation of $g$ very difficult.
However, it is relatively easy to describe combinatorially all the
possible bifurcations. Then, equipped with our Thurston-like
realization result, we will be able to prove easily the existence
of such bifurcations. For instance we will classify all the
hyperbolic components $H$ that contain a path converging to $g$
and that along the path the algebraic and geometric limits
coincide. Conversely, given a hyperbolic component $H$,  we will
apply our technique to determine all the boundary \gf\ maps $g$
that are path-accessible from $H$ with the same properties.

{\noindent\bf Statements}

All branched coverings, homeomorphisms in this paper are
orientation preserving. Let $G:\cbar\to \cbar$ be an orientation preserving branched
covering with degree $\deg G\ge 2$. Its {\em postcritical set} is
defined to be
$$\P_G:=\text{closure}\{G^n(c)|\ n>0, c\text{ a critical
point of }G \}.$$ Denote by $\P_G'$ the accumulation set of
$\P_G$.

We say that $G$ is {\em \pf}\ if every critical point has a finite
orbit (i.e. $\P'_G=\emptyset$). We say that $G$ is a  {\bf
sub-hyperbolic \sr}\ if $\P_G'$ is finite (or empty); and in case
 $\P_G'\ne \emptyset$,
the map $G$ is holomorphic in a neighborhood of $\P_G'$ and every
periodic point in $\P_G'$ is either attracting or
super-attracting.

Two sub-hyperbolic semi-rational maps $G_1$ and $G_2$ are called
{\em \ce}, if there is a pair $(\phi,\psi)$ of homeomorphisms of
$\cbar$, and a neighborhood $U_0$ of $\P_{G_1}'$ such that: \begin{enumerate}\item[
(a)] $\phi\circ G_1=G_2\circ \psi$; \item[
(b)] $\phi$ is holomorphic in $U_0$; \item[
(c)] the two maps $\phi$ and $\psi$ are equal on $\P_{G_1}$, thus on $\P_{G_1}\cup \overline{U_0}$ (by the isolated zero theorem);\item[
(d)] the two maps $\phi$ and $\psi$ are isotopic to each other rel $P_{G_1}\cup\overline{ U_0}$, i.e., there is a continuous map $H:[0,1]\times \cbar\to \cbar$ such that each $H(t,\cdot)$
is a homeomorphism of $\cbar$, $H(0,\cdot)=\phi$, $H(1,\cdot)=\psi$, and $H(t,z)\equiv \phi(z)$
for any $t\in [0,1]$ and any $z\in \P_{G_1}\cup \overline {U_0}$.\end{enumerate}

Given a sub-hyperbolic semi-rational map $G$, we consider the
problem of whether  there is a rational map \ce\ to it.

Thurston gave a combinatorial criterion of the same problem for
\pf\ branched coverings, based  on the absence of Thurston
obstructions (see \S \ref{obstructions} and Theorem \ref{Thurston}
below). We prove here:

\REFTHM{hyperbolic} Let $G$ be a sub-hyperbolic \sr\ with
$\P_G'\ne \emptyset$. Then  $G$ is \ce\ to a rational map $g$ if
and only if $G$ has no Thurston obstruction. In this case the
rational map $g$ is unique up to M\"obius conjugation. \ENDTHM

The necessity of having no Thurston obstruction, and
the unicity of the rational map $g $, are known to be true for a
wider class of maps. See \cite{Mc2} (or Theorem \ref{mcm} below)
and  [C1].

Thus it remains only to prove the existence part here: i.e. to
show that if $G$ is unobstructed then it is c-equivalent to a
rational map.

In the process of proving the theorem, we introduce the concept of repelling systems
over disconnected Riemann surfaces of finite type,  and those \cm. We develop
a corresponding Thurston-like theory, including the notions of c-equivalence
to holomorphic models, Thurston obstructions, and then a theorem
saying that such a system without obstructions is c-equivalent to a holomorphic model
(see Theorems \ref{more-general} and \ref{sb} for detailed statements).

The general strategy of the proof of Theorem \ref{hyperbolic} can be then described
as follows:
we define $K_G$, its {\em filled Julia set} relative to $\P_G'$,
to be the set of points not attracted by the cycles in $\P'_G$,
i.e. \REFEQN{KK} K_G:=\{z\in \cbar\ |\
\overline{\bigcup_{n>0}\{G^n(z)\}}\cap \P_G'=\emptyset\}\ .\ENDEQN

Step 0. We show  that up to a change of representatives in the
c-equivalence class of $G$, we may assume that $G$ is
quasi-regular (Lem. \ref{BB}).

Given now $G:\cbar\to \cbar$ an unobstructed \qr\ sub-hyperbolic semi-rational map.
\begin{enumerate}
\item there is a restriction $G|_{L_1}:L_1\to L_0$ in a neighborhood of $K_G$ which is an unobstructed \rs\ (\text{Lem.} \ref{un-obstructed}).
\item there is
a sub-\rs\ $F$ of $G|_{L_1}$ that is both unobstructed and \cm \ (\text{Thms.} \ref{Unobstructed}.(B),\ref{Repellor})
\item
This \rs\  \cm\  has no boundary obstructions nor renormalized obstructions (\text{Lem.} \ref{Prepare}).
\item
Any $F$ with properties  in Step 3 is c-equivalent to a holomorphic model (\text{Thm.} \ref{sb}).
\item
 $G|_{L_1}:L_1\to L_0$ is c-equivalent to a holomorphic model (\text{Thm.} \ref{Unobstructed}.(A)).
 \item  $G:\cbar\to \cbar$ is c-equivalent to a rational map (\text{Prop.} \ref{final}).
\end{enumerate}

Steps 1-3 consist of detailed study of Thurston obstructions for repelling systems, as well as combinatorics of puzzle neighborhoods of $K_G$.

Step 4 (Theorem \ref{sb}) is the core part of this work. It is proved using Gr\"otzsch inequalities, and, in the presence of renormalizations,
Thurston's original Theorem, together with a form of reversed Gr\"otzsch inequality.

Steps 5-6 are standard applications of Measurable Riemann Mapping Theorem.

Steps 2-5 together lead to a Thurston-like theorem for repelling systems
(see Theorem \ref{more-general} for a precise statement), which is of independent interest.

Notice that we do not take the approach, as one might have attempted to do, of iteration in an
infinite dimensional Teichm\"uller space.

{\noindent\bf Remarks.}
Actually the no-obstruction condition in both our theorem \ref{hyperbolic} and
Thurston's original one, is in general a difficult condition to
verify, as there are infinitely many candidate obstructions.
Therefore in order to apply them effectively, further efforts are
often needed. In the \pf\ setting, many methods have been
developed to overcome this difficulty. But in the case at hand,
our result would have been left unsatisfactory if no further
criteria have been given. Fortunately, what we have actually
proved (see Theorems \ref{sb} and \ref{infinite} below) does provide more effective
criteria. More precisely we will decompose the dynamics into
several renormalization pieces that are in fact \pf\ maps,
together with a transition matrix that records the gluing data.
This decomposition is not entirely trivial and presents some
interests even for rational maps. Our proof shows then that in order to
be \ce\ to a rational map it amounts only to check Thurston's
condition for the renormalizations (thus back to the \pf\ setting),
and for the gluing data, which is only one eigenvalue to
calculate.

We obtain also a combination result that is very practical to use.
We will show in Theorem \ref{infinite} that for any finite
collection $f_i$ of rational maps with connected Julia set $J_i$
(\pf\ or not), together with a compatible (unobstructed) gluing
data $D$, one can glue the $f_i$'s on neighborhoods of $J_i$
together following $D$, to obtain a rational map $g$, so that each
$f_i$ appears as a renormalization of $g$.

For a similar decomposition-gluing approach, we recommend
\cite{Pi}. The topology of Julia components for hyperbolic
rational maps has been well understood. See \cite{PT}.

Theorem \ref{hyperbolic} was already announced in [CJS], together
with a sketch of the main ideas of the proof. Numerous details
there were however missing, and sometimes erroneous. The
presentation here will be totally different. In particular the concept of repelling systems and the related Thurston-like theory are new.
This will lead also to two stronger
and easier to use results: Theorems \ref{sb} and \ref{infinite}.

Along the proof we will provide  numerous supporting diagrams and
pertinent examples.

{\noindent\bf Organization}.

The paper is organized as follows: In \S \ref{Conclusion} we prove
Step 0 and Step 6 above. We introduce the concept of repelling system and show how it appears as a restriction near $K_G$ of a global map $G$.

In \S \ref{repel} we first recall the definition of Thurston obstructions and state
Thurston's original theorem. We then develop the corresponding concepts for repelling systems and state a Thurston-like theorem in this setting (Theorem \ref{more-general}).
Assuming it we prove Theorem \ref{hyperbolic} (we just need to do Step 1 above).

In \S \ref{againagain}
we introduce the concepts of constant complexity
repelling systems and the specific obstructions associated to them. We state our Thurston-like
theorem, Theorem \ref{sb}, in this setting.
Assuming this we complete Steps 2-5 above and proves Theorem \ref{more-general}.

In \S \ref{bridge}-\ref{I} we give the proof of Theorem \ref{sb}.

In the final section \S \ref{comb} we state
Theorem \ref{infinite}.

{\noindent\bf Acknowledgment}. We would like to express our thanks
to Xavier Buff, John Hubbard, Yunping Jiang,  Carsten Lunde Petersen, Kevin Pilgrim, Mary
Rees, Mitsuhiro Shishikura and Dennis Sullivan for helpful
comments.

\section{Reducing to restrictions near $K_G$}\label{Conclusion}

Let $G:\cbar\to\cbar$ be a sub-hyperbolic \sr\ with $\P_G'\ne
\emptyset$, i.e. $G$ is a branched covering such that the cluster
set $\P_G'$  of its postcritical set is finite and non-empty, that
$G$ is holomorphic in a neighborhood of $\P_G'$, and that every
periodic cycle in $\P_G'$ is either attracting or superattracting.

Our objective is to show that if $G$ has no Thurston obstruction
then $G$ is \ce\ to a certain rational map.

\subsection{Making the map quasi-regular}

\REFLEM{BB} Let $G$ be a sub-hyperbolic \sr\ with
$\P_G'\ne\emptyset$. Then $G$ is c-equivalent to a \qr\
sub-hyperbolic \sr. \ENDLEM

\beginp
Consider $G$ as a branched covering from $\cbar$ onto $\cbar$.
There is a unique complex structure $\XXX'$ on $\cbar$ such that
$G: (\cbar,\XXX')\to \cbar$ is holomorphic (see \cite{DD}, section
6.1.10). The uniformalization theorem provides thus a conformal
homeomorphism $\xi: (\cbar,\XXX')\to \cbar$. Set $R:=G\circ
\xi^{-1}$. Then $R:\cbar\to \cbar$ is a branched covering,
holomorphic with respect to the standard complex structure,
therefore a rational map.

Let $U\subset\cbar$ be a finite union of quasi-discs with mutually disjoint closures, such that $\P'_G\subset U$,
$G^{-1}(U)\supset\overline{U}$ and $G$ is holomorphic in a
neighborhood of $\overline U$ with respect to the standard complex
structure. Then the new structure $\XXX'$ is compatible with the
chart $G^{-1}(U)$, so $\xi(U)$ is also a finite disjoint union of
quasi-discs. Set $L:=\cbar\smm U$. Then there is a quasi-conformal
homeomorphism $\eta:L\to \xi(L)$ such that $\eta=\xi$ on
$(\partial L)\cup (\P_G\cap L)$ and $\eta$ is isotopic to $\xi$
rel $(\partial L)\cup (\P_G\cap L)$ (see Lemma \ref{open}). Set
$\zeta=\eta^{-1}\circ \xi$ on $L$ and $\zeta=id$ on $U$. Then
$\zeta$ is isotopic to the identity rel $\overline{U}\cup\P_G$. So
$G\circ \zeta^{-1}$ is c-equivalent to $G$. But $G\circ
\zeta^{-1}=R\circ \eta$ on $L$, with $\eta$ quasi-conformal and
$R$ holomorphic. One sees that $G\circ \zeta^{-1}$ is
quasi-regular in $L$, thus quasi-regular on $\cbar$. \qed

\subsection{Repelling system as restriction}\label{relation}

{\noindent\bf Definition 1}.  For two subsets $E_1,E_2$ of
$\cbar$, we use the symbol $E_1\subset\subset E_2$ if the closure
of $E_1$ is contained in the interior of $E_2$. We use also
$E^c:=\cbar\smm E$ to denote the complement of $E$ in $\cbar$. If
$E\subset \cbar$ is compact connected, each complementary
component of $E$ (i.e. a component of $E^c$) is a disc.

We will cover $K_G$ by a suitable puzzle $\LLL$ so that in
particular $G^{-1}(\LLL)\subset\subset\LLL$, just as in
Branner-Hubbard's study of cubic polynomials with disconnected
Julia set (see \cite{BH}).  The restriction $G|_{G^{-1}(\LLL)}$,
considered as a dynamical system, leads naturally to the concept
of repelling systems. Such dynamical systems can be also
considered as generalization of Douady-Hubbard's polynomial-like
mappings (\cite{DH3}) in two aspects: the domain of definitions
will be several components each with several boundary curves (this
is necessary as we are dealing with rational maps), and the
dynamics will be quasi-regular branched coverings without
necessarily analyticity.

{\noindent\bf Definition 2}. We say that  $S\subset \cbar$ is a
{\em (quasi-circle) \rc} if it is
$\cbar$ minus finitely many (might be zero) open
quasi-discs $D_l$ so that their closures $\overline{D_l}$ are mutually disjoint.
Therefore $\partial S$ is a disjoint union of finitely many (might be zero) quasi-circles.
We say that $\LLL=S_1\sqcup\cdots \sqcup S_k$ is a
{\em \ps}, if each $S_i$ a \rc, and two distinct $S_i$ and $S_j$
are either contained in distinct copies of the Riemann sphere, or are mutually disjoint. Each $S_i$ is also
called an {\em $\LLL$-piece}.

{\noindent\bf Definition 3}. We say that a map $F:\E\to \LLL$
is a {\bf (\qr) \rs}, if:\\
 $\E\subset\subset\LLL$ are two nested \ps s and $F:\E\to \LLL$ is a \qr\ proper mapping;\\
more precisely if
$\LLL=S_1\sqcup\cdots \sqcup S_k$ is a  \ps, $\E:=\bigcup_{i,j\in \{1,\cdots, k\}, \de\in \La_{ij} }E_{ij\de}$
and $F:\E\to \LLL$ is a map such that, \begin{enumerate}
\item[(TOP)] for each pair $(i,j)$, the index set $\La_{ij}$ is finite or empty, and,  for any $\de\in \La_{ij}$,
the set $E_{ij\delta}$ is a \rc\ compactly contained in the interior of $ S_i$.
Furthermore for each $i$ the $E_{ij\delta}$'s
are mutually disjoint for all possible choices of $j$ and $\delta$.
\item[
(DYN)] $F|_{E_{ij\de}}:E_{ij\delta}\to S_j$ is a quasi-regular proper mapping for all possible $i,j$ and $\delta$.\end{enumerate}

This notion generalizes Douady-Hubbard's polynomial-like maps in both the complexity of the
domains and the regularity of the map (recall that $F:\E\to \LLL$ is polynomial-like if
both $\E$ and $\LLL$ are connected and simply connected hyperbolic Riemann surfaces with
$\E\subset\subset \LLL$ and if $F$ is an analytic proper mapping). Notice that due to the
disconnectedness of the domains the number of preimages $\#F^{-1}(b)$ need not be constant when we let $b$ vary in $\LLL$.

To a \rs\ $F:\E\to \LLL$ we associate its
{\bf postcritical set $P_F$} and {\bf filled Julia set $K_F$} by:
$$\P_F:=\text{closure}\{F^n(c)\in \LLL\ |\ c\text{ a critical
point of }F,\ n>0,\ c,F(c),\cdots,F^{n-1}(c)\in \E \};$$
$$K_F:=\{z\in \E\ | \ F^n(z)\in \E,\ \forall n>0\}. $$

The set $\P_F$ might be empty, in this case we say that $F$ is unbranched. One may construct examples for which $K_F$ is empty (for example $\LLL=S_1\sqcup S_2$, $\E=E_{12}\subset\subset S_1$
and $F(E_{12})=S_2$), although we will be only interested in the case that $K_F\ne \emptyset$, with either $\P_F=\emptyset$ or not. We have
$F^{-1}(K_F)=K_F$, $F(\P_F\cap K_F)\subset\P_F\cap K_F$.

We say furthermore that $F:\E\to \LLL$ is {\em \pf} if $\P_F$ is
finite or empty (this is equivalent to say $\#\P_F\cap K_F<\infty$).
In particular we say that $F:\E\to \LLL$ is an {\em
annuli covering} if each component of $\LLL$ and $\E$ is a closed
annulus, and $F$ is unbranched, i.e. $\P_F=\emptyset$.

A \rs\ $F:\E\to \LLL$ is {\em  holomorphic} if $F$ is holomorphic in
the interior of $\E$.

\medskip
The following {\em restriction principle} provides the most fundamental examples of
the above concepts:

\REFLEM{restriction} Let $G$ be a \qr\ sub-hyperbolic \sr\ with
$\P'_G\ne \emptyset$. Then there exists a \ps\ neighborhood
$L_0$ of $K_G$ such that, setting $L_1=G^{-1}(L_0)$, $$L_1\subset\subset L_0,\ \text{ and }\ G|_{L_1}: L_1\to L_0\ \text{is a \pf\ \rs}.$$ \ENDLEM

\beginp Note that for $S\subset\cbar$ a \rc,
a necessary and sufficient condition for each component of $G^{-1}(S)$ to be
a \rc\ is that $\partial S$ does not contain any critical value of $G$.

One can find an open set $U_0$ which is the union of a quasi-disc
neighborhood for each point of $\P_G'$ so that these quasi-discs
have disjoint closures, that $\partial U_0$ is disjoint from
$\P_G$, that $G$ is holomorphic in a neighborhood of
$\overline{U_0}$, and that $G(U_0)\subset\subset U_0$.

Set $L_0=\cbar\smm U_0$. Then $K_G\subset\subset G^{-1}(L_0)\subset\subset L_0$,
$\P_G\cap L_0$ is finite (or empty). This $L_0$ satisfies the requirement of the lemma.\qed

Note that one may also set $L_n=G^{-n}(L_0)$ for $n\in \N$, to produce a sequence of
\rs s $G|_{L_{n+1}}: L_{n+1}\to L_n$ satisfying the requirement of the lemma.

\medskip

{\noindent \bf Examples}.
\begin{enumerate}
\item[{\bf A.}] Let $\E\subset\subset \LLL$ be two closed
quasi-discs, and $F:\E\to \LLL$ be a holomorphic proper map. Then $F$ is a polynomial-like map in the
sense of Douady-Hubbard, $K_F$ is simply the filled Julia set, and
$\P_F$ is the postcritical set.
\item[
{\bf B.}]  $\LLL$ is a closed quasi-disc, $\E$ is the
union of finitely many disjoint closed quasi-discs contained in
the interior of $\LLL$, and $F$ maps each $\E$-piece
quasi-conformally onto the larger disc $\LLL$. In this case
$\P_F=\emptyset$ and $K_F$ is the non-escaping set of $F$. If $F$
is also holomorphic, the filled Julia set $K_F$ is a Cantor set.
This happens when $F$ is $z^2+c$ for large $c$ and $\LLL$ is a disc
in $\C$ bounded by an equipotential such that $0\in \LLL$ but
$c\not\in \LLL$.
\item
[{\bf C.}] By convention we may consider $\E=\LLL=\cbar$ as
\ps s and a \qr\ \pf\ branched covering $F:\cbar\to\cbar$ as a
\rs.
\end{enumerate}

More important classes of examples are provided by Lemma \ref{restriction}, by the annuli coverings (see below).
See also \S\ref{Examples}.

{\noindent\bf Marking.} Let $F:\E\to \LLL$ be a \pf\ \rs. We say that it is {\em marked},
if it is equipped with a marked set $\P$, satisfying that \\
$$\P_F\subset\P\subset (\LLL\smm \partial \LLL);\ \#\P<\infty \text{ and }\ F(\P\cap \E)\subset\P\ .$$

If not explicitly mentioned, we will consider $F$ to be marked by its postcritical set $\P_F$.

\medskip Motivated  by Thurston's theory, we give the following:

{\noindent\bf Definition 4}. We say that two marked \pf\ repelling systems
$(\E\stackrel{F}{\to}\LLL,\P)$ and $(\E'\stackrel{R}{\to}\LLL',\P')$ are {\bf c-equivalent}, if there is a pair of quasi-conformal
homeomorphisms $\Phi,\Psi:\LLL\to \LLL'$ such that
\REFEQN{C-class}\left\{\begin{array}{l}\Psi(\E)=\E'\ \text{and}\ \Psi(\P)=\P'\\
 \Psi|_{\partial \LLL\cup \P}= \Phi|_{\partial\LLL\cup \P}\vspace{0.1cm} \\
 \Psi \text{ is isotopic to $\Phi$ rel }\partial\LLL\cup \P\\
\Phi\circ F\circ \Psi^{-1}|_{\E'}=R
\end{array}\right. \ \text{ in particular}\ \begin{array}{rcl}
\LLL\supset \E & \overset{\Psi}{\longrightarrow} & \E'\subset \LLL' \vspace{0.09cm} \\
                                    F\downarrow && \downarrow R \vspace{0.09cm}\\
                    \LLL & \underset{\Phi}{\longrightarrow} & \LLL'
                    \end{array}\ \text{commutes}.\ENDEQN
We say that $(\E\stackrel{F}{\to}\LLL,\P)$ is {\bf c-equivalent to a holomorphic model}, if there
is a holomorphic $(\E'\stackrel{R}{\to}\LLL',\P')$ c-equivalent to it.

See \S\ref{Examples} for examples. We have the following criterion:

\REFLEM{5'} A marked \pf\ \rs\ $(\E\stackrel{F}{\to}\LLL,\P)$ is c-equivalent to a holomorphic
model iff there is a pair $(\Theta,\mu)$ such that: \\
(a) $\Theta:\LLL\to \LLL$ is a quasi-conformal homeomorphism with
$\Theta|_{\partial\LLL\cup\P}= id$ and $\Theta$ is isotopic to the
identity  rel $\partial\LLL\cup\P$. \\
(b) $\mu$ is a Beltrami differential on $\LLL$ with
$\|\mu\|_\infty<1$ and $(F\circ
\Theta^{-1})^*(\mu)=\mu|_{\Theta(\E)}$. \ENDLEM

\beginp Assume that $(\E\stackrel{F}{\to}\LLL,\P)$ is c-equivalent to a holomorphic
$(\E'\stackrel{R}{\to}\LLL',\P')$. Let $(\Phi,\Psi)$ be the pair of
quasi-conformal maps given by Definition 4. Set
$\Theta=\Phi^{-1}\circ\Psi$. Then $\Theta$ satisfies the required
isotopic conditions. Let $\mu$ be the Beltrami coefficient of
$\Phi$. Then \Ref{C-class} means exactly that $(F\circ
\Theta^{-1})^*(\mu)=\mu|_{\Theta(\E)}$.

Conversely assume the existence of the pair $(\Theta,\mu)$. By the
Measurable Riemann Mapping Theorem, there is a quasi-conformal map
$\Phi$ defined piecewisely on $\LLL$ with Beltrami coefficient $\mu$. Set $\Psi= \Phi\circ \Theta$. Then for $\E':=\Psi(\E)$, $\P'=\Psi(\P)$, $\LLL':=\Phi(\LLL)$ and $R:=\Phi\circ
F\circ \Psi^{-1}:\E'\to \LLL'$, we know
that $(\E'\stackrel{R}{\to}\LLL',\P')$ is a holomorphic \rs\ c-equivalent to
$(\E\stackrel{F}{\to}\LLL,\P)$. \qed

The following result relates repelling systems to our main
interest (Theorem \ref{hyperbolic}) through restriction:

\REFPROP{final} Let $G$ be a quasi-regular sub-hyperbolic \sr\ with
$\P_G'\ne \emptyset$. If there is a \ps\ neighborhood $\LLL$ of $K_G$
with $G^{-1}(\LLL)\subset\subset \LLL$ and $\partial \LLL\cap \P_G=\emptyset$ such that
$G|_{G^{-1}(\LLL)}:G^{-1}(\LLL)\to \LLL$, as a \pf\ \rs, is c-equivalent to a holomorphic model, then $G$ is
\ce\ to a rational map.\ENDPROP

\beginp Set $\E=G^{-1}(\LLL)$ and $F=G|_{G^{-1}(\LLL)}$. By assumption $F$  is c-equivalent to a holomorphic model (with marked set $\P_F$, which is
equal to $\P_G\cap K_G$). By Lemma \ref{5'} there
is a pair  $(\Theta,\mu)$, with $\Theta$ a quasi-conformal map of
$\LLL$ satisfying $\Theta|_{\partial\LLL\cup(\P_G\cap\LLL)}=id$ and
$\Theta$ isotopic to the identity rel
$\partial\LLL\cup(\P_G\cap\LLL)$, with $\mu$ a Beltrami differential
on $\LLL$ such that $\|\mu\|_\infty<1$ and $(G\circ
\Theta^{-1})^*\mu=\mu|_{\Theta(\E)}$.

Choose $U_0$ an open neighborhood of $\P'_G$ disjoint from $\LLL$ so
that $G^{-1}(U_0)\supset\overline{U_0}$ and $G$ is holomorphic on
$G^{-1}(U_0)$. Set $L_0=U_0^c$ and $L_n=G^{-n}(L_0)$. As $L_n$
forms a decreasing sequence of sets shrinking down to $K_G$, there
is an integer $N\ge 0$ such that $L_N\subset \LLL$. So every orbit
passing through $L_0\smm \E$ stays there for at most $N+1$ times
before being trapped by $U_0$.

Extend the map $\Theta$ to a quasi-conformal map of $\cbar$ by
setting $\Theta:=id$ on $\cbar\smm\LLL$, then $\Theta$ is
quasi-conformal and isotopic to the identity rel $\P_G$. Set
$G_1=G\circ\Theta^{-1}$. Then $G_1$ is again quasi-regular, and is
holomorphic on $\Theta(U_0)=U_0$. Clearly, each $G_1$-orbit passes
through $L_0\smm\Theta(\E)$ at most $N+1$ times.

Extend now $\mu$ outside $\LLL$ by $\mu=0$. Let $\Phi_1:\cbar\to
\cbar$ be a global integrating map of this extended $\mu$. Set
$G_2:=\Phi_1\circ G_1\circ \Phi_1^{-1}$. Then $G_2$ is again
quasi-regular, and is holomorphic in the interior of
$\Phi_1\circ\Theta(\E)$ and in $\Phi_1(U_0)$. Elsewhere each $G_2$-orbit
passes at most $N+1$ times.

One can now apply the Shishikura principle: we spread out the
Beltrami differential $\nu_0\equiv 0$ using iterations of $G_2$ to
get an $G_2$-invariant Beltrami differential $\nu$. Note that
$\nu=0$ on $\Phi_1(U_0)$, and $\|\nu\|_\infty<1$.  Integrating $\nu$
by a \qc\ homeomorphism $\Phi_2$ (necessarily holomorphic on
$\Phi_1(U_0)$), we get a new map $R:= \Phi_2\circ G_2\circ  \Phi_2^{-1}$
which is a rational map and is \ce\ to $G_2$, therefore to $G$. See
the following diagram.

$\ds \begin{array}{cclclclc} (\cbar,\E) &
\stackrel{\Theta}{\longrightarrow}&\cbar&
\stackrel{\Phi_1}{\longrightarrow} & \cbar &
\stackrel{\Phi_2}{\longrightarrow} & \cbar &
\\
G\downarrow && \downarrow G_1 &&\downarrow G_2&&\downarrow R\\
(\cbar, \LLL) & \underset{id}{\longrightarrow} & \cbar &
\underset{\Phi_1}{\longrightarrow} & \cbar &
\underset{\Phi_2}{\longrightarrow} & \cbar &
\end{array} $
\qed

\section{Thurston-like theory for repelling systems}\label{repel}

We are thus interested in whether a given \rs\ (for example a
restriction of $G$ near $K_G$ as in Lemma \ref{restriction}) is
\ce\ to a holomorphic model. We will see that, similar to
Thurston's theory, the answer is yes if the map has no obstructions
that are similar to Thurston's original obstructions.

\subsection{Gr\"otzsch inequality and Thurston obstructions}\label{obstructions}

Thurston obstructions are in fact closely related to the
Gr\"otzsch inequality on moduli of annuli. The best way to
understand them is to start from real models.

{\noindent\bf 1. Slope obstruction.} Suppose we want to make a
tent map $f$ on $[0,1]$ with folding point $c$ and with $f(c)>1$,
with left slope $d_1$ and right slope $-d_2$. This is possible iff
$d_1^{-1}+d_2^{-1}<1$. More generally, suppose we have $k$
disjoint closed intervals $I_1\sqcup\cdots\sqcup I_k$ in $\R$, on
which we have a topological dynamical system with the following
combinatorics:

For each pair $(i,j)$, there are finitely many (might be zero)
intervals $I_{ij\delta}$, for $\de$ running through some finite or
empty index set $\La_{ij}$ (depending on $(i,j)$), such that \begin{enumerate}\item[
(TOP)] $I_{ij\delta}$ is a sub-interval of $I_i$, and the
$I_{ij\delta}$'s are mutually disjoint for all possible $i,j$ and
$\delta$. \item[
(DYN)] $f:I_{ij\delta}\to I_j$ is a homeomorphism  for all possible
$i,j$ and $\delta$.\end{enumerate}

The question we ask is: given a collection of slopes (in absolute
value) $(d_{ij\de})_{ij\de}$, is there a collection of
$(I_j,I_{ij\de}$ and $f:I_{ij\delta}\to I_j)$ so that each
$f:I_{ij\delta}\to I_j$ is affine with slope (in absolute value)
$d_{ij\de}$, for every possible multi-index $(i,j,\de)$?

Let us search at first necessary conditions. Assume that such a
piece-wise affine map $f$ exists. Then (DYN) implies
$|I_{ij\de}|=\frac{|I_j|}{d_{ij\de}}$ whereas  (TOP) implies
$\sum_{j,\de}|I_{ij\de}|< |I_i|$. Put these together we get
\REFEQN{inequality}
\sum_j\left(\sum_{\delta\in \La_{ij}}\frac1{d_{ij\delta}}\right) |I_j| < |I_i|\
,\  i=1,\cdots,k \ .\ENDEQN Let $D=(a_{ij})$ denotes the  {\em
transition matrix} defined by $a_{ij}:=\sum_{\delta\in\La_{ij}}  1/d_{ij\delta}$ (it is a
non-negative matrix, with, by convention, $a_{ij}=0$ if $\La_{ij}=\emptyset$). Then the necessary condition
\Ref{inequality} can be reformulated as: $Dv<v$, where
$v:=(|I_i|)$ is a vector with strictly positive entries.

It is quite easy to check that this necessary condition is also
sufficient. Therefore the answer to the above question is: such an
affine dynamical system exists if and only if the transition
matrix $D$ admits a vector $v$ with strictly positive entries so
that $Dv<v$, or equivalently (see Lemma \ref{matrix}), if and only
if  the leading eigenvalue $\la(D)$ of $D$ from Perron-Frobenius
theorem is strictly less than $1$.

Once this is done, the following 'complexification' will become
easy:

{\noindent\bf 2. Gr\"otzsch obstruction for annuli coverings.} Now we may think the
intervals $I_j$ are thin tubes, and the subintervals $I_{ij\de}$
as essential sub-annuli. More precisely,
let $\A=A_1\sqcup\cdots \sqcup A_k$ be a  \ps\  with each $A_i$ a closed annulus. For each pair $(i,j)$, let
$(A_{ij\delta})_{\delta}$ be finitely many (might be zero)
 sub-annuli of $A_i$ such that
\begin{enumerate}
\item[(TOP)]  $A_{ij\delta}$ is an
essential sub-annulus of $A_i$. Furthermore for each $i$ the $A_{ij\delta}$'s
are mutually disjoint for all possible choices of $j$ and $\delta$.
\item
[(DYN)] $f:A_{ij\delta}\to A_j$ is a quasi-regular covering of
degree $d_{ij\delta}$  for all possible $i,j$ and $\delta$.
\end{enumerate}

Set $\E=\bigsqcup_{ij\de} A_{ij\de}$ and consider $f:\E\to \A$ as an
 (un)repelling system with empty postcritical set (an annuli covering). The
question is: is $f:\E\to \A$ \ce\ to a holomorphic model?

Assume that $f$ is already holomorphic. Denote by $|A_*|$ the
modulus (rather than the length) of the interior of the annulus
$A_*$. Then $|A_{ij\de}|=|A_j|/d_{ij\de}$ (due to (DYN)) and
$\sum_{j,\de}|A_{ij\de}|<|A_i|$ (due to Gr\"otzsch inequality and
(TOP)). Therefore the leading eigenvalue $\la(D)$ of the
transition matrix $D$ is less than $1$. We have, naturally:

\REFLEM{annulicase}  An annuli covering $f:\E\to \A$ is \ce\ to a
holomorphic model if and only if $\la(D)<1$. \ENDLEM
\beginp $\Longrightarrow$. Assume $f:\E\to \A$  is \ce\ to a holomorphic $R:\E'\to\A'$.
Then the two maps have the same transition matrix $D$. By the argument above $\la(D)<1$.

$\Longleftarrow$. This will be done in Lemma \ref{Annuli}.\qed

This Lemma is not needed in the proof of our main result. But it
helps the understanding of Thurston obstructions and its proof
will shed light to our more complicated situation.

{\noindent\bf 3. Thurston obstruction} for a pair $(h,\P)$.  Let
$h:\cbar\to\cbar$ be a  branched covering, and $\P\subset\cbar$ a
closed marked  set containing $\P_h$ and $h(\P)$.
For example we may take $h$ to be a sub-hyperbolic semi-rational map $G$ and $\P=\P_G$.

A Jordan curve $\g$ in $\cbar\smm \P$ is said {\em
null-homotopic} (resp. {\em peripheral}) within $\cbar\smm \P$ if one of its
complementary component contains zero (resp. one) point of $\P$; and is said {\em non-peripheral} within $\cbar\smm \P$
otherwise, i.e. if each of its two
complementary components contains at least two points of $\P$.

We say that $\Gamma=\{\gamma_1, \cdots, \gamma_k\}$ is a {\em
multicurve of  within $\cbar\smm \P$}, if each $\g_i$ is a Jordan curve
in $\cbar\smm\P$ and is non-peripheral within $\cbar\smm \P$, and the
$\g_j$'s are mutually disjoint and mutually non-homotopic within $\cbar\smm \P$.

Each multicurve $\G$ induces a {\em $(h,\P)$-transition matrix}
$D_\G$ together with its {\em leading eigenvalue} $\la(D_\G)$ as
follows: Let $(\gamma_{ij\delta})_{\de\in \La_{ij}}$ denote the collection of the components of
$h^{-1}(\gamma_j)$ homotopic to $\gamma_i$  within $\cbar\smm \P$, with $\La_{ij}$
some finite or empty index set depending on $ij$.
 Then $h:\gamma_{ij\delta}\to \gamma_j$
is a topological covering of a certain degree $d_{ij\delta}$. The
transition matrix $D_\G=(a_{ij})$ is defined by $a_{ij}=\sum_{\de\in \La_{ij}} 1/d_{ij\de}$
(and $a_{ij}=0$ if $\La_{ij}=\emptyset$).

We say that a multicurve $\G$ is $(h,\P)$-{\em stable} if every
curve of $h^{-1}(\gamma)$, with $\gamma\in\G$,  is either null-homotopic
or peripheral within $\cbar\smm \P$,  or is homotopic within $\cbar\smm \P$ to a curve in $\G$.
This implies that for any $m>0$, every curve of $h^{-m}(\gamma)$,
$\gamma\in\G$ is either null-homotopic or peripheral within $\cbar\smm \P$,  or is
homotopic within $\cbar\smm \P$ to a curve in $\G$.

We say that a multicurve $\Gamma$ is {\em a Thurston obstruction
for $(h,\P)$} if it is $(h,\P)$-stable and $\la(D_\G)\ge 1$. In
the particular case $\P=\P_h$  we say simply that $\G$ is a
Thurston obstruction for $h$.

In case that $\P$ is finite (in particular $h$ is \pf)
we say that two such pairs $(h,\P)$, $(\widetilde h,\widetilde \P)$ are c-equivalent
if there is a pair of homeomorphisms
$(\phi,\psi):\cbar\to\cbar$ such that $\phi(\P)=\widetilde \P$, that $\phi$ is isotopic to $\psi$
rel $\P$, and that $\phi\circ h\circ \psi^{-1}=\widetilde h$.

\REFTHM{Thurston} (Marked Thurston theorem). Let $h$ be a \pf\
branched covering of $\cbar$ with $\deg h\ge 2$. Assume that the
signature of its orbifold is not $(2,2,2,2)$, or more particularly $\P_h$ contains periodic
critical points or at least five points. Let $\P$ be a finite marked
set containing $\P_h$ and $ h(\P)$. If $(h,\P)$
has no Thurston obstructions, then $(h,\P)$ is c-equivalent to a unique
rational map model. More precisely there are homeomorphisms
$(\phi,\psi):\cbar\to\cbar$ such that $\phi$ is isotopic to $\psi$
rel $\P$ and that $f:=\phi\circ h\circ \psi^{-1}$ is a rational map.
And the conformal conjugacy class of the pair $(f,\phi(\P))$ is unique.

Furthermore, if $h$ is \qr,  both $\phi$ and $\psi$ can be taken
to be quasi-conformal.\ENDTHM

Here we omit the  definition of  orbifold and its signature (see
e.g. \cite{DH1}). We mention only that if $\P_h$ contains periodic
critical points, or at least 5 points, then the signature of its
orbifold is not $(2,2,2,2)$. This is enough for our purpose here.

\noindent {\bf Remark}. Our statement is slightly different than
the original Thurston Theorem, where $\P=\P_h$. But the arguments
in \cite{DH1} can be easily adapted to prove this more general
form. In case $h$ is \qr, we may replace $\phi$ by a
quasi-conformal map $\phi'$ isotopic to $\phi$ rel $\P$. This is
possible since $\P$ is finite (see Lemma \ref{open}). Lifting this
isotopy will provide us a quasi-conformal map $\psi'$ isotopic to
$\psi$ rel $\P$ such that $\phi\circ h\circ \psi^{-1}=\phi'\circ
h\circ \psi'^{-1}$.

Conversely, we have the following result of McMullen (\cite{Mc2}):

\REFTHM{mcm}  Let $f$ be a rational map with $\deg f\ge 2$, and
let $\P$ be a closed subset such that $f(\P)\subset \P$ and
$\P_f\subset \P$. Let $\G$ be a $(f,\P)$-multicurve whose
transition matrix is denoted by $D$. Then $\la(D)\le 1$. If
$\la(D)=1$, then either $f$ is \pf\ whose orbifold has signature
$(2,2,2,2)$; or $\G$ includes a curve that is contained in the
Siegel discs or Herman rings of $f$. \ENDTHM

Again this form is slightly stronger than McMullen's original
version. But the proof goes through without any trouble.

\subsection{Thurston obstructions for repelling systems.}

 Let $(\E\stackrel{F}{\to}\LLL,\P)$ be a marked \pf\ \rs,
 in other words $F:\E\to\LLL$ is a \qr\ branched covering among two nested \ps s, and $\P\subset \LLL\smm \partial \LLL$ is a finite set containing $\P_F$ and
 $F(\P\cap \E)$.  (In case $\LLL=\cbar$  we are back to Thurston's setting).

Two Jordan curves
in $\LLL\smm \P$ are {\em homotopic} if they are both contained in a common $\LLL$-piece $S$
and are homotopic to each other within $S\smm \P$.

A Jordan curve $\g\subset \LLL\smm \P$ is said {\em
null-homotopic} (resp. {\em peripheral}) within $\LLL\smm \P$ if it bounds an open
disc $D$ so that $D\subset \LLL$ and $D\cap \P=\emptyset$ (resp. $\#D\cap \P=1$);
 is said
 {\em non-peripheral} within $\LLL\smm \P$ otherwise (this is equivalent to say that $\g$ is contained
 in $S\smm \P$ for some component $S$ of $\LLL$, and either $\g$ bounds no disc in $S$, or
 $\g$ bounds a disc $D$ in $S$ containing at least
 two points of $\P$). For example if $\g$ is a boundary curve of an $\LLL$-piece $S$, and if
 $S$ is not a closed disc, then $\g$ is non-peripheral.

We say that $\Gamma=\{\gamma_1, \cdots, \gamma_k\}$ is a {\em
multicurve} within $\LLL\smm \P$, if each $\g_i$ is a Jordan curve
in $\LLL\smm \P$ and is non-peripheral within $\LLL\smm \P$, and the
$\g_j$'s are mutually disjoint and mutually non-homotopic within $\LLL\smm\P$.

Each multicurve $\G$ induces an {\em $(F,\P)$-transition matrix}
$W=W_\G$ together with its {\em leading eigenvalue} $\la(W_\G)$ as
follows: Let $(\gamma_{ij\delta})_{\de\in \La_{ij}}$ denote the collection of the components of
$F^{-1}(\gamma_j)$ homotopic to $\gamma_i$  within $\LLL\smm \P$, with $\La_{ij}$
some finite or empty index set depending on $ij$.
 Then $F:\gamma_{ij\delta}\to \gamma_j$
is a topological covering of a certain degree $d_{ij\delta}$. The
transition matrix is defined by $$W_\G=(b_{ij}),\quad b_{ij}=\sum_{\de\in \La_{ij}} 1/d_{ij\de}$$
(with $b_{ij}=0$ if $\La_{ij}=\emptyset$).

We say that a multicurve $\G$ within $\LLL\smm \P$ is $(F,\P)$-{\em stable} if every
curve of $F^{-1}(\gamma)$, with $\gamma\in\G$,  is either null-homotopic
or peripheral within $\LLL\smm \P$,  or is homotopic within $\LLL\smm \P$ to a curve in $\G$.

We say that a multicurve $\Gamma$ within $\LLL\smm \P$ is {\em a Thurston obstruction
for $(F,\P)$} if it is $(F,\P)$-stable and $\la(W_\G)\ge 1$. See \S\ref{Examples} for examples.

The following principle will be used  frequently, and is a direct consequence of the fact that $F(\P\cap \E)\subset \P$ and
that $F$ is a covering over $\LLL\smm \P$:

 \medskip

{\noindent\bf Basic pull-back principle}. \begin{enumerate}\item
 Let $D$ be an open Jordan disc contained in $\LLL\smm \P$ with $\partial D\cap \P=\emptyset$. Then every
component of $F^{-1}(D)$ is again an open disc and is contained in $\E\smm \P$.
Each curve in $F^{-1}(\partial D)$ is the boundary of a component of $F^{-1}(D)$.\item
 Let $A$ be an open annulus contained in $\LLL\smm \P$. Then every
component of $F^{-1}(A)$ is again an open annulus and is contained in $\E\smm \P$.\item
 Let $D$ be an open Jordan disc contained in $\LLL$ with $\partial D\cap \P=\emptyset$ such that $D$ contains a unique
point of $\P$. Then
every component of $F^{-1}(D)$ is again an open disc, is contained in $\E$, and contains at most one point of $\P$. Each curve in $F^{-1}(\partial D)$ is the boundary of a component of $F^{-1}(D)$.\end{enumerate}

The following is an easy consequence:

\REFLEM{peripheral} Let $(\E\stackrel{F}{\to}\LLL,\P)$ be a marked \pf\ \rs.  For any peripheral (resp. null-homotopic) curve
$\g\subset\LLL\smm\P$, each curve in $F^{-1}(\g)$ is either peripheral
or null-homotopic (resp. is null-homotopic).\ENDLEM

 We will prove:

\REFTHM{more-general} (Thurston theorem for marked repelling systems). Let $(\BBB\stackrel{G}{\to}\MMM,\Q)$ be a marked \pf\ \rs. We assume in addition that
no $\MMM$-piece is homeomorphic to $\cbar$, in other words we require $\partial S\ne \emptyset$
for every $\MMM$-piece $S$. If $(G,\Q)$
has no Thurston obstructions, then   $(G,\Q)$ is c-equivalent
to a holomorphic model map.
\ENDTHM

\subsection{Proof of Theorem \ref{hyperbolic} using Theorem \ref{more-general}}\label{nono}

\REFLEM{un-obstructed} Assume that $G:\cbar\to \cbar$ is a \qr\ sub-hyperbolic
\sr\ with $\P'_G\ne \emptyset$  and without Thurston obstructions.
Then there are puzzle surfaces $L_1,L_0$ such that
$$K_G\subset\subset L_1\subset \subset L_0,\ G^{-1}(L_0)=L_1$$ and that, the restriction $G|_{L_1}: L_1\to L_0 $,
marked by $\P_G\cap L_0$, is  a  \pf\ \rs\ without Thurston obstructions.\ENDLEM
\beginp One can find an open set $U_0$ which is
the union of a quasi-disc neighborhood for each point of $\P_G'$
so that these quasi-discs have disjoint closures, that $\partial
U_0$ is disjoint from $\P_G$, that $G$ is holomorphic in a
neighborhood of $\overline{U_0}$, and that $G(U_0)\subset\subset
U_0$. Set $L_0=U_0^c$. Topologically $L_0$ is the sphere minus
finitely many (open) holes. Set $L_1=G^{-1}(L_0)$, $H=G|_{L_1}$ and  $\Q:=\P_G\cap L_0$.
So $(H,\Q)$ is a marked \pf\ \rs.
The assumption $\P'_G\ne\emptyset$ implies that $\partial L_0\ne \emptyset$.

We will show now: under the assumption that $G:\cbar\to \cbar$ has no Thurston obstructions,
the  \rs\ $(H,\Q)$ has no Thurston obstructions.

 Assume at first that $L_0$ is a closed disc containing at most one point of $\Q$.
 In this case $\partial L_0$ is a single curve and is null-homotopic or peripheral within
 $L_0\smm \Q$. And there is no multicurve within $L_0\smm \P$ as every curve in $L_0\smm \P$
 is either null-homotopic or peripheral. Consequently $(H,\Q)$ has no Thurston obstructions.

 Next assume that $L_0$ is a closed annulus disjoint from $\Q$. Then there is only one class of non-peripheral Jordan curves within $L_0\smm \Q=L_0$, namely that of
 a  boundary curve $\g$ of $L_0$. But such a $\g$ is also non-peripheral within $\cbar\smm \P_G$ as each of the two discs of $\cbar\smm L_0$ contains points of $\P'_G\subset \P_G$.
 The curves in $G^{-1}(\g)$ are either peripheral within $L_0$, or homotopic to
 $\g$ within $L_0$. Therefore $\{\g\}$ is stable for both $(\cbar\stackrel{G}{\to}\cbar,\P_G)$ and $(H,\Q)$
 and the corresponding two transition matrices are identical. By the assumption that $(\cbar\stackrel{G}{\to}\cbar,\P_G)$ has no Thurston obstructions, the corresponding
 leading eigenvalue is less than $1$. Therefore $\{\g\}$ is not an Thurston obstruction
 for $(H,\Q)$. And $(H,\Q)$
 has no obstructions.

 In the remaining case, $L_0$ is a \rc, with $$\#(L_0\cap \Q)+\#\{\text{boundary curves of }L_0\}\ge 3.$$ In particular each of its boundary curves is non-peripheral
 within $L_0\smm \Q$.

 Let now $\G$ be a multicurve within $L_0\smm \Q$. In other words, \begin{enumerate}\item[
a)] each curve in $\G$ is non-peripheral within $L_0\smm\Q$,\item[
b)] the curves in $\G$ are mutually disjoint,\item[
c)] the curves in $\G$ are mutually non-homotopic within $L_0\smm\Q$.\end{enumerate}

We want to show that $\G$ is also a multicurve within $\cbar\smm \P_G$, i.e. $\G$ satisfies a),b) and c) with $L_0\smm\Q$ replaced by $\cbar\smm \P_G$.
By a), for each curve $\g$ in $\G$, either both discs of $\cbar\smm\g$ contains a component of  $\cbar\smm L_0=U_0$ (therefore infinitely many points of $\P_G$); or one disc of $\cbar\smm\g$
is contained in $L_0$ and contains at least two points of $\Q\subset \P_G$, while the other disc contains all components of $U_0$ (therefore infinitely many points of $\P_G$).
In both cases each component of $\cbar\smm \g$ contains at least two points of $\P_G$. So $\g$ is non-peripheral within $\cbar\smm \P_G$.

By b) the curves in $\G$ are mutually disjoint.

By c), given any two curves $\g,\g'$ of $\G$, we have $\g,\g'\subset L_0\smm \Q$ and the open annulus $A(\g,\g')$ bounded by
$\g,\g'$ intersects either $U_0$ or $\Q\subset \P_G$ (or both).
In the former case $A(\g,\g')$ contains
a component of $U_0$. Therefore in both
cases $A(\g,\g')$ intersects $\P_G$. So $\g,\g'$ are also non-homotopic within $\cbar\smm \P_G$.

This arguments implies that $\G$ is also a multicurve within $\cbar\smm \P_G$.

Assume now that $\G$ is a  multicurve within $L_0\smm \Q$ and is $(H,\Q)$-stable. In other words, for any $\g\in \G$ and
any curve $\de$ of $G^{-1}(\g)$, either $\de$ bounds a disc that is contained in $L_0$
and that contains at most one point of $\Q=L_0\cap \P_G$, or $\de$ is homotopic
within $L_0\smm \Q$ to a curve $\g'$ in $\G$. Thus either $\de$ bounds a disc
that contains at most one point of $\P_G$, or it is homotopic
within $\cbar\smm  \P_G$ to  $\g'$.

This shows that $\G$ is a multicurve within $\cbar\smm \P_G$ that is also $(\cbar\stackrel{G}{\to} \cbar,\P_G)$-stable. The two transition matrices
(by $(H,\Q)$ and by $(\cbar\stackrel{G}{\to} \cbar,\P_G)$)
are identical, therefore have the same leading eigenvalue $\la$.

By the assumption that $(\cbar\stackrel{G}{\to} \cbar,\P_G)$ has no Thurston obstructions,
we know that $\la<1$. So $\G$ is not a Thurston obstruction for $(H,\Q)$.

Therefore $(H,\Q)$
 has no obstructions.
\qed

Assuming Theorem \ref{more-general}, we may now give the

{\noindent\em Proof of Theorem \ref{hyperbolic}} (the existence
part). Let $G$ be sub-hyperbolic \sr\ with $\P_G'\ne \emptyset$
and without Thurston obstruction. We may assume in addition that
$G$ is globally \qr, up to a change of representatives in its
c-equivalence class (by Lemma \ref{BB}).

We may then apply Lemma \ref{un-obstructed} to $G$ to show that
it has a restriction  near $K_G$ which is a \pf\ \rs\  without
Thurston obstructions, therefore is c-equivalent to a holomorphic model by Theorem
\ref{more-general}. We may then apply Proposition \ref{final} to conclude that
$G$ is c-equivalent to a rational map.\qed

\section{Reduction to a sub-repelling system}\label{againagain}

Our main objective here is:

 \REFTHM{Unobstructed} Let   $(\BBB\stackrel{G}{\to}\MMM,\Q)$ be  a marked \pf\ \rs\ such that every $\MMM$-piece has a non-empty boundary.

Let $\E,\LLL$ be two \ps\ neighborhoods of $K_G$ satisfying:
\begin{enumerate}\item[
{\rm (*)}] $K_G\subset\subset \E\subset\subset\LLL\subset\subset \MMM$, $\E=G^{-1}(\LLL)$, $\Q\cap \partial
\LLL=\emptyset$;
\item[
{\rm
(**)}] for every $\LLL$-piece $S$,
 and for the $\MMM$-piece $S_0$ containing $S$ in the interior, and for the copy $\cbar$ of the Riemann sphere containing $S_0$ (therefore $S$),  every (disc-like) component of $\cbar\smm S$ contains either components of $\partial S_0$ or points of $ \Q$ (or both);
\end{enumerate}
 Let $F=G|_{\E}:\E\to\LLL$ be the sub-\rs\   marked by $\Q\cap \LLL$.
  Then,
 \begin{enumerate}
\item[{\rm  (A)}] If $(F,\Q\cap \LLL)$ is c-equivalent to a holomorphic model so is $(G,\Q)$.
 \item[{\rm (B)}]
If $(G,\Q)$ has no Thurston obstructions, so does $(F,\Q\cap \LLL)$.
$$\begin{array}{ccc}
G \text{\ has no ob.} &  & G\sim_c \text{ hol. model}\vspace{0.2cm}\\
\Downarrow {\rm (B)} && \Uparrow {\rm (A)}\vspace{0.2cm}\\
F \text{\ has no ob.} & \underset{?}{\Longrightarrow} & F\sim_c \text{ hol. model.}
  \end{array}$$\end{enumerate}\ENDTHM

Once the theorem is proved the problem of c-equivalence to a holomorphic model for $(G,\Q)$ is reduced
the problem for a suitable sub-\rs.

\subsection{How to get a stable multicurve}\label{News}

The following criterion is very useful:

\REFLEM{non-stable} A marked \pf\ \rs\ $(\E\stackrel{F}{\to}\LLL,\P)$ has a Thurston obstruction
if and only if there is a multicurve $\G'$ within $\LLL\smm \P$ (not necessarily $(F,\P)$-stable) with $\la(W_{\G'})\ge 1$.
\ENDLEM
\beginp We will need to produce an $(F,\P)$-stable multicurve with the leading eigenvalue
of its transition matrix greater than one.

Let $\G$ be a multicurve within $\LLL\smm \P$, i.e. the curves in $\G$ are non-peripheral, mutually disjoint and mutually non-homotopic within $\LLL\smm \P$. The action of $(F,\P)$ induces  a directed graph $\La_\G$ as follows: the vertices are
 the curves in $\G$. And there is an edge directing from $\de$ to $\g$ (maybe $\de=\g$) if
 $\de$ is homotopic to a curve in $F^{-1}(\g)$ within $\LLL\smm \P$.

 If $A\subset \G$ is a sub-multicurve then $W_{A}$ is the submatrix of $W_\G$ corresponding to the entries of $A$,
 denoted by $W_\G|_{A}$, and $\La_{A}$ is the corresponding subgraph of $\La_\G$.
 In this case $\la(W_{A})\le \la(W_\G)$ (by Lemma \ref{compare}).

 We will say that a multicurve $A$ within $\LLL\smm \P$ is {\em irreducible} if any two vertices of its graph $\La_{A}$ can be
 connected by following successively the directed edges. It is elementary that if $\la(W_\G)>0$ then
 there is an irreducible sub-multicurve $A\subset \G$ with $\la(W_{A})= \la(W_\G)$.

Each multicurve $\G$ defines a pulled-back multicurve $\G_1$ as follows:
The curves in $F^{-1}(\G):=\bigcup_{\g\in \G}F^{-1}(\g)$ are already mutually disjoint.
But some of them might be peripheral, or homotopic to another within $\LLL\smm \P$.
Pick one representative in each homotopic class of the non-peripheral curves in
$F^{-1}(\G)$. Together they form a new multicurve, $\G_1$.

In general a curve in $\G_1$ might not be homotopically disjoint from a curve in $\G$.

Saying that $\G$ is stable is equivalent to say that every curve in $\G_1$ is
homotopic to a curve in $\G$.

In case that every vertex of $\La_\G$ is the departure of an edge (for example when
$\G$ is irreducible), then every curve in $\G$ is
homotopic to a curve in $\G_1$. One obtains then
successive pulled back multicurves $\G_2,\cdots$, such that curves in $\G_i $ are homotopic
to curves in  $ \G_{i+1}$. As $\LLL\smm \P$ is
topologically finite, there is $N$ such that $\#\G_N=\#\G_{N+1}$. Consequently $\G_N$ is $(F,\P)$-stable.

Assume now $\G'$ is a multicurve within $\LLL\smm\P$, not necessarily stable, such that $\la(W_{\G'})\ge 1$.

There there is an irreducible sub-multicurve $A\subset \G'$ with $\la(W_A)=\la(W_{\G'})$.
Pulling back successively $A$ we get a multicurve $A_N$ that is $(F,\P)$-stable
and contains $A$. And
$1\le \la(W_{\G'})=\la(W_{A})\le \la(W_{A_N}).$ \qed

\subsection{Proof of Theorem \ref{Unobstructed}}
{\em Proof of Theorem \ref{Unobstructed}}.\\
(A). This part can be proved similarly as Proposition \ref{final}. We will omit the details here.
\\
(B). Note that
Condition (**) assures a certain minimality of $\LLL$, so that $\LLL$ does not have trivial holes in $\MMM$.

Set $\P=\Q\cap \LLL$, the marked set of $F$.

At first we prove the  following facts that will be used frequently in the sequel:
\begin{enumerate}
\item[(a)] Two Jordan curves $\g,\g'$
in $\LLL\smm\P$, homotopic within $\LLL\smm \P$, are also  homotopic
within $\MMM\smm \Q$.
\item[(b)]  For a Jordan curve $\g$
in $\LLL\smm\P$, it is null-homotopic within $\LLL\smm \P$ iff it is null-homotopic
within $\MMM\smm\Q$.
\item[(c)] A peripheral curve within $\LLL\smm \P$  is also peripheral within $\MMM\smm\Q$.
\end{enumerate}

{\noindent\em Proof.} (a). There is an $\LLL$-piece $S$ containing both $\g,\g'$
and $\g$ and $\g'$ are homotopic within $S\smm \P$. But $S$ is contained in an $\MMM$-piece $S_0$, and $\P\cap S=\Q\cap \LLL\cap S=\Q\cap S$. So $S\smm \P\subset S_0\smm \Q$ and $\g$ and $\g'$ are homotopic within $S_0\smm \Q$, therefore within $\MMM\smm\Q$.

(b). Again let $S$ (resp. $S_0$) be the $\LLL$-piece (resp. $\MMM$-piece) containing $\g$.

If $\g$ is null-homotopic within $\LLL\smm \P$ then it bounds a disc $D$ contained in
$S\smm \P$. But $S\smm \P\subset S_0\smm \Q$. So $D\subset S_0\smm \Q$
and $\g$ is null-homotopic within $S_0\smm \Q$, therefore within $\MMM\smm \Q$.

Conversely if $\g$ is null-homotopic within $\MMM\smm \Q$, it bounds an open  disc $D$ contained in $S_0\smm \Q$. If $D$ is not contained in $S$, and as $\g=\partial D\subset S$, we must have that $D$ contains a component of $\cbar\smm S$. By Condition (**)   we know that $D$ intersects $\partial S_0\cup \Q$.
This contradicts the fact that $D\subset S_0\smm \Q$. Therefore $D\subset S$.
But $D\cap \Q=\emptyset$. So $D\subset S\smm \Q=S\smm \P$.
Therefore $\g$ is also null-homotopic within $S\smm \P$, hence within $\LLL\smm\P$.

(c). By definition $\g$ is peripheral within $\LLL\smm \P$ if it bounds a disc $D$ that is contained
in an $\LLL$-piece $S$ such that $\#D\cap \P=1$. In this case $D$ is also contained in the
$\MMM$-piece $S_0$ that contains $S$, and $\#D\cap \Q=\#D\cap (Q\cap S)=\#D\cap \P=1$.
By definition again $\g$ is also peripheral within $\MMM\smm \Q$.
\qed

Assume now that $(G,\Q)$ has no Thurston obstructions. We will prove that $(F,\P)$ has
no obstructions either.

Let $T$ be a multicurve within $\LLL\smm \P$, i.e.:\\
i) each curve in $T$ is non-peripheral within $\LLL\smm \P$;\\
ii) the curves in $T$ are mutually disjoint;\\
iii) the curves in $T$ are mutually non-homotopic within $\LLL\smm \P$.

Change the representatives within the same homotopy classes if necessary, we may assume
in addition:\\
iv) a curve in $T$ is either equal to a boundary curve of $\LLL$, or, is disjoint from $\partial \LLL$ and is not homotopic to a curve in $\partial \LLL$ within $\LLL\smm \P$.

 When considering
homotopy within $\MMM\smm \Q$ (which contains $\LLL\smm \P$), there are two new phenomena:

1. Some of the curves in $T$ may now become peripheral (but never null-homotopic) within
$\MMM\smm \Q$ (Figure \ref{Peripheral} shows how this may happen).

2. Some of the curves in $T$ may now become homotopic to each other within $\MMM\smm \Q$.

The following two parts treat each phenomenon separately:

{\noindent\bf 1. Curves in $T$ that become peripheral within
$\MMM\smm \Q$}.

\begin{figure}[htbp]
\begin{center}
\input{peripheral.pstex_t}
\end{center}
\vspace{-0.5cm} \caption{\it The curve $\g$ is not peripheral within $\LLL\smm \P$ but is
peripheral within $\MMM\smm \Q$, as it bounds a disc $D(\g)$ that is contained in $\MMM$
and contains a unique point of $\Q$.}
\label{Peripheral}
\end{figure}
We now consider homotopy within $\MMM\smm \Q$. By (b) each curve in
$T$ is non-null-homotopic within $\MMM\smm \Q$, is therefore either peripheral or
non-peripheral within $\MMM\smm \Q$. We thus decompose $T$ into
$T=Z\sqcup X$, with $Z$ (resp. $X$) denoting the collection of curves in $T$ that
are peripheral (resp. non-peripheral) within $\MMM\smm \Q$. Denote by $W_Z, W_X$ the
$(F,\P)$-transition matrix of $Z$, $X$ respectively.

{\noindent\bf Lemma I.} The $(F,\P)$-transition matrix $W_T$ has the
following block decomposition (where $O$ denotes a rectangle zero-matrix):
$\ds W_T= \left(\begin{array}{cc} W_X & O  \\ * & W_Z \end{array}\right).$

 Proof. Let $\g\in Z$ and $\beta$ be a curve of $F^{-1}(\g)=G^{-1}(\g)$. We just need to show that if $\beta$ is homotopic within $\LLL\smm\P$ to a curve $\g'$ in $ T$, then $\g'\in Z$.

 By definition of $ Z$ the curve $\g$ bounds an open disc $D(\g)$ contained in $\MMM$
and containing a unique point of $\Q$. Therefore, applying the Basic pull-back
principle to $D(\g)$, we know that each component of $G^{-1}(D(\g))$ is disc-like,
is contained in $\BBB=G^{-1}(\MMM)$, and contains at most one point of $\Q$. Let $\beta$
be a curve in $G^{-1}(\g)$. Then $\beta$ is the boundary of a component of  $G^{-1}(D(\g))$. So $\beta$ is either
null-homotopic or peripheral (within $\MMM\smm \Q$). On the other hand, $\beta$ is homotopic within $\LLL\smm\P$ to $\g'\in T$ by assumption. And the homotopy can
be taken within $\MMM\smm \Q$ by (a). Consequently $\g'$ is either
null-homotopic or peripheral within $\MMM\smm \Q$.

But no curves in $T$ are null-homotopic  within $\MMM\smm \Q$ (by (b) and by the
definition of $T$), so both $\beta$ and $\g'$ are peripheral within $\MMM\smm \Q$.
Therefore $\g'\in Z$. \qed

{\noindent\bf Lemma II}. Each curve of $Z$ is a boundary curve of $\LLL$.

Proof. Let $\g\in T$ that is not a boundary curve of $\LLL$. We just need to show that $\g$
is necessarily non-peripheral within $\MMM\smm\Q$.

There is an $\LLL$-piece $S$, an $\MMM$-piece $S_0$ and a Riemann sphere $\cbar$ such that
$\g$ is contained in the interior of $S$ and $ S\subset\subset S_0\subsetneq \cbar$.
Now $\cbar\smm \g$ has two disc components $D_1,D_2$.

By i) either one disc, say $D_1$, is contained in $S$ and
contains at least two points of $\P\subset \Q$, or both $D_1,D_2$ intersect $\cbar\smm S$.

In the former case, $D_1\subset S\subset S_0$, so $D_2$ contains $\partial S_0$, which by assumption
is non-empty. This implies that $\g$ is non-peripheral within $S_0\smm \Q$.

In the latter case, as $\g\subset S$, each $D_i$  contains a component
of $\cbar\smm S$. By Condition (**), each $D_i$ contains either a curve in $\partial S_0$ or
points in $ \Q$ (or both). Assume by contradiction that $\g$ is peripheral within $S_0\smm \Q$.
Then one of $D_1,D_2$, denoted by $D$, contains a unique point of $\Q$ and no boundary component of $S_0$. We have $D\subset S_0$. But $\g$ is not peripheral within $S\smm \P$
by i). So $D$ is not contained in $S$. Therefore $D$ contains components of $\cbar\smm S$.
 Let $\De$ be one component of $\cbar\smm S$ contained in $D$. Then $\De$ is
bounded by a curve $\de$ which is also a boundary curve of $S$. By Condition (**) again $\De$
must intersect $\partial S_0\cup \Q$. But $\De\subset D\subset S_0$.
So $\De\cap \partial S_0=\emptyset$. On the other hand, $D\cap \Q$ consists of a single point, denoted by $a$. So $a\in \De$. This shows that $D$ contains a unique component of $\cbar\smm S$. By iv) $\g\cap \partial S=\emptyset$, so $D$ contains
a unique boundary curve of $S$, which must be $\de$. Furthermore the annulus $A(\g,\de)$ does not contain $a$. Therefore $A(\g,\de)\subset S$ and $$A(\g,\de) \subset D\smm\{a\}=D\smm \Q\subset D\smm \P.$$ So $A(\g,\de)\subset S\smm \P$.
This means that the boundary curve $\de$ of $S$ is homotopic to $\g$ within $S\smm\P$. This is not possible by iv).\qed

{\noindent\bf Lemma III.} There is a power $q\ge 1$ such that $(W_Z)^q=O$.
Therefore $\la(W_Z)=0$ and $\la(W_X)=\la(W_T)$.

Proof. Let $\GGG$ be the union of $Z$ with the curves in $\partial \LLL$ that are peripheral within $\LLL\smm \P$ (these added curves are disjoint from curves in $Z$ by i) and iv) ). By Claim II every curve in $\GGG$
is a boundary curve of $\LLL$.

By definition a curve $\g\in \partial \LLL$ is peripheral if it is the boundary of an $\LLL$-piece
$D(\g)$ (therefore $D(\g)$ is a disc) with $\#D(\g)\cap \P=1$. Note that $D(\g)\subset \LLL\subset \MMM$.

 Let $\g\in\GGG$. Then $\g$ bounds a disc $D(\g)$ which
is contained in $\MMM$ and which
contains a unique point $a(\g)$ of $\Q$.

We decompose $\GGG$ into $\GGG_p\sqcup \GGG_0$ according to $a(\g)$ is periodic or not.

By Basic pull-back principle, for all $k\ge 1$,  the
components of $G^{-k}(D(\g))$ are all disc-like, and each is bounded by exactly one curve of   of $G^{-k}(\g)$.

Assume $\g\in \GGG_0$, so that
 $a(\g)$ is not periodic (as $G$ is \pf, the orbit of $a(\g)$ is either preperiodic or eventually
 escapes $\BBB$). Then there is an integer $k(\g)\ge 1$ such that $G^{-k(\g)}(a(\g))$ contains no
points of $\Q$. But $G^{-k}(D(\g))\cap \Q\subset G^{-k(\g)}(a(\g))$. So $G^{-k(\g)}(D(\g))\cap \Q=\emptyset$. Thus the curves in  $G^{-k(\g)}(\g)$
are null-homotopic within $\MMM\smm \Q$ and hence are
null-homotopic within $\LLL\smm\P$ by (b). Therefore for all $k\ge k(\g)$,  the curves in  $G^{-k}(\g)$ are all
null-homotopic within $\LLL\smm\P$.

There is therefore a common integer $k$, so that for every $\g\in \GGG_0$, the curves in  $G^{-k}(\g)$ are all
null-homotopic within $\LLL\smm\P$.

Let now $\g\in \GGG_p$, i.e. with $a(\g)$ periodic. Set $a=a(\g)$. Denote by $p$ its period.
This implies in particular that the orbit of $a$ does not escape $\MMM$, so $a\in K_G=K_F\subset \LLL$.

Denote by $\{\eta_1,\cdots,\eta_m\}$  the curves in
$\GGG$ homotopic to $\g$ within $\MMM\smm \Q$, i.e. each $\eta_j$ bounds a disc $D(\eta_j)$ which is contained in $\MMM$ with $D(\eta_j)\cap \Q=\{a\}$. As $\eta_i\cap \eta_j=\emptyset$ for $i\ne j$,
we have either $D(\eta_i)\subset D(\eta_j)$ or $D(\eta_i)\supset D(\eta_j)$. We numerate the $\eta_j$'s
in the increasing order, i.e. such that $D(\eta_j)$ contains $D(\eta_{j-1})$ and $\eta_{j-1}$. The smallest disc $D(\eta_1)$ must be contained in $\LLL$, since there is no
curve in $\partial\LLL$ separating $a\in\LLL$ from $\eta_1\subset\partial\LLL$.
Therefore $\overline D(\eta_1)$ is an $\LLL$-piece, and $\eta_1$ is peripheral within $\LLL\smm \P$.

Fix $j\in \{1,\cdots,m\}$. The components of $G^{-p}(D(\eta_j))$  are all disc-like, with one of them, denoted by $D(\beta_j)$, containing $a$, and the others containing a preimage of $a$ that is not periodic.

Therefore $G^{-p}(\eta_j)=F^{-p}(\eta_j)$ has a unique
component $\beta_j$, which is the boundary of $ D(\beta_j)$, homotopic to $\g$ (and to $\{a\}$) within $\MMM\smm \Q$, and the other
components are either null-homotopic or homotopic within $\LLL\smm\P$ to a
curve in $\GGG_0$.

{\noindent \bf Claim}. $\overline D(\beta_1)\subset D(\eta_1)$ and $\overline{D}(\beta_j)
\subset D(\eta_{j-1})$ for $j=2,\cdots,m$.

Proof. At first the enlarged collection of curves $\{\eta_1,\cdots,\eta_m,\beta_1,\cdots, \beta_m\}$ are mutually disjoint. This is clearly true between the $\eta_j$'s and between
de $\beta_j$'s. But $\beta_j\cap \eta_i=\emptyset$ as well, as $\eta_i\subset \partial \LLL$,
$\beta_j\subset F^{-p}(\LLL)$ and $F^{-p}(\LLL)$ is contained in the interior of $\LLL$.

Therefore the discs $D(\eta_i)$, $D(\beta_j)$ are nested in a certain order.

We prove now that $\beta_j\subset D(\eta_j)$ for
$j=1,\cdots,m$.

We prove it at first for $\beta_1$. Note that $\beta_1$ is the
boundary of $E_1$, the component of $F^{-p}(\overline D(\eta_1))$ containing $a$.
But $F^{-p}(\LLL)$ is contained in the interior of $\LLL$, so
$E_1$ is contained in the interior of an $\LLL$-piece, which must be $\overline D(\eta_1)$,
i.e. $\beta_1\subset D(\eta_1)$.

Assume by contradiction that there is a minimal integer $j\ge 2$  such that
$\beta_j\not\subset D(\eta_j)$. Then $\beta_j\cap
\overline{D}_j=\emptyset$ (due to again, $F^{-p}(\LLL)\cap \partial \LLL=\emptyset$).
The annulus $A(\eta_j,\eta_{j-1})$ is contained in $D(\eta_j)$ therefore in $\MMM$.
So its inverse images by $G$ are well defined.
By the Basic pull-back principle the components of $G^{-p}(A(\eta_j,\eta_{j-1}))$
are all annuli. One of them must be
$A(\beta_j,\beta_{j-1})$, which contains
$A(\eta_j,\eta_{j-1})$ as a sub-annulus. Therefore
$$G^m(A(\eta_j,\eta_{j-1})),\ m=0,\cdots,p-1$$ are all contained in $\BBB=G^{-1}(\MMM)$.
Set $A'=\bigcup_{m=0}^{p-1} G^m(A(\eta_j,\eta_{j-1}))$. Then $G(A')\subset A'$. Trivially either
$\overline A(\eta_j,\eta_{j-1})$ is an $\LLL$-piece or there is a point $z\in
A(\eta_j,\eta_{j-1})\cap (\BBB\smm \LLL)$. The former case is not possible, as components of $F^{-p}(\LLL)$ are compactly contained in $\LLL$. In the latter case, the $G$-orbit
of $z$ never escapes $A'$ which is a subset of $\BBB$, so $z\in K_G$. This is again impossible as $K_G=K_F\subset \LLL$ and $z\notin \LLL$.

This proves that $\beta_j\subset D(\eta_j)$ for all
$j=1,\cdots,m$. It follows that $\overline D(\beta_j)\subset D(\eta_j)$.

Fix $j= 1,\cdots,m$. Denote
 by $S_j$ the $\LLL$-piece containing $\eta_j$ as a boundary curve. Thus $S_1=D(\eta_1)$.
We want to show now for $j\ge 2$ either $S_j\cap D(\eta_j)=\emptyset$ or $S_j= \overline A(\eta_j, \eta_{j-1})$.

Assume $S_j\cap
D(\eta_j)\ne \emptyset$. We have
$S_j\subset\overline{D}_j\smm\{a\}$. Then $\eta_{j-1}\subset \partial S_j$ since no
curve in $\partial\LLL$ separates $\eta_j$ from $\eta_{j-1}$. So $S_j\subset \overline A(\eta_j, \eta_{j-1})\subset \MMM\smm\Q$. But $S_j$ can not have
other boundary curves due to Condition (**).
Therefore  $S_j=\overline A(\eta_j, \eta_{j-1})$.

 It follows $S_2=S_3=\overline A(\eta_2, \eta_3)$, and more generally $S_j=S_{j+1}=
 \overline A(\eta_j, \eta_{j+1})$ for any even number $j$ with $2\le j< m$.

 Fix $j=1,\cdots,m$. We have $\beta_j\subset F^{-p}(\LLL)\subset \LLL$. Let $S'$ be the $\LLL$-piece containing $\beta_j$. We want to show that $S'$ is one of $S_i$.

 If $\beta_j\subset D(\eta_1)$ then $S'=D(\eta_1)=S_1$. Otherwise $S'$ has a boundary
 component $\eta$ separating $a$ from $\beta$. Therefore $\eta$ is one of $\eta_i$ and
 $S'=S_i$.

 Let now $j$ be an even number with $2\le j<m$. We know that $\beta_j\subset D(\eta_j)$
 and $\beta_j\subset S_i$ for some $i$. Therefore $S_i\subset \overline D(\eta_{j-1})$.
  But $\beta_j\cap \eta_{j-1}=\emptyset$. So $\beta_j\subset D(\eta_{j-1})$.
  Furthermore, $\overline A(\beta_j,\beta_{j+1})$ is a component of $F^{-p}(S_j)$,
  so must be contained entirely in the $\LLL$-piece $S_i$. Therefore $\beta_{j+1}\subset
  \overline D(\eta_{j-1})$ and consequently $\beta_{j+1}\subset D(\eta_{j-1})$.

  This ends the proof of the claim: $\overline D(\beta_j)\subset D(\eta_{j-1})$ for any $j\ge 2$.

 Note that in each component of $G^{-p}(D(\eta_m))$ the $G^p$-preimages of the
 discs $D(\eta_i)$, $D(\beta_j)$ are nested in the same order. There is therefore $q_j$,
 such that for any $n\ge q_j$,  all
curves of $G^{-np}(\eta_j)$ are either null-homotopic or
homotopic within $\overline D(\eta_1)\smm \{a\}$, therefore within $\LLL\smm \P$, to $\eta_1$, which is a curve in $\partial\RRR$.

Combining these arguments together, we find a $q$, such that every
curve in $F^{-q}(\GGG)$ is either null-homotopic or
homotopic within $\LLL\smm \P$ to a curve in $\partial\RRR$. So $(W_Z)^q=O$.

Therefore
$\la(W_Z)=0$ and $\la(W_X)=\la(W_T)$. \qed

{\noindent\bf 2. Curves in $X$ that are homotopic to each other within $\MMM\smm \Q$}. Now we decompose
$X=T\smm Z$ into $X_1\sqcup\cdots\sqcup X_k$ according to the homotopy
within $\MMM\smm \Q$, i.e., two curves in $X$ are homotopic within $\MMM\smm \Q$ if
and only if they belong to some subset $X_i$. Pick one curve
$\g_i$ in each $X_i$ and set $\G:=\{\g_1, \cdots, \g_k\}$.
Clearly $\G$ is a multicurve  within $\MMM\smm\Q$ .

Set $D_{\G}:=(b_{ij})$ the $(G,\Q)$-transition matrix  of $\G$. Set $W_X:=(a_{\de\be})$. By definition $$
b_{ij}=\sum_{\al\in\GGG_{ij}} \frac 1{\deg(G:\al\to\g_j)}\quad
\text{and}\quad a_{\de\be}=\sum_{\al\in\FFF_{\de\be}}\frac
1{\deg(F:\al\to\be)}, $$ where $\GGG_{ij}$ is the collection of
curves in $G^{-1}(\g_j)$ homotopic to $\g_i$ within $\MMM\smm \Q$;
$\FFF_{\de\be}$ is the collection of curves in $F^{-1}(\be)$
homotopic to $\de$  within $\LLL\smm\P$. We claim that,
$$\forall\ i,j\in\{1,\cdots,k\},\ \forall \be\in X_j,\quad \sum_{\de\in X_i} a_{\de\be}\le b_{ij}\ .$$

Assume at first $\be=\g_j$. We have
$$\begin{array}{lcl}
\bigcup_{\de\in X_i}\FFF_{\de\g_j}&=&\{\eta\in F^{-1}(\g_j),
\ \exists\ \de\in X_i\ \text{such that $ \eta,\de$ are homotopic within $\LLL\smm\P$}  \}\\
&\subset &\{\eta\in F^{-1}(\g_j),
\ \exists\ \de\in X_i\ \text{such that $ \eta,\de$ are homotopic within $\MMM\smm\Q$}  \}\\
&= &\{\eta\in G^{-1}(\g_j),
\ \exists\ \de\in X_i\ \text{such that $ \eta,\de$ are homotopic within $\MMM\smm\Q$}  \}\\
&= &\{\eta\in G^{-1}(\g_j),
\text{$ \eta,\g_i$ are homotopic within $\MMM\smm\Q$}  \}\\
&=& \GGG_{ij}\ ,
\end{array}
$$
where the inclusion is due to (a), and the second equality is due to
$F^{-1}(\g_j)=G^{-1}(\g_j)$. Therefore
$$
\sum_{\de\in X_i} a_{\de\g_j}=\sum_{\de\in X_i} \sum_{\al\in\FFF_{\de\g_j}}\frac
1{\deg(F:\al\to\g_j)}=\sum_{\al\in \bigcup_{\de\in X_i}\FFF_{\de\g_j}}\frac
1{\deg(F:\al\to\g_j)}$$ $$\le \sum_{\al\in \GGG_{ij}}\frac
1{\deg(G:\al\to\g_j)}=b_{ij}\ .
$$

 This implies the claim for
$\be=\g_j$.

When $\be\neq\g_j$, replace $\g_j$ by $\be$ in $\G$.
The replacement does not change the transition matrix $D_{\G}$. So
the claim is still true.

Applying now Corollary \ref{mp1} from linear algebra, we have $\la(W_X)\le \la(D_{\G})$.

 But $\la(D_{\G})<1 $ as $(G,\Q)$
has no Thurston obstructions. Consequently $\la(W_T)=\la(W_X)<1$. So $(F,\P)$ has no Thurston obstructions.
\qed

\section{Constant complexity under pullback}\label{cm}

We are now searching for repelling systems with some specific properties such that on one hand  we are
capable to solve their problem of c-equivalence to holomorphic models, and on the other hand
they appear as sub-systems of any \rs. This leads to an important class
of repelling systems: those \cm. We will prove at first that every \rs\  contains a sub-system that is \cm.
We then introduce two particular types of obstructions for this class of maps, and state our Thurston-like theorem in this setting, Theorem \ref{sb}: a \rs\ \cm\ without these specific obstructions is
c-equivalent to a holomorphic model. The proof of Theorem \ref{sb} will be postponed to the next sections.
We conclude the present section with a proof of Theorem \ref{more-general} using Theorem
 \ref{sb}.

\subsection{Definitions}

Let $(\E\stackrel{F}{\to}\LLL,\P)$ be a marked \pf\ \rs.
Let $S$ be an $\LLL$-piece.
We say that  $S$ is {\em simple} if either $S$ is annular with $S\cap\P=\emptyset$, or $S$ is disc-like with $\#S\cap\P\le 1$. Otherwise we say that $S$ is {\em complex}, i.e. if $\#\{\text{curves in }\partial S\}+\#(S\cap \P)\ge 3$.

More generally, let
$E\subset \LLL$ be a \rc. We say that $E$  is {\em simple} if $E$ is contained in  either a  closed annulus $A$ in
$\LLL\smm \P$; or in a closed disc $D$ in $\LLL$,
such that $D$  contains zero or one point of $\P$ in its interior, and that $\partial D\cap \P=\emptyset$. Otherwise we say $E$
is {\em complex}.

Constant complexity means roughly that under the pull-back by $F$,
both the number and the (homotopic) shapes of the complex
$\LLL$-pieces remain stable. More precisely:

{\noindent\bf Definition 5.} Let $(\E\stackrel{F}{\to}\LLL,\P)$ be a marked \pf\ \rs. We say that $(F,\P)$ is
{\bf \cm}, if every
complex $\LLL$-piece $S$, if any, contains an $\E$-piece
$E_S$ such that $E_S\cap \P=S\cap \P$ and that the components of $S\smm E_S$ are either annular or disc-like (this implies in particular both $\P_F,\P$ are contained in $K_F$).

Such $E_S$ is said to be
{\em parallel} to $S$. One way to obtain a parallel subsurface of $S$ is as follows:
first thicken the boundary of $S$ (without touching $\P$) to
reduce $S$ to a sub-bordered surface $E'$, then dig a few open holes
compactly contained in  $ interior(E')\smm \P$, the result is a bordered surface
$E$ parallel to $S$ (see Figure \ref{structure}).

For an example one may take $g(z)=z^2-1$. Cut off a suitable
neighborhood of $\P_g=\{\infty,0,-1\}$ to obtain a puzzle
neighborhood $\LLL$ of the Julia set so that
$g^{-1}(\LLL)\subset\subset \LLL$. In this case $\LLL$ has only one
piece $S$, which has three boundary curves, and  $E_S=g^{-1}(\LLL)$
has four boundary curves. Now $g: g^{-1}(\LLL)\to\LLL$ is a repelling system \cm. For details and further examples, see \S
\ref{examples}.

\subsection{Achieving constant complexity via restriction}\label{connection}

\REFTHM{Repellor} Let $(\BBB\stackrel{G}{\to} \MMM,\Q)$ be a marked \pf\ \rs with $\partial \MMM\ne \emptyset$. Then there are  two \ps\ neighborhoods $\E,\LLL$ of $K_G$ satisfying:
\begin{enumerate}\item[
{\rm (*)}] $K_G\subset\subset \E\subset\subset\LLL\subset\subset \MMM$, $\E=G^{-1}(\LLL)$, $\Q\cap K_G=\Q\cap \LLL$;
\item[
{\rm (**)}] for every $\LLL$-piece $S$,
 and for the $\MMM$-piece $S_0$ containing $S$ in the interior, and for the copy $\cbar$ of the Riemann sphere containing $S_0$ (therefore $S$),  every (disc-like) component of $\cbar\smm S$ contains either components of $\partial S_0$ or points of $ \Q$ (or both);
\item[
{\rm (***)}] the sub-\rs\  $F=G|_{\E}:\E\to\LLL$,  marked by $\Q\cap K_G$,
is \cm.\end{enumerate}
 \ENDTHM

 To prove the theorem, we need the following process together with its two properties:

{\noindent\bf Hole-filled-in process}. Let $S_0$ be an $\MMM$-piece. It is contained in a Riemann sphere $\cbar$. Let $E\subset S_0$ be a \rc. The {\it filled-in} of $E$, denote by $\widehat E$, is
defined to be the union of $E$ with all components of $E^c=\cbar\smm E$ disjoint from $\partial S_0\cup \Q$. Clearly, $\widehat E\subset S_0$ and $\partial\widehat
E\subset\partial E$.

{\noindent\bf Monotonicity Property}. Let $E_1\subset\subset E_2$ be
 two nested bordered connected surfaces in $\MMM$. Then $\widehat E_1\subset\subset\widehat
E_2$.

{\noindent\it Proof}. This property is easier to understand from
their complements. There is an $\MMM$-piece $S_0$ containing both $E_1$ and $E_2$.
Note that $E_2^c$ is a disjoint union of discs
while $(\widehat E_2)^c$ is the union of some of them which meet
$\partial S_0\cup \Q$. Because $E_1^c\supset\supset E_2^c$, these discs
are compactly contained in $E_1^c$, which can not be thrown away
under filled-in process of $E_1$ since they meet $\partial S_0\cup \Q$.
So $(\widehat E_1)^c\supset\supset(\widehat E_2)^c$ and hence $\widehat
E_1\subset\subset\widehat E_2$.

{\noindent\bf Pull-back property}. Let $S\subset \MMM$ be a
\rc\ with $\partial S\cap \Q=\emptyset$, and let $E_1$ be a component of $G^{-1}(S)$. Let
$\tilde E_1$ be the component of $G^{-1}(\widehat S)$ containing
$E_1$. Then $E_1$ is again a \rc, and $\tilde E_1 \subset\widehat E_1$. See the following diagram:
$$\begin{array}{ccrcl} && \La_0\subset \widehat S & \stackrel{\text{fill}}{\supset} & S\\
&& G\uparrow && \uparrow G \\
\La_1\subset \widehat E_1=fill(E_1) & \stackrel{\text{conclusion}}{\supset} & \tilde E_1& \supset & E_1
  \end{array}$$

{\noindent\it Proof}. Let $\La_0$ be the $\MMM$-piece containing $S$ and $\La_1$ be the $\MMM$-piece containing $E_1$. Assume $\La_0\subset \cbar_0$ and $\La_1\subset \cbar_1$. Denote by $\tilde E_i$ ($1\le i\le k$) the
components of $G^{-1}(\widehat S)$. Noticing that $\widehat S\smm S$ is a
union of disjoint open discs, is contained in $\La_0$ and is disjoint from $ \Q$. Set $V=G^{-1}(\widehat
S)\smm G^{-1}(S)=G^{-1}(\widehat S\smm S)$. Then $V$ is also a disjoint union
of discs and $\overline V \cap (\partial \MMM\cup \Q)=\emptyset$ (by the Basic pull-back principle). These discs are contained in
$G^{-1}(\widehat S)=\cup\tilde E_i$ and hence can not separate any
$\tilde E_i$, i.e. $\tilde E_i\smm V$ is also connected for $1\le
i\le k$. Therefore $E_1=\tilde E_1\smm V$. Note that each
component of $\tilde E_1\cap V$ is the union of some components of
$V$ which are discs contained in $\La_1\smm \Q$. They are also components
of $\La_1\smm E_1$. By definition, $\tilde E_1=E_1\cup(\tilde
E_1\cap V)\subset\widehat E_1$.

\medskip

{\noindent\it Proof of Theorem \ref{Repellor}}.

{\noindent\bf Choice of $N'$ to stabilize the postcritical set}.
Clearly there is an integer $N_0\ge 0$ such that for all $n\ge N_0$,
we have $L_n\cap\Q=K_G\cap\Q$,  in other words  every critical point of $G$ in
$L_{n}$ is actually in $K_G$ and is eventually periodic. For convenience we will
choose $N'\ge \max\{1,N_0\}$.

{\noindent\bf Choice of $N''$ to stabilize  the homotopy classes
of the boundary curves}. Note that $\partial L_n\cap
\partial L_{m}=\emptyset$ if $n\neq m$. For any integer $m\ge 0$, we
consider the homotopy classes within $\MMM\smm \Q$ of the Jordan
curves in $\bigcup_{k=0}^m \partial L_k$. The number of these
homotopy classes is weakly increasing with respect to $m$, but is
uniformly bounded from above, as $\Q\cup\partial \MMM$
has only finitely many connected components. There is
therefore an integer $N''\ge N'$ such that for any $n\ge N''$,
every boundary curve of $L_n$ is either null homotopic or
homotopic, within $\MMM\smm \Q$, to a curve in
$\bigcup_{k=0}^{N''-1}\partial L_k$.

{\noindent\bf Filled-in for $L_n$}. For any two $L_n$-pieces,
their filled-in are either disjoint or one contains another. Let
$\LLL_n$ be the union of the filled-in of all the $L_n$-pieces. Then
each $\LLL_n$-piece is the filled-in of an $L_n$-piece; for
every $\LLL_n$-piece $S$, each complementary component of $S$
contains points of $\partial \MMM\cup \Q$ (i.e. $\LLL_n$ satisfies Property (**)). Note that the total number of $\LLL_n$-pieces might be less than that
of $L_n$-pieces, as some $L_n$-piece might be hidden in the hole of another, thus
disappears in the filled-in process.

 It is
easy to check from the definition that if $S$ is a $G$-complex $L_n$-piece, then $\widehat S$ is a $G$-complex
$\LLL_n$-piece.

Assume that $E$ is a component of $G^{-1}(S)$ where $S$ is an
$L_n$-piece. Let $\tilde E$ be the component of $G^{-1}(\widehat S)$
containing $E$. Then $\tilde E\subset\widehat E\subset\subset\LLL_n$ by
the pull-back property and the monotonicity property of filled-in,
as $E$ is an $L_{n+1}$-piece and hence is contained in an
$L_n$-piece. Combining with the fact that each $\LLL_n$-piece is the
filled-in of an $L_n$-piece, we have
$G^{-1}(\LLL_n)\subset\LLL_{n+1}\subset\subset \LLL_n$. Note that
$G^{-1}(\LLL_n)\cap\Q=\LLL_n\cap\Q=K_G\cap \Q$ for all $n\ge N'$.

{\noindent\bf Choice of $N$ to stabilize the number and the shape
of the complex pieces}. From now on we assume $n\ge N''$. We claim
that for $k\ge 1$, each non-null-homotopic (within $\MMM\smm\Q$) curve $\gamma$ on
$G^{-k}(\partial\LLL_n)$ is homotopic to a curve on $\partial\LLL_n$
within $\MMM\smm\Q$. Note that $\gamma\subset
G^{-k}(\partial\LLL_n)\subset G^{-k}(\partial L_n)=\partial
L_{n+k}$. By the stability of the homotopy classes of boundary
curves, there is an integer $m$ with $m<N''\le n$ so that $\gamma$ is
homotopic (within $\MMM\smm\Q$) to a curve $\beta$ on $\partial L_m$.
Because $L_{n+k}\subset\subset L_n\subset\subset L_m$, there
exists a curve $\alpha$ on $\partial L_n$ so that $\alpha$
separates $\beta$ from $\gamma$. So $\alpha$ is also homotopic
(within $\MMM\smm\Q$) to $\gamma$. Let $S$ be the $L_n$-piece containing
$\alpha$. The fact that $\alpha$ is not null-homotopic implies
that $\widehat S$ is a $\LLL_n$-piece and $\alpha\subset\partial\widehat S$.
The claim is proved.

Let $S$ be a $G$-complex  $\LLL_n$-piece. Assume that $E_1$
and $E_2$ are $G^{-1}(\LLL_n)$-pieces in $S$ and that $E_1$ is
$G$-complex. Then there is a closed curve $\gamma$ on
$\partial E_1$ such that $\gamma$ separates $E_1\smm\gamma$ from
$E_2$. If $\gamma$ is null-homotopic within $\MMM\smm\Q$, then $E_2$ is
simple and $E_2\cap\Q=\emptyset$. Assume that $\gamma$ is
non-null-homotopic within $\MMM\smm\Q$. From the above claim, we have a
curve $\alpha$ on $\partial\LLL_n$ such that $\alpha$ is homotopic
to $\gamma$ within $\MMM\smm\Q$. Moreover, $\alpha$ can be taken on
$\partial S$ since $E_1\subset S$. Now the closed annulus enclosed
by $\gamma$ and $\alpha$, denote by $A(\gamma, \alpha)$, is
disjoint from $\partial \MMM\cap \Q$. It contains either $E_1$ or $E_2$ because
$\gamma$ separates $E_1\smm\gamma$ from $E_2$. Because $E_1$ is
$G$-complex, we see that $E_2\subset A(\gamma, \alpha)$
and hence $E_2$ is $G$-simple and
$E_2\cap\Q=\emptyset$.

The above argument shows that $S$ contains at most one
 $G^{-1}(\LLL_n)$-piece that is $G$-complex. In case that $S$ contains a
$G$-complex  $G^{-1}(\LLL_n)$-piece $E_S$, other
$G^{-1}(\LLL_n)$-pieces in $S$ are simple and disjoint from $\partial \MMM\cap \Q$.
Combining with the fact that $G^{-1}(\LLL_n)\cap\Q=\LLL_n\cap\Q$,
we see that $E_S\cap\Q=S\cap\Q$.

We can show now that each
 component of $E^c_S$ contains at most one component of
$\partial S$. Let $D$ be a component of $E_S^c$ and
$\gamma=\partial D\cap\partial E_S$. If $\gamma$ is null-homotopic
within $\MMM\smm\Q$, then $D$ contains no component of $\partial S$ since
each closed curve in $\partial S$ is non-null-homotopic rel
$\Q$. Now assume that $\gamma$ is non-null-homotopic rel
$\Q$. Then there is a curve $\beta$ on $\partial S$ homotopic to $\gamma$ within $\MMM\smm \Q$. Therefore the closed annulus $A(\gamma, \beta)$ bounded by $\gamma$ and $\beta$
is contained in $\MMM\smm \Q$.

If $\beta\subset
D^c$, then $E_S\subset A(\gamma, \beta)$. This contradicts to the
fact that $E_S$ is $G$-complex. So $\beta\subset D$. Thus
$A(\gamma, \beta)\subset S$ since $S$ is the filled-in of an
$L_n$-piece. Therefore no other components of $\partial S$ is
contained in $D$. This implies that components of $S\smm E_S$ are either
annular or disc-like.

Let $s_n$ be the number of $G$-complex  $\LLL_n$-pieces. Let
$t_n$ be the number of $G$-complex  $G^{-1}(\LLL_n)$-pieces.
Then $t_n\le s_n$ since each $G$-complex $\LLL_n$-piece
contains at most one $G$-complex $G^{-1}(\LLL_n)$-piece. We
claim that $s_{n+1}\le t_n$.

Let $\widehat E_1$ and $\widehat E_2$ be distinct $G$-complex
$\LLL_{n+1}$-pieces, where $E_1$ and $E_2$ are $L_{n+1}$-pieces.
Then $E_1$ and $E_2$ are also $G$-complex  by the
definition. Note that $L_{n+1}\subset
G^{-1}(\LLL_n)\subset\LLL_{n+1}$. We have two distinct
$G^{-1}(\LLL_n)$-pieces $\tilde E_1$ and $\tilde E_2$ such that
$E_i\subset\tilde E_i\subset\widehat E_i$ ($i=1,2$). Again $\tilde
E_1$ and $\tilde E_2$ are also $G$-complex. So
$s_{n+1}\le t_n$.

Now we have $s_{n+1}\le s_n$. There is therefore an integer $N\ge
N''$ such that $s_n\equiv s_N$  for $n\ge N$.

Define now a new \rs\  $F:\E\to \LLL$ to be $G|_{G^{-1}(\LLL_N)}: G^{-1}(\LLL_N)\to
\LLL_N$. It is \pf\ with $\P_F=\P_G\cap K_G\subset\Q\cap K_G$, and $K_F=K_G$.
Furthermore  $(F,\Q\cap K_G)$ is \cm.
 \qed

\subsection{Boundary curves and complex pieces}

We now turn to the study of properties of constant complexity maps.

\REFLEM{stable} Let $(\E\stackrel{F}{\to}\LLL,\P)$ be a marked \pf\ \rs\ \cm.\begin{enumerate}\item
 For any $n\ge 0$, any curve in $F^{-n}(\partial\LLL)$ is either
null-homotopic or homotopic to a curve in $\partial\LLL$ within $\LLL\smm \P$. \item
 For any complex $\LLL$-piece $S$, there is a unique $\E$-piece $E_S$ parallel to $S$,
and $F(E_S)=:S'$ is again a complex $\LLL$-piece. \item
 $F_*:S_1\to S_2$ if $F(E_{S_1})=S_2$ is a well defined map
from the set of complex $\LLL$-pieces into itself. Every such $\LLL$-piece is eventually periodic under $F_*$. \item
 For any complex $\LLL$-piece $S$ and any integer $m\ge 1$, there
is a unique  $F^{-m}(\LLL)$-piece $E$ in $S$
parallel to $S$. Moreover, $F^m(E)$ is again a complex $\LLL$-piece.\end{enumerate}
\ENDLEM

Before proving it we will decompose $\LLL$ following its topology
and its intersecting property with $\P$. Let $S$ be a $\LLL
$-piece. We say that $S$ is of $\AAA$-type if $S\cap\P=
\emptyset$ and $S$ has exactly two boundary curves; is of
$\OOO$-type if $S\cap\P= \emptyset$ and $S$ has exactly one
boundary curve; is of $\RRR$-type if $S\cap\P$ is a single point
and $S$ has exactly one boundary curve; is of $\CCC$-type if
$\#(S\cap\P)+\#\{\text{boundary curves}\}\ge 3$ (see table \Ref{S}).
Note that an $\LLL$-piece is complex iff it is a $\CCC$-piece.

We now decompose $\LLL$ into $\CCC\sqcup\RRR\sqcup \A\sqcup\OOO$
with $\CCC$ the union of $\CCC$-type pieces, $\RRR$ the union of
$\RRR$-type pieces, $\AAA$ the union of $\AAA$-type pieces and
$\OOO$ the union of $\OOO$-type pieces. It may happen that some
sets among $\OOO,\A,\RRR, \CCC$ are empty. See examples in \S
\ref{examples}.

\REFEQN{S}
\text{\begin{tabular}{|l|c|c|c|}
\multicolumn{4}{c}{\bf Classification of $\LLL$-pieces $S$}\vspace{0.1cm}\\
\hline
 & 1 boundary curve & 2 boundary curves & $\ge 3$ boundary curves \\
$\cap\P\backslash\text{shape}$ & (disc) & (annulus) & (pants, pillowcase \\
&&& without corners, etc.)\\
\hline
$S\cap \P=\emptyset$ &  $\OOO$-type &  $\AAA$-type & $\CCC$-type \\
$\g\subset\partial S$ is & null-homotopic &non-peripheral &non-peripheral \\
&&$\partial S=\partial_-S\sqcup \partial_+S$&\\
\hline \!\!\begin{tabular}{c}
$S\cap\P\ne\emptyset$ \\
     \\ $\g\subset\partial S$ is
\end{tabular}
&\!\!\begin{tabular}{c|c} $\#S\cap\P=1$ & $\#S\cap\P>1$ \\
$\RRR$-type & $\CCC$-type
    \\ peripheral & non-periph. \end{tabular}  &
\!\!\begin{tabular}{c} \\ $\CCC$-type \\
non-peripheral \end{tabular}
 & \!\!\begin{tabular}{c} \\ $\CCC$-type \\
non-peripheral\end{tabular} \\
\hline
\end{tabular}}\ENDEQN

{\noindent\em Proof of Lemma \ref{stable}.} (1). Due to the basic pull-back principle we just need to
prove it for $n=1$. Let $\g$ be a boundary curve of $\LLL$. Then
$\g$ is a boundary curve of some $ \LLL$-piece $S$.

If $S$ is of $\OOO$-type, then all components of $F^{-1}(S)$ are
discs in $\LLL$ and are disjoint from $\P$. Therefore all curves
in $F^{-1}(\g)$ are null-homotopic.

Recall that by the definition \cm, each $\CCC$-piece $S'$ contains
a unique complex $\E$-piece $E_{S'}$ and $E_{S'}$ is parallel to
$S'$.

If $S$ is of $\AAA$-type or $\RRR$-type, then each component $E$
of $F^{-1}(S)$ is contained in \REFEQN{ss}
\OOO\sqcup\AAA\sqcup\RRR\sqcup\bigcup_{S':\
\CCC-\text{piece}}S'\smm E_{S'}\ . \ENDEQN But each component of
$S'\smm E_{S'}$ for a $\CCC$-piece $S'$ is either an annulus or a
disc, and is contained in $\LLL\smm\P$. So each boundary curve of
$E$, in particular each curve of $F^{-1}(\g)$, is either
null-homotopic or homotopic to a curve in $\partial\LLL$.

Finally if $S$ is of $\CCC$-type, then a component $E$ of
$F^{-1}(S)$ is either equal to $E_{S'}$ for some $\CCC$-piece
$S'$, or is contained in \Ref{ss}. In any case each boundary curve
of each of $E$, in particular each curve of $F^{-1}(\g)$, is
either null-homotopic or homotopic to a curve in
$\partial\LLL$.

(2). The existence of $E_S$ is given by the definition \cm. Its uniqueness follows from
the fact that components of $S\smm E_S$ are annular or disc-like and are disjoint from $\P$.
We know that  $S':=F(E_{S})$ is again an
$\LLL$-piece. It must be also an $\CCC$-piece since each component of the $F$-preimage of a simple
$\LLL$-piece is also simple.

(3). Clearly $F_*$ is well defined due to (2). Since the number of $\CCC$-pieces are
finite, each of them is eventually periodic under $F_*$.

(4). Let $S$ be a $\CCC$-piece. We have seen from (3) that
$S':=F(E_S)$ is again a $\CCC$-piece. By definition \cm, we know
that $S'=E_{S'}\sqcup(\sqcup_iA_i\sqcup_j D_j)$, with $A_i$ annuli
and $D_i$ discs, and that $(\sqcup_iA_i\sqcup_j D_j)\cap
\P=\emptyset$, and that one component of $\partial A_i$ is
contained in $\partial S'$.

There are no complex $F^{-2}(\LLL)$-pieces in $\sqcup_iA_i\sqcup_j
D_j$. By the basic pull-back principle, $E^2_S:=E_S\cap
F^{-1}(E_{S'})$ is connected. It is reduced from $E_S$ after
thickening the boundary and then cutting off a few disjoint holes
(without touching $\P\cap E_S$). This implies that $E^2_S$ is
parallel to $S$. So $E^2_S$ is complex. Clearly, there is no other
complex $F^{-2}(\LLL)$-piece in $E_S$.

Inductively, for any integer $m\ge 1$, there is a unique complex
$F^{-m}(\LLL)$-piece $E^m_S$ in $S$ and $E^m_S$ is parallel to $S$.
Moreover, $F^m(E^m_S)$ is again a $\CCC$-piece. \qed

\subsection{The boundary multicurve}

Let $(\E\stackrel{F}{\to}\LLL,\P)$ be a marked \rs.

We consider now  the boundary
curves $\g$ of $\LLL$ that are non-peripheral (within $\LLL\smm\P$) (i.e. either $\g$ is the boundary of a disc piece $D$, with $\#D\cap \P\ge 2$, or $\g$ is
a boundary curve of a non-disc piece). This set might be empty, or some of the curves might be
homotopic to each other (for example the two boundary curves of an annular component of $\LLL\smm \P$). In any case we give the

{\noindent\bf Definition 6}. A {\em boundary multicurve} $Y$ of
$(F,\P)$ is a collection of curves in $\partial\LLL$ representing all the
homotopy classes within $\LLL\smm \P$ of the non-peripheral curves in
$\partial\LLL$.

Then {\em boundary transition matrix}\ $W_Y=(a_{ij})$
is defined by \REFEQN{WY} a_{ij}=\sum_\al\frac1{\deg(F:\al\to
\g_j)},\ENDEQN where the sum is taken over all the Jordan curves
(if any) of $\al\subset F^{-1}(\g_j)$ that are
homotopic to $\g_i$ (within $\LLL\smm\P$).

We will say that $(F,\P)$ has a {\bf boundary
obstruction} if $Y\ne\emptyset$ and $\la(W_Y)\ge 1$.

In general $Y$ is not $(F,\P)$-stable, or even worse, we might have $Y\ne\emptyset$ and $W_Y$ equal to the
zero matrix.

Assume from now on that $(F,\P)$ is also \cm. Then $Y$ is $(F,\P)$-stable by Lemma \ref{stable}.
It can be described more explicitly using the above classification of $\LLL$-pieces:

A closed curve in $\partial\LLL$ is null-homotopic iff it is
contained in $\partial\OOO$, is peripheral iff it is contained
in $\partial\RRR$. Two closed curves in $\partial\LLL$ are
homotopic iff they are the boundary curves of an $\AAA$-piece. Label
by $+$ and $-$ the two boundary curves for each $\AAA$-piece $S$,
i.e. $\partial S=\partial_{+}S\sqcup\partial_{-}S$. Let
$\partial_{+}\AAA$ and $\partial_{-}\AAA$ denote the union of the
corresponding boundary curves of every $\AAA$-pieces. Then the
boundary multicurve $Y$ can be taken to be the collection of
closed curves in $\partial\CCC\cup\partial_{+}\AAA$.

Polynomials with all critical points escaping to $\infty$ provide
examples, when restricted to a suitable neighborhood of the Julia
set, of repelling systems with $Y=\emptyset$. For an annuli covering
$F:\E\to\A$ with $\A=A_1\sqcup\cdots\sqcup A_k$ (see Lemma
\ref{annulicase}), the set $Y$ consists of $k$ boundary curves of
$\A$, one in each $A_i$, and the $(F,\P)$-transition matrix $W_Y$
coincides with the transition matrix $D$ defined in \S
\ref{obstructions}. For further examples see \S \ref{examples}.

\subsection{Renormalizations and renormalized obstructions}

A repelling system \cm\ has another, somewhat more important property: it admits renormalizations
and they behave like \pf\ branched coverings of $\cbar$.

{\noindent\bf Definition 7}. A marked \pf\ \rs\ $(E\stackrel{H}{\to}S, \tilde \P)$ \cm\ is of {\em
Thurston type} if both $E$ and
$S$ are connected, and $$\#(S\cap \tilde \P)+\#\{\text{boundary curves of }S\}\ge 3.$$ In other words, $S,E\subset \cbar $ are quasi-circle \rc s, $E$ is compactly contained in the interior of $S$, the components of $S\smm E$ are either annular or disc-like, $H:E\to S$ is an
orientation preserving branched covering, with a finite (or empty) postcritical set $\P_H$
which is contained in $E$, the set
$\tilde\P\subset E$ is a finite (or empty) set containing both $H(\tilde\P)$ and $\P_H$,
and $\#(S\cap \tilde \P)+\#\{\text{boundary curves of }S\}\ge 3$.

Again, an example can be provided by
the map $g(z)=z^2-1$, with $S$ equal to $\cbar$ minus a suitable
neighborhood of $\infty,0,-1$, and with $E=g^{-1}(S)$.

Let $(\E\stackrel{F}{\to}\LLL,\P)$ be a marked \pf\ \rs\ \cm. Assume $\CCC\ne\emptyset$. By Lemma
\ref{stable}, we have a map $F_*$ defined on the collection of
$\CCC$-pieces by $F_*(S_1)=S_2$ if $F(E_{S_1})=S_2$ where
$E_{S_1}$ is the unique $\E$-piece in $S_1$ parallel to $S_1$. Assume that
$S$ is an $\CCC$-piece that is $p$-periodic under $F_*$. Let $E$ be the unique
$F^{-p}(\LLL)$-piece in $S$ parallel to $S$. Then $F^p(E)=S$. Denote $H=F^p|_{E}$.
Then $H:E\to S$ is a \rs\ satisfying that
$\P_H\subset\P_F\cap S\subset \P\cap S$, and that, marked by $\P\cap S$, the marked system  $(H,\P\cap S)$ is of Thurston type.

{\noindent\bf Definition 8}. We will call the marked \rs\  $(E\stackrel{H}{\to} S,\P\cap S)$  a {\bf renormalization} of
$(F,\P)$.  We say that $(F,\P)$ has a {\bf
renormalized obstruction} if it has a renormalization $(E\stackrel{H}{\to} S,\P\cap S)$ that has a Thurston
obstruction.

\REFLEM{Prepare} Let $(\E\stackrel{F}{\to}\LLL,\P)$ be a marked \pf\ \rs\ \cm. If $(F,\P)$ has no Thurston obstructions, then it has no
boundary obstructions nor renormalized obstructions.
\ENDLEM
\beginp As $(F,\P)$ has no Thurston obstructions, we have $\la(W_T)<1$ for the transition matrix $W_T$ of every
multicurve $T$ in $\LLL\smm \P$, in particular for $T$ equal to the boundary multicurve. Therefore $(F,\P)$ has no boundary obstructions.

It remains to show that any renormalization $H:E\to S$ marked by $\tilde\P:= \P\cap S$ has no Thurston
obstructions.  Assume by contradiction that $\la(W_\G)\ge 1$, for some $(H,\tilde\P)$-stable multicurve $\G$, with $(H,\tilde\P)$-transition matrix
$W_\G$.

Let $S_0(=S), S_1,\cdots, S_{p-1}$ be the $F_*$-periodic cycle of $S$. Let $E_i$ be the unique $\E$-piece in $S_i$
parallel to $S_i$, $i=0,1,\cdots, p-1$. Set $\G_p:=\G$. Define inductively, for $i=p-1,p-2,\cdots, 0$, the multicurve $\G_i\subset F^{-1}(\G_{i+1})\cap E_i$ representing the homotopy
classes (within $S_i\smm \P$) of the non-peripheral curves in $F^{-1}(\G_{i+1})\cap E_i$. By stability of $\G$, each curve of $\G_0$ is homotopic to a curve in $\G_p$.

Consider $F':\bigcup_iE_i\to \bigcup_iS_i$ as a \rs.
Set $\G'=\G_1\cup\cdots \cup \G_{p}$. It is a multicurve within $\bigcup_iS_i\smm \P$,
therefore within $\LLL\smm \P$.
Denote by $W_{\G'}$ its $(F,\P)$-transition matrix and by $W'$ its $(F',\P)$-transition matrix. Then the $p$-th power $(W')^p$ restricted to $\G_p$ is equal
to $W_\G$. Therefore
$$1\le \la(W_\G)\le \la((W')^p)=\la(W')^p\ .$$
But each entry of $W'$ is less than or equal to the corresponding entry of $W_{\G'}$.
Therefore $\la(W')\le \la(W_{\G'})$. This implies that $\la(W_{\G'})\ge 1$.
This contradicts the assumption that $(F,\P)$ has no Thurston obstructions (by Lemma \ref{non-stable}).
\qed

We can now state our Thurston-like result in this setting, whose proof will occupy Sections \ref{bridge}-\ref{I}:

\REFTHM{sb} Let $(\E\stackrel{F}{\to}\LLL,\P)$ be a marked \pf\ \rs\ \cm. Assume that $(F,\P)$ has no boundary obstructions nor renormalized obstructions. Then
$(F,\P)$
 is c-equivalent to a holomorphic
model.\ENDTHM

\subsection{Proof of Theorem \ref{more-general} using Theorems \ref{sb}}\label{Fin}

{\noindent\em Proof of Theorem \ref{more-general}}.
Let $(\BBB\stackrel{G}{\to} \MMM,\Q)$ be a marked \pf\ \rs\ without Thurston obstructions.
 We will prove that $(G,\Q)$ is c-equivalent to a holomorphic model.

As first we apply Theorem \ref{Repellor}  to $(G,\Q)$ to show that
it has  a \pf\ \rs\ restriction $F:\E\to \LLL$ near $K_G$ which, marked by $\P:=\Q\cap K_G$, is  \cm, and satisfies the conditions in Theorem \ref{Unobstructed}.
So we may apply Theorem \ref{Unobstructed}.(B)  to show that $(F,\P)$ has no
Thurston obstructions. Lemma \ref{Prepare} then leads $(F,\P)$ to the setting
of Theorem
\ref{sb}, i.e. $(F,\P)$ is \cm, and has no boundary obstructions nor renormalized obstructions. Now we may apply Theorem \ref{sb} to conclude that $(F,\P)$
is c-equivalent to a holomorphic model. Finally we conclude for $(G,\Q)$
using  Theorem \ref{Unobstructed}.(A).\qed

Note that it could be  more difficult to check the condition of
Theorem \ref{hyperbolic} and Theorem \ref{more-general}, namely $G$ has no Thurston obstructions.
Whereas Theorem \ref{sb} turns it into the problem of checking
the leading eigenvalue of $W_Y$ for a single multicurve $Y$,
and then the absence of Thurston obstructions for \pf\ branched
coverings (arising from the renormalizations), to which there is a huge literature (see for example
the references in [ST]). This form is particularly suitable for
combination of rational maps, i.e. starting with \pf\ rational maps
(thus already holomorphic) as the renormalizations and glue them
suitably together.

\section{C-equivalence to holomorphic models}\label{bridge}

From now on we concentrate on the proof of Theorem \ref{sb}: a marked \pf\ \rs\ \cm\
without boundary obstructions nor renormalized obstructions is c-equivalent to a holomorphic
model. In this section we will prove the theorem for the non-renormalizable case.
In this case only Gr\"otzsch inequalities are needed, but not the original Thurston theorem.

Recall that from Definition 4 and Lemma \ref{5'} that a marked \rs $(\E\stackrel{F}{\to}\LLL,\P)$
is c-equivalent to a holomorphic model, if there is a marked \rs $(\E'\stackrel{R}{\to}\LLL',\P')$
with $R$ holomorphic, and two quasi-conformal homeomorphisms $\Theta:\LLL\to \LLL$, $\Phi:\LLL\to \LLL'$  with
$$\left\{\begin{array}{l} \Phi\circ \Theta(\E)=\E',\ \Phi\circ \Theta(\P)=\P'\vspace{0.2cm}\\
R\circ \Phi\circ \Theta|_{\E} = \Phi\circ F \vspace{0.2cm}\\
\Theta \text{\ is isotopic to the identity rel\ } \P\cup \partial \LLL.\end{array}
  \right.$$

\subsection{Examples}\label{Examples}

{\noindent\bf Example 1}.\quad
\begin{minipage}{4.5cm}
{\mbox{}\hspace{-0.8cm}\epsfig{width=4.5cm,height=4cm,figure=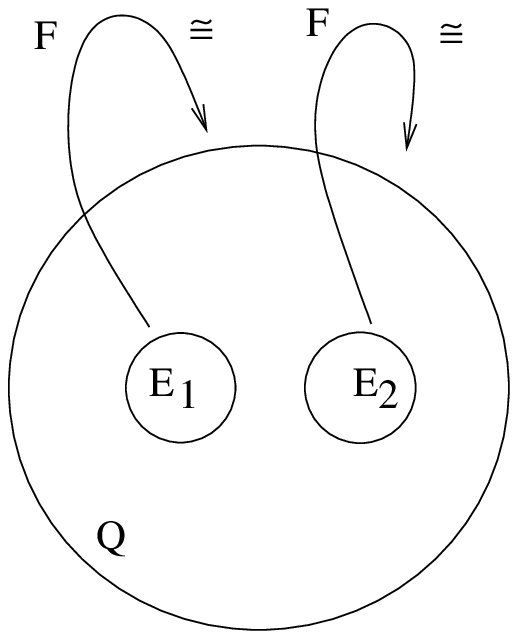}}%
\end{minipage}
\begin{minipage}{8.5cm}
{ $\LLL=Q$, $\E=E_1\cup E_2$, with $Q,E_1,E_2$ quasi discs.}%
{\ \\ $F:E_i\to Q$ are quasi-conformal homeomorphisms}%
\\ $\P=\emptyset$.
\end{minipage}
\mbox{}\\[-0.5cm]\hspace*{-1.2cm}\mbox{}

{\noindent\bf Example 2}.\quad
\begin{minipage}{6cm}
{\mbox{}\hspace{-0.4cm}\epsfig{width=6cm,height=3cm,figure=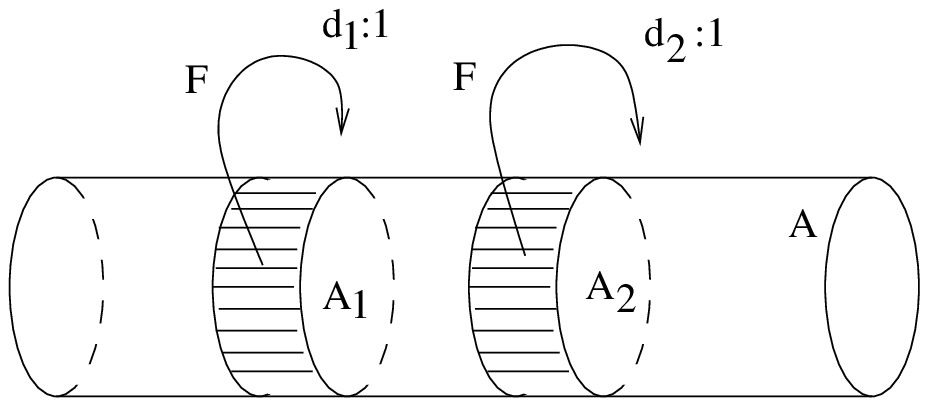}}%
\end{minipage}
\begin{minipage}{6.5cm}
{$\LLL=A$, $\E=A_1\cup A_2$, \\with $A,A_1,A_2$ closed annuli. }%
{\ \\ $F:A_i\to A$ are degree $d_i$\\ quasi-regular coverings.}%
\\ $\P=\emptyset$.
\end{minipage}

\REFLEM{examples} The map $F$ in Example 1 is always c-equivalent to a holomorphic
model, whereas in Example 2 it is so if and only if $\frac1{d_1}+\frac1{d_2}<1$.
\ENDLEM
\beginp Let $F: E_1\cup E_2\to Q$ be as in Example 1.
\begin{enumerate}\item
 Construct at first the model map $R$ by setting $E_i'=E_i$, $Q'=Q$ and by choosing $R:E_i'\to Q'$ any
conformal homeomorphism.
\item Set $\Phi=id: Q\to Q$. \item Set $\left\{ \begin{array}{lcl}
\theta|_{E_i} & =& R^{-1}\circ F \\
\theta|_{\partial Q} & = & id.
\end{array}
 \right.$\quad Then $R\circ \Phi\circ \theta|_{E_i} = \Phi\circ F$. \item Extend $\theta$ as a homeomorphism of $Q$.  \item One checks easily that this $\theta$ is isotopic to the identity rel $\partial Q$.\end{enumerate}

 Let now $F:A_1\cup A_2\to A$ be a map in Example 2. Note that $A$ has a unique multicurve $\G$, up to homotopy,
 with $\G$  consisting
 of a boundary curve of $A$. Its $F$-transition matrix has only one entry, which is $\frac1{d_1}+\frac1{d_2}$.
 Therefore $F$ has a Thurston obstruction if and only if $\frac1{d_1}+\frac1{d_2}\ge 1$.

 If $F$ is c-equivalent to a holomorphic $R:A_1'\cup A_2'\to A'$. Then by Gr\"otzsch inequality $\frac1{d_1}+\frac1{d_2}<1$.

 Conversely assume  $\frac1{d_1}+\frac1{d_2}<1$.  \\
{\epsfig{width=10.5cm,height=3.5cm,figure=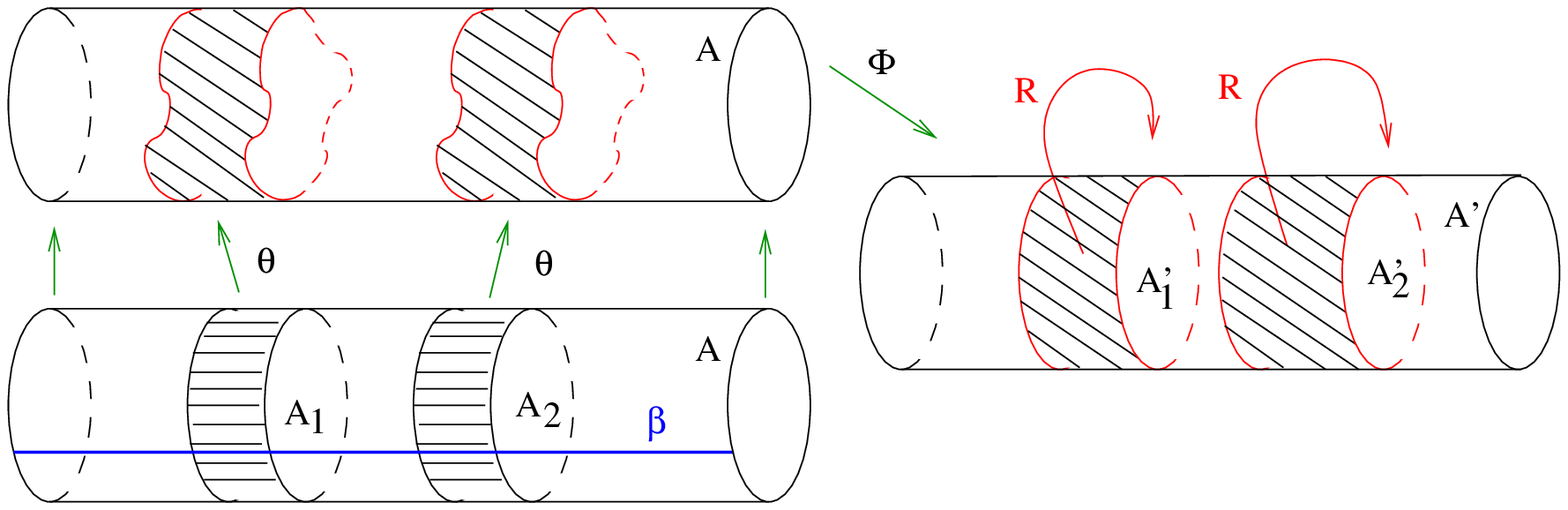}}
\begin{enumerate}\item[a.]
a. Construct at first a round modulus $A'$ of modulus, say $v$. And let $\Phi: A\to A'$ be a qc-homeomorphism.
\item[b.]
Construct then two disjoint essential round
submoduli $A_1',A_2'$ in $A'$ of moduli $v/d_1$ and $v/d_2$ respectively, and displaced in the same order
as the $A_i$'s in $A$. This is possible precisely because $\sum \frac1{d_i}<1$.
Choose $R: A_i'\to A'$ a holomorphic covering of degree $d_i$, matching the boundary correspondence as $F$.
\item[c.] Pull-back $A_i'$ by $\Phi$
\item[d.] Set $\left\{ \begin{array}{lcl}
\theta|_{A_i}:  & =& (\Phi^{-1}\circ R\circ \Phi)^{-1}\circ F : A_i\,{ \longrightarrow}\, \Phi^{-1}(A_i)\\
\theta|_{\partial A} & = & id
\end{array}
 \right.$
\item[e.] Extend $\theta$ as a qc-homeomorphism $A\to A$. Then { $R\circ \Phi\circ \theta|_{A_i}=\Phi\circ F$.}
\item[f.] Via {\em Dehn twist} on $A\smm (A_1\cup A_2)$ to modify the extension so that $\theta({ \beta})\sim {
\beta}$.\end{enumerate}
 This guarantees that $\theta$ is isotopic to the identity rel $\partial
A$.
\qed

\medskip

{\noindent\bf Example 3}. Let $L$ be a pair of trousers bounded by three quasi-circles $\g_0,\g_{-1},\g_*$.
Let $E\subset L$ be a surface bounded by four quasi-circles $\beta_0,\beta_{-1},\beta_*,\beta_1$,
with $\beta_1$ bounding a complementary disc of $E$ that is entirely contained in $L$, and with each
other $\beta_i$ bounding a complementary disc of $E$ that contains the corresponding $\g_i$.
Let $H:E\to L$ be a quasi-regular covering of degree $2$. Again $\P_H=\P=\emptyset$. And $H$ is of Thurston type.

The boundary
multicurve $Y$ is simply $\{\g_0,\g_{-1},\g_*\}$. For example if we require $H:\beta_*\to \g_*$
to be degree $2$, $H:\beta_{\pm 1}\to \g_0$ of degree $1$ and $H:\beta_{0}\to \g_{-1}$ of degree $2$,
then the transition matrix $W_Y$ is
$\left(\begin{matrix} 0&\frac12&0 \\ 1 & 0&0\\0&0&\frac12\end{matrix}\right)$ with leading eigenvalue $1/\sqrt{2}$.
Such a map behaves like $z^2-1$.

On the other hand, if we require instead $H:\beta_*\to \g_*$
to be degree $2$, $H:\beta_{\pm 1}\to \g_1$ of degree $1$ and $H:\beta_{0}\to \g_{0}$ of degree $2$,
then the transition matrix $W_Y$ is $\left(\begin{matrix} \frac12&0&0\\0&1&0\\0&0&\frac12
                                       \end{matrix}\right)$ with leading eigenvalue $1$.
This $H$ can be constructed more explicitly as follows (suggested by X. Buff):
Let $g(z)=z^2$. Let $L$ be $\cbar$ minus a small round disc of radius $\ep$ around each of the three points $0,\infty$ and $1$.
Let $E'=g^{-1}(L)$. As $1$ is a repelling fixed point of $g$, $E'$ is not compactly contained
in $L$ and $g:E'\to L$ is not a \rs.
Let now $\eta: \overline D(1,2\ep)\to \overline D(1,2\ep)$ to be a homeomorphism
fixing pointwisely the boundary and center, mapping the boundary component of $E'$
into the interior of $L$. Extend $\eta$ elsewhere by identity.
Set $H=g\circ \eta^{-1}: \eta(E')\to L$.

Note that every non-peripheral curve is homotopic to one of the curves in $Y$. So $H$ has a Thurston
obstruction if and only if $\la(W_Y)\ge 1$.

 \epsfig{width=11.3cm,height=6cm,figure=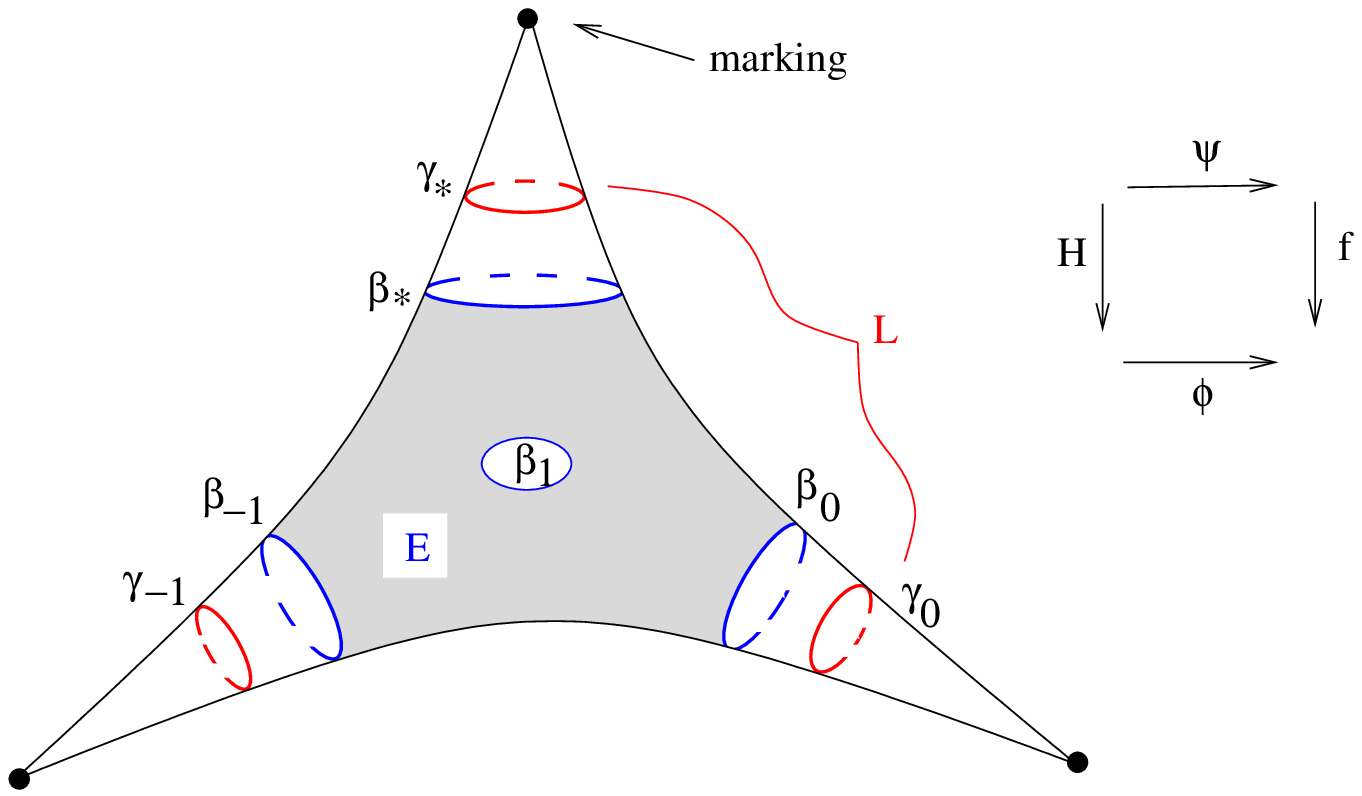}\vspace{-0.1cm}
\REFLEM{three}
This $H$ is c-equivalent to a holomorphic model  if and only if
the $H$-transition matrix of $Y=\{\g_0,\g_{-1},\g_*\}$ has a leading eigenvalue $\la<1$.\ENDLEM
\beginp Mark one point in each complement component of $L$.
Denote the marked set by $\P$. Extend $H$ as a quasi-regular branched cover $\widehat H$ of $\cbar$ such that
the critical values are in $ \P$ and $\widehat H(\P)\subset \P$, and that $\widehat H$ is
holomorphic outside $L$. In particular $\widehat H$ is \pf.

Assume at first the $H$ is c-equivalent to a holomorphic model. Then by Lemma \ref{5'} there is a
quasi-conformal automorphism $\theta: L\to L$, isotopic to the identity rel $\partial L$, and a
Beltrami differential $\mu$ supported on $L$ with $\|\mu\|_\infty< 1$ , such that $(H\circ \theta^{-1})^*\mu=\mu|_{\theta(E)}$.

Proceed now at in the proof of Proposition \ref{final}. Extend $\theta$ to be equal to the identity outside $L$.
Set $\widehat H_1=\widehat{H}\circ \theta^{-1}$.
Extend $\mu$ outside $L$ by $\mu=0$. Let $\phi_1:\cbar\to
\cbar$ be a global integrating map of this extended $\mu$. Set
$\widehat H_2:=\phi_1\circ \widehat H_1\circ \phi_1^{-1}$. Then $\widehat H_2$ is again
quasi-regular, and is holomorphic in the interior of
$\phi_1\circ\theta(E)$ and in $\phi_1(L)$. Elsewhere each $\widehat H_2$-orbit
passes at most once.
We spread out the
Beltrami differential $\nu_0\equiv 0$ using iterations of $\widehat H_2$ to
get an $\widehat H_2$-invariant Beltrami differential $\nu$. Note that
$\nu=0$ on $\phi_1(L)$, and $\|\nu\|_\infty<1$.  Integrating $\nu$
by a \qc\ homeomorphism $\phi_2$ (necessarily holomorphic on
$\phi_1(L)$), we get a new map $f:= \phi_2\circ \widehat H_2\circ  \phi_2^{-1}$
which is a rational map and is \ce\ to $\widehat H_2$, therefore to $\widehat H$. See
the following diagram.
$$\begin{array}{cclclclc} (\cbar,E) &
\stackrel{\theta}{\longrightarrow}&\cbar&
\stackrel{\phi_1}{\longrightarrow} & \cbar &
\stackrel{\phi_2}{\longrightarrow} & \cbar &
\\
\widehat H\downarrow && \downarrow \widehat H_1 &&\downarrow \widehat H_2&&\downarrow f\\
(\cbar, L) & \stackrel{id}{\longrightarrow} & \cbar &
\stackrel{\phi_1}{\longrightarrow} & \cbar &
\stackrel{\phi_2}{\longrightarrow} & \cbar &
\end{array} $$
Set $\phi=\phi_2\circ \phi_1$. Then $\phi(\partial L)$ is contained in the attracting basins of $f$. But $f$ is
\pf, each attracting cycle must be superattracting. It follows that the boundary multicurve
of the \rs $f|_{\phi\circ \theta(E)}$ has leading eigenvalue $<1$. Therefore $\la(W_Y)<1$ for our map $H$.

Conversely, assume that $\la(W_Y)<1$ for $H$.
The fact that $\#\P=3$ implies that $(\widehat H,\P)$ has no Thurston obstructions. By Thurston theorem  there are $\Phi=id$, $\Psi$ homeomorphism, and $f$ a rational map, such that
$f\circ \Psi=\Phi\circ \widehat H$, and that $\Psi$ is isotopic to $\Phi$ rel $\P$.

Let $v$ be a positive vector so that $W_Yv<v$. Change $\Phi$ within its isotopy class
to $\phi$ so that $\phi(\g)$, for any $\g\in Y$, is an equipotential in a Fatou component $\Sigma$ of $f$. Further, the annular component of $\Sigma\smm \phi(\g)$ is of modulus $v(\g)$. Let $\psi$ be the homeomorphism isotopic to $\Psi$ so that $f\circ \psi=\phi\circ \widehat H$.
Now $W_Yv<v$ assures that $\psi(E)$ is compactly contained in $\phi(L)$.

Use Lemma \ref{topology} to modify $\psi$ on $L\smm E$ so that $\psi$ is isotopic to $\phi$ rel $\cbar\smm L$, and $f\circ \psi|_E=\phi|_L\circ H$.
\qed

{\noindent\bf Example 4}. The following $F:\bigcup E_*\to S\sqcup A$ is an unbranched \rs\ \cm.\\
\begin{minipage}{11cm}%
{\epsfig{width=11cm,height=5.5cm,figure=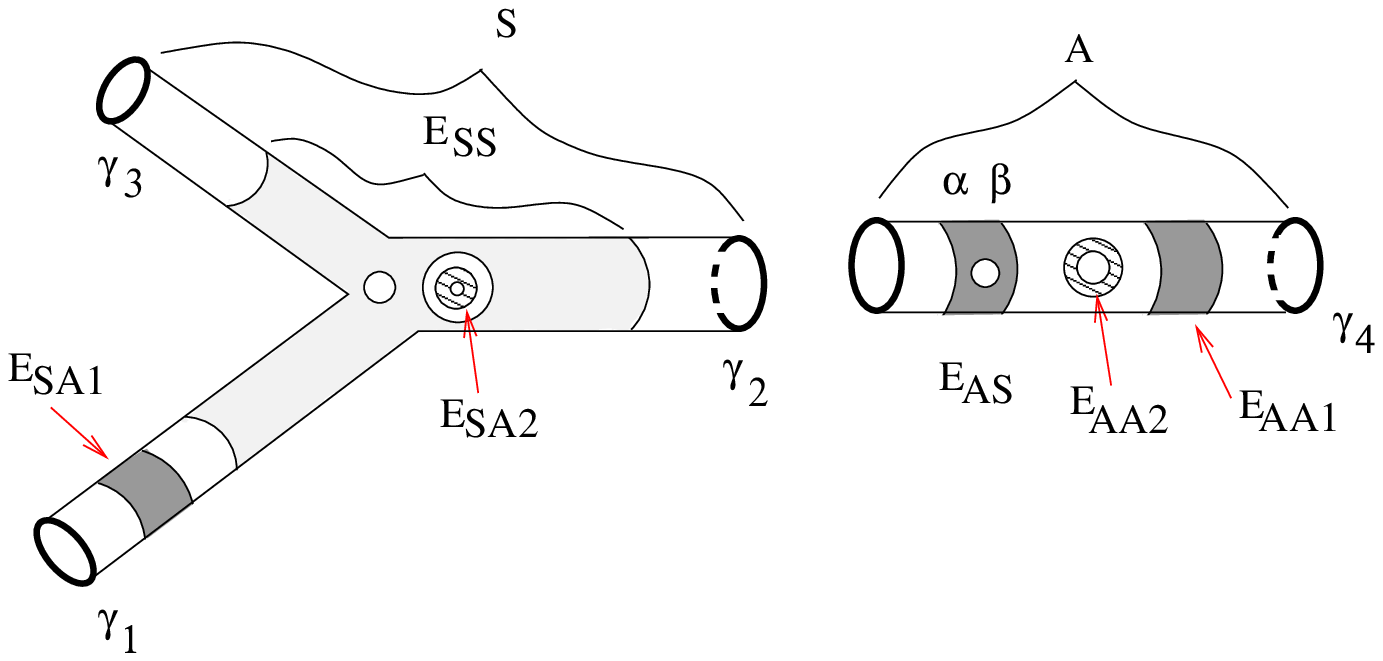}}%
\end{minipage}

\mbox{}\\[-2.4cm]\hspace*{-2.4cm}\mbox{}

\hspace*{6cm}
\begin{minipage}{11.5cm}
{$\LLL=S\sqcup A$,\ $S$ is a pair of pants, $A$ is an annulus}%
{\\ $F:E_{SS}\to S,\ E_{AS}\to S$ are coverings}%
{\mbox{}\\ $F:E_{*Aj}\to A$ are coverings, $\P=\emptyset$.}
\end{minipage}

\medskip

In this case the boundary multicurve $Y$ consists of the three boundary curves $\g_1,\g_2,\g_3$ of $S$ together
with one of the two boundary curves, named by $\g_4$ of $A$. And $F$ has only one renormalization $(E_{SS}\stackrel{F}{\to} S)$.

\REFLEM{four} If the $F$-transition matrix $W_Y$ of $Y$ has leading eigenvalue less than $1$,
then  $F$ is c-equivalent to  a holomorphic model.\ENDLEM
{\noindent\em Sketch of the proof}. Let $v=(v_1,v_2,v_3,v_4)$ be a vector with each $v_i>0$ such that $W_Yv<v$. Set $H=F|_{E_{SS}}$. Now the boundary multicurve is simply $\{\g_1,\g_2,\g_3\}$, the
its $H$-transition matrix $D_*$ is a submatrix of $W_Y$, therefore $D_*u<u$ for $u=(v_1,v_2,v_3)$.
 We may then, as in  Example 3, construct  $\phi,\psi$ making $H$ c-equivalent to a holomorphic model $f$, so that $\phi(\g_i)$ has the potential prescribed by $v_i$, $i=1,2,3$.

  But the
presence of $E_{AS}$ requires further control of $\phi,\psi$.

Fix $M>0$ a large number. Modify again $\phi,\psi$ (but not $f$) according to
the vector $Mu$.

Let now $A'$ be a round annulus of modulus $Mv_4$.

Fill in the hole in $E_{AS}$ to get $\hat E_{AS}$ that is an essential annulus is $A$.
Denote by $\al,\beta$ the boundary curves of $\hat E_{AS}$ .

 Assume, say, $F$ maps $\al$ to $\g_1$ with degree $d_1$, and $\beta$
 to $\g_2$ with degree $d_2$. Then a candidate $\hat E'_{AS}$ has a modulus
 bounded by $\frac{Mv_1}{d_1}+\frac{Mv_2}{d_2} + C$ with a constant $C$ independent of $M$
 (due to Lemma \ref{Inequality}). Choose $M$ so that $C<M\inf_i (v_i-(W_Yv)_i)$.

Set then $E'_{SS}=f^{-1}(S')=\psi(E_{SS})$.

Now $(W_Y(Mv))_i+C<Mv_i$  guarantees that one can insert two disjoint essential annuli $\hat E'_{AS}, E'_{AA1}$ in
$A'$, and insert an essential annulus $E'_{SA1}$ in the
corresponding annulus of $S'\smm E'_{SS}$, together with a model holomorphic covering $R: E'_{*Aj}\to A', E'_{*S}\to S'$
with the same degree as $F$ and the same boundary correspondence.
Then the construction of $\theta$ is similar.

Details are provided in  \S\ref{I}.\qed

\subsection{A criterion for c-equivalence to holomorphic models}

The following simple remark turns concept of c-equivalence
into a more practical form:

\REFLEM{struc} Let $(\E\stackrel{F}{\to}\LLL, \P)$ be a marked
\pf\ \rs. Then $(F,\P)$ is \ce\ to a holomorphic model if and only
if: for each $\LLL$-piece $S$, there is a pair of quasi-conformal
homeomorphisms $\theta_S: S\to S$, $\phi_S: S\to S'\subset\cbar_{S'}$ (here
we consider $\cbar_{S'}$ as a distinct copy of the Riemann
sphere),
such that: \\
(a) $\theta_S$ is the identity on $\partial S\cup\P$ and is
isotopic to the identity rel $\partial S\cup\P$. \\
(b) For every $\E$-piece $E$ contained in $S$, and for $\wtq=F(E)$
(which is again an $\LLL$-piece), the composition
$R_E:=\phi_{\wtq}\circ F\circ \theta_S^{-1}\circ \phi_S^{-1}$ is
holomorphic in the interior of $\phi_S\circ\theta_S(E)$. See
\Ref{e}.\ENDLEM

The proof is almost straight forward. One just need to set $\theta=\Phi^{-1}\circ \Psi$
in the definition of c-equivalence, and set $\theta_S=\theta|_S$, $\phi_S=\Phi|_S$
and $R_E=R|_{\phi_S\circ\theta_S(E)}$.

Therefore to prove that $(F,\P)$ is
c-equivalent to a holomorphic model, it amounts to construct
 $\phi_S,\theta_S$ and $R_E$ for each $\LLL$-piece $S$  and each $\E$-piece $E$ in $S$ so that they satisfy
(a) and (b). In practical the maps $\phi_S, R_E$  will  be
constructed first. One constructs then each $\theta_E=\theta_S|_E$, and finally glue the various
$\theta_E$'s together to get $\theta_S$. See the following schema:

\REFEQN{e}
\forall\ \E\text{-piece}\ E\subset S,\
\begin{array}{rcccl}\multicolumn{5}{l}{\text{\bf Order\ of\ the\ construction}\vspace{0.2cm}} \\S &\underset{{\bf 4.}}{\stackrel{\theta_S}{\longrightarrow}}& S
&\underset{{\bf 1.}}{\stackrel{\phi_S}{\longrightarrow}} & S'\subset \cbar_{S'}\\
\cup&&\cup&&\cup \\
 E &\underset{{\bf 3.}}{\stackrel{\theta_S|_E=\theta_E}{\longrightarrow}}& \underset{{\bf 2.}}{\widetilde E} &
 \stackrel{\phi_S}{\longrightarrow} & E' \\
F \downarrow &&&&\downarrow\,_{\bf 1.}\  R_E\text{\rm\ holomorphic}\\
 \wtq&\stackrel{id}{\longrightarrow}&\wtq& \underset{{\bf 1.}}{\stackrel{\phi_{\wtq}}{\longrightarrow}}
 &\wtq'\subset \cbar_{\wtq'} \end{array}\ENDEQN

\subsection{Annuli coverings}\mylabel{baby}

 Let now $f:\E\to \A$ be an annuli covering. In other words,
$\A=A_1\sqcup \cdots\sqcup A_k$,  $\E=\bigsqcup_{ij\de} E_{ij\de}$, with each $A_i$ (and $E_{ij\de}$) an annulus,  and with $E_{ij\de}\subset\subset A_i$ for $i,j\in\{1,\cdots, k\}$ and $\de$ in
some finite or empty index set $\La_{ij}$ depending on $(i,j)$, and  $f:E_{ij\de}\to A_j$ is a
quasi-regular covering of certain degree, denoted by $d_{ij\de}$. Recall that the transition matrix $D$ is defined by
$$D=(a_{ij}),\ a_{ij}=\sum_{\de\in \La_{ij}}\frac1{d_{ij\de}}\ .$$
We will prove the following more concrete form of Lemma \ref{annulicase}:

\REFLEM{Annuli} For the annuli covering $f:\E\to \A$ as above, assume that there is a vector $v=(v_1,\cdots, v_k)$ with $v_i>0$ for any $i$ such that
$Dv<v$, i.e. \REFEQN{6}\forall\ i=1,\cdots,k,\quad \sum_{j,\de}\frac{
v_j}{d_{ij\de}} <v_i\ .\ENDEQN
Then $f:\E\to \A$ is c-equivalent to a holomorphic model $R:\E'\to\A'$ with
$\A'=A'_1\sqcup \cdots\sqcup A'_k$ and $\mod(A'_i)=v_i$.\ENDLEM
 Here the modulus of a closed annulus
will mean the modulus of its interior as an open annulus.
Now Lemma \ref{matrix} relates $\la(D)<1$ to the existence of such vector $v$. And Lemma \ref{annulicase} follows.

\medskip

{\noindent\em Proof of Lemma \ref{Annuli}}.

We consider each $A_i$ as embedded in a
distinct copy $\cbar_i$ of the Riemann sphere. Take one more copy
$\cbar'_i$  for each
$i=1,\cdots,k$.

For each $i\in \{1,\cdots,k\}$, we will construct a pair
$(\theta_i,\phi_i)$ such that (see \Ref{theta}):

(a) $\theta_i: A_i\to A_i$ is a \qc\ map, with $\theta_i|_{\partial A_i}=id$ and with $\theta_i$ isotopic to the identity rel $\partial A_i$.
The set $\theta_i(E_{ij\de})$ is denoted by $\widetilde E_{ij\de}$ for each
possible $(j, \de)$.

(b) $\phi_i:A_i\to A_i'$ is a \qc\ homeomorphism. The set  $\phi_i(\widetilde E_{ij\de})$ is denoted by $E_{ij\de}'$
for each
possible $(j, \de)$.

(c) For each multi index $ij\de$, the map $R_{ij\de}:=\phi_j\circ f\circ
(\phi_i\circ \theta_i)^{-1}|_{E_{ij\de}'}$ is holomorphic in the
interior of $E_{ij\de}'$.

\REFEQN{theta}\begin{array}{lcccl}A_i
&\stackrel{\theta_i}{\longrightarrow}&
A_i &\stackrel{\phi_i}{\longrightarrow}& A_i'\subset \cbar_i'\vspace{0.1cm}\\
\cup && \cup && \cup \\
E_{ij\de}&\stackrel{\theta_i=\theta_{ij\de}}{\longrightarrow}&\widetilde
E_{ij\de}
&\stackrel{\phi_i}{\longrightarrow} &E_{ij\de}'\vspace{0.1cm}\\
\downarrow f &&&&\downarrow R_{ij\de} \text{\ holomorphic}\\
A_j&\stackrel{id}{\longrightarrow}&A_j&\stackrel{\phi_j}{\longrightarrow}
& A_j'\subset \cbar_j'\end{array}\ENDEQN
We will follow the order of construction as indicated by \Ref{e}.
Once this is done we can apply Lemma \ref{struc} to conclude that
$f:\E\to \A$ is c-equivalent to a holomorphic model.

{\noindent\bf 1. Definition of} $(\phi_i,A_i',E_{ij\de}',R_{ij\de})$:
 For every $i\in \{1,\cdots,k\}$, choose $A_i'\subset\cbar'_i$ a closed
round annulus of modulus $v_i$, and let
$\phi_i:A_i\to A'_i$ be a \qc\ homeomorphism.

 Fix an
index $i$. For every possible choice of $(j,\de)$, the lower diagram of \Ref{theta}
indicates that $R_{ij\de}: E'_{ij\de}\to A_j'$ is a covering isomorphic to
$\phi_j\circ f:E_{ij\de}\to A'_j$, therefore is
an annuli covering of degree $d_{ij\de}$. But $R_{ij\de}$ is
supposed to be holomorphic and $\mod(A'_j)=v_j$. This imposes that
$E_{ij\de}'$ must be a sub-annulus of $A'_i$ with modulus
$v_j/d_{ij\de}$.

Choose now a closed
round essential annulus $E_{ij\de}'$ in $A_i'$ such that (1)
$E_{ij\de}'\cap \partial A_i'=\emptyset$; (2)
$\mod(E_{ij\de}')=v_j/d_{ij\de}$ and the ($E_{ij\de}'$)'s are mutually disjoint for all possible indices
$(j,\de)$ (this is possible precisely because of
\Ref{6}); (3) the
($E_{ij\de}'$)'s are displaced in $A'_i$ in the same order as the
($E_{ij\de}$)'s in $A_i$.

Choose now $R_{ij\de}: E_{ij\de}'\to A_j'$  a holomorphic covering
of degree $d_{ij\de}$ among the two round annuli, so that it
permutes the boundary curves in the same way as $f:E_{ij\de}\to
A_j$.

More precisely this can be done through {\em boundary labeling}:
for each $A_i$ choose a labelling by $+$ and $-$ for its  two
boundary curves. This induces a labelling by $\pm$ on the boundary
curves of each essential sub-annulus $E_{ij\de}$, so that
$\partial_-E_{ij\de}$ separates $\partial_-A_i$ to
$\partial_+E_{ij\de}$. Now use each $\phi_i$ to transport these
labellings to $\partial A_i'$ which then induces a labelling on each  $\partial E_{ij\de}'$. The covering
$f:E_{ij\de}\to A_j$ maps $\partial_-E_{ij\de}$ to one of
$\partial_\pm A_j$. We choose $R_{ij\de}$ so that it sends
$\partial_-E'_{ij\de}$ to $\phi_j(f(\partial_-E_{ij\de}))$, the
corresponding boundary component of $A_j'$.

{\noindent\bf 2. Definition of $\widetilde
E_{ij\de}$}: For any multi-index $(i,j,\de)$, set
$\widetilde E_{ij\de}:=\phi_i^{-1}(E_{ij\de}')$ (there are a
priori two ways to label its boundary curves, one as an essential
sub-annulus of $A_i$, one transported by $\phi_i^{-1}$ of the
labeling of $\partial E_{ij\de}'$, but these two labellings
actually coincide).

{\noindent\bf 3. Definition of $\theta_{ij\de}$}. For any multi-index $(i,j,\de)$,
let $\theta_{ij\de}: E_{ij\de}\to\widetilde E_{ij\de}$ be a (choice of a) lift
of the \qc\ map $\phi_j:A_j\to A_j'$ via the two quasi-regular
coverings of the same degree: $f|_{E_{ij\de}}$ and $R_{ij\de}\circ
\phi_i|_{\widetilde E_{ij\de}}$. It is a \qc\ map and preserves
the boundary labelling.
$$\begin{array}{rclcl} E_{ij\de} & \stackrel{\theta_{ij\de}}{\dashrightarrow} & \widetilde E_{ij\de}&\stackrel{\phi_i}{\longrightarrow} &  E'_{ij\de}\\
f\downarrow &&&& \downarrow R_{ij\de} \\
A_j && \underset{\phi_j}{\longrightarrow} && A'_j
  \end{array}$$

{\noindent\bf 4. Definition of $\theta_i$}. Fix an index $i$.
Define $\theta_i: A_i\to A_i$ to be a \qc\ map such that
$\theta_i|_{E_{ij\de}}=\theta_{ij\de}$ and $\theta_i|_{\partial
A_i}=id$. It exists always, because all the boundary curves are
quasi-circles and all $\theta_{ij\de}$ are quasi-conformal maps
preserving the boundary labelling (see Lemma \ref{open}).

The map $\theta_i$ satisfies all the required properties, except
possibly the one about their homotopy type.

{\noindent\bf 4'. Adjustment of the homotopy type of $\theta_i$}.
As the lower commuting diagram in \Ref{theta} only requires
information on $\theta_{ij\de}$, we will modify each $\theta_i$
without changing its value on the  $E_{ij\de}$'s.

Fix an index $i$. Choose an arc $\beta$ connecting the two
boundaries of $A_i$. Then $\theta_i(\beta)$ is again an arc in
$A_i$ with the same end points. Precompose $\theta_i$ with a \qc\
repeated Dehn twist supported in the interior of $A_i\smm
\bigcup_{(j,\de)}E_{ij\de}$ if necessary we can ensure that
$\theta_i(\beta)$ is homotopic to $\beta$ (rel $\partial A_i$).
After this adjustment $\theta_i$ is well isotopic to the identity
rel $\partial A_i$. \qed

\subsection{The non-renormalizable case}

{\em Proof of Theorem \ref{sb} in case $\CCC=\RRR=\emptyset$ and $\P=\emptyset$}.

Let $F:\E\to \LLL$ be a \rs. Assume that $\P_F=\P=\emptyset$ and that every $\LLL$-piece
is either annular or disc-like. Then $F$ is \cm\ and there is nothing to
renormalize. Furthermore the boundary multicurve $Y$ is simply the
collection of one boundary curve in each annular piece of $\LLL$. Assume now
$\la(W_Y)<1$. We want to prove that $F:\E\to \LLL$ is c-equivalent to
a holomorphic model. Set
$$\LLL=\A\sqcup \OOO,\quad \A=A_1\sqcup\cdots\sqcup A_k,\quad \OOO=O_{k+1}\sqcup\cdots \sqcup O_m$$
so that each $A_i$ is annular and each $O_j$ is disc-like.

{\noindent\bf 0. The vector $v$}. Choose a vector $v\in \R^k$ with every entry positive such that $W_Y(v)<v$.

{\noindent\bf 1. Definition of $(\phi_S,S', E',R_E)$}. Consider each $\A\sqcup\OOO$-piece $S$
as a subset of a distinct copy
$\cbar_S$ of the Riemann sphere.
If
$S=A_i$, define as above $\phi_S:S\to S'$ to be a quasi-conformal homeomorphism so that $S'$ is an
round annulus with modulus $v_i$. If $S=O_i$, set simply
$\phi_S=id$ and $S'=S$.

For any $\E$-piece $E$,  there are two $\A\sqcup\OOO$-pieces $S$
and $\wtq$ (possibly the same) such that $E\subset\subset S$ and
$F(E)=\wtq$. As $F$ is a covering, we know that $E$ is an annulus
(resp. disc) if and only if $\wtq$ is an annulus (resp. disc).

We decompose $\E$ into $\E^\A\sqcup \E^{0,\A}\sqcup \E^{0,0}$ as
follows:

$\bullet$ $\E^\A$ consists of $\E$-pieces that are essential
sub-annuli in $\A$. We numerate as above these pieces  by
$E_{ij\de}$ so that $E_{ij\de}\subset A_i$ and $F(E_{ij\de})=A_j$.
Label by $\pm$ the boundary curves of each $A_i$, give each
$\partial E_{ij\de}$ the induced labelling.

$\bullet$ $\E^{0,\A}$ consists of the remaining annular
$\E$-pieces. Such a piece $E$ is a closed annulus contained in
some $\LLL$-piece $S$, with either $S=O_i$ or $S=A_i$. In the latter
case $E$ is contained in a component $B$ of
$S\smm\bigcup_{(j,\de)}E_{ij\de}$, and $E$ is non-essential in
$B$. See Figure \ref{structure}. In any case both boundary curves
of $E$ are $F$-null-homotopic. Furthermore $E$ has one
complementary component $\De_E$ that is entirely contained in $S$.
We will call $\De_E$ the hole of $E$. We will label $\partial E$
so that $\partial_-E=\partial \De_E$. Therefore $\partial_+E$
denotes the outer boundary of $E$.

$\bullet$ Finally $\E^0$ consists of the disc pieces of $\E$. Such
a piece $E$ may be contained in $\A$ or in $\OOO$, but $F(E)$ is
always an $\OOO$-piece, and $F:E\to F(E)$ is a homeomorphism.

The general strategy is quite simple to explain: we should at
first construct $(E',R_E)$ for all the $\E^\A$-pieces as in Lemma
\ref{Annuli}. For the remaining $\E$-pieces, as they may be
nested in each others holes, we should construct $(E',R_E)$
inductively from outer pieces to inner pieces.

{\noindent\bf 1.1. Construction of $(E',R_E)$ for each
$\E^\A$-piece $E$}. For all possible $(i,j,\de)$, define
$E'_{ij\de}$ to be a round annulus in $A_i'$ exactly as in the
proof of Lemma \ref{Annuli} above, in particular the various
$E'_{ij\de}$ for a given $i$ are displaced in $A_i'$ exactly in
the same order as the $E_{ij\de}$'s in $A_i$ and
$\mod(E'_{ij\de})=v_j/\deg(F:E_{ij\de}\to A_j)$. This is where we
have used the assumption $\la(W_Y)<1$. Let $R_{ij\de}:
E'_{ij\de}\to A_j'$ be a holomorphic covering of degree
$\deg(F:E_{ij\de}\to A_j)$, and permuting the boundary labelling
in the same way as $F:E_{ij\de}\to A_j$.

{\noindent\bf 1.2. Construction of $(E',R_E)$ for each
$\E^{0,\A}$-piece $E$}. Let now $E$ be an $\E^{0,\A}$-piece. Then
$E$ is contained in a set $B$ which is either some $O_i$ or one
component of $A_i\smm\bigcup_{(j,\de)}E_{ij\de}$ for some $A_i$.
Furthermore $F$ maps $E$ onto some $A_j$ as a covering. Denote its
degree by $d_E$. Denote by $B'$ the corresponding $O_i'$ or the
corresponding component of $A'_i\smm  \bigcup_{j\de}E'_{ij\de}$.
We want to set $E'$ to be an innessential annulus in $B'$ with
modulus $v_j/d_E$, so that there is a holomorphic covering
$R_E:E'\to A_j'$ of degree $d_E$.

However one might run into some moduli difficulty if we do so
randomly, as various pieces of $\E^{0,\A}$ in $B$ may be nested in
each others holes. The correct way to do this is to place $E'$ one
by one from outside to inside. More precisely, numerate the
$\E^{0,\A}$-pieces in $B$ by $E_{Bj\de}$, with $\de$ running in
some index set (depending on $j$) such that $F(E_{Bj\de})=A_j$ and
it is a covering of degree $d_{Bj\de}$. Define then a layer
(depth) function $l(E)$ on the set of $(E_{Bj\de})_{j\de}$ as
follows: set $l(E)=1$ if $E$ is outermost, i.e. not contained in
the hole of any other $\E^{0,\A}$-pieces. Set inductively $l(E)=m$
if $E$ is contained immediately in the hole of a $E_{Bj\de}$ with
$l(E_{Bj\de})=m-1$.

Start now with an $E_{Bj\de}$ so that $l(E_{Bj\de})=1$. Choose
$E'_{Bj\de}\subset\subset B'$ to be any round inessential annulus
of modulus $v_j/d_{Bj\de}$. Label its outer boundary curve by $+$.
Choose then $R_{Bj\de}:E'_{Bj\de}\to A_j'$ a holomorphic covering
of degree $d_{Bj\de}$, so that it permutes the boundary labelling
in the same way as $F:E_{Bj\de}\to A_j$.

Construct similarly $(E',R_E)$ for every layer $1$ piece in $B$,
and be sure that the various $E'$'s are mutually disjoint.

Now we should construct $(E',R_E)$ for the next layer
$\E^{0,\A}$-pieces in $B$. Proceed this layer by layer. As each
time we are supposed to find finitely many disjoint annuli non
mutually nested of prescribed moduli in the hole of some
previously constructed $E'$, the construction is always
realizable.

Do this construction for every component $B$ of $\LLL\smm \E^\A$.

{\noindent\bf 1.3. Construction of $(E',R_E)$ for each
$\E^0$-piece $E$}. Assume $E\subset S$ and $F(E)=\wtq$. We should choose
 a closed quasi-disc $E'$ in $S'$ disjoint from
the previously constructed pieces, together with a conformal map
$R_E:E'\to \wtq'$. There is no difficulty here and we omit the
details.

To recapitulate we may extend the layer function $l(E)$ on all
$\E$-pieces so that $l(E)=0$ for $\E^\A$-pieces and $l(E)=+\infty$
for $\E^0$-pieces, and then construct $(E',R_E)$ following the
natural order of the layer function.

{\noindent \bf 2-3. Definition of $(\widetilde E,\theta_E)$}. This
is done exactly as in Lemma \ref{Annuli}, by setting
$\widetilde E=\phi_S^{-1}(E')$ for $S$ the $\LLL$-piece containing
$E$, and $\theta_E: E\to \widetilde E$ as a lift of
$\phi_{\wtq}:\wtq \to \wtq'$ via $F|_E$ and $R_E\circ
\phi_S|_{\widetilde E}$, where $\wtq=F(E)$.

{\noindent\bf 4. Definition of $\theta_S$}. Fix an $\LLL$-piece $S$.
We claim that we can define a \qc\ map $\theta_S: S\to S$ so that
$\theta_S=\theta_E$ on each $\E$-piece $E$ contained in $S$ and
$\theta_S=id$ on $\partial S$. Clearly the extension can be chosen
so that $\theta_S$ is a homeomorphism, as the $\theta_E$'s for all
possible $E$ preserve the boundary labelling. But all the boundary
curves are quasi-circles and all the $\theta_E$'s are
quasi-conformal. One can then apply Lemma \ref{open} to make the
extension quasi-conformal.

{\noindent\bf 4'. Adjustment of the homotopy type of
$\theta_{A_i}$}. Clearly $\theta_S$ for $S$ an $\OOO$-piece is
already isotopic to the identity rel $\partial S$. However for $S$
an $\A$-piece one might have to precompose $\theta$ with a repeated
Dehn twist supported on the interior of $S\smm\E$ as in Lemma
\ref{Annuli}. After that $\theta_S$ is also isotopic to the
identity rel $\partial S$. \qed

\section{Holomorphic model of a renormalization cycle}\label{prepare}

Let $(\E\stackrel{F}{\to} \LLL,\P)$ be a marked \pf\ repelling system \cm\
without boundary obstructions nor renormalized obstructions. We will prove here that a renormalization cycle of $(F,\P)$ is c-equivalent to a holomorphic model, which satisfies in addition some prescribed moduli properties.

We always denote by $\D$ the unit disc.
A {\em marked disc} is a pair $(\De,a)$ with $\De$ an open
hyperbolic disc in $\cbar$ and a marked point $a\in \De$.  An {\em
equipotential} $\g$ of $(\De,a)$ is a Jordan curve
that is mapped to a round circle under a conformal representation $\chi:
\De\to\D$ with $\chi(a)=0$. The {\em potential} of such a $\g$ is
defined to be $\kappa(\g):=\text{mod} A(\partial \De,\g)$, the
modulus of the annulus between $\partial \De$ and $\g$. These
notions do not depend on the choice of $\chi$. The map $\kappa$ maps the set of equipotentials bijectively
onto the open interval $]0,+\infty[$. For example in the
marked disc $(\D,0)$, the circle $\{|z|=e^{-v}\}$ is an
equipotential with potential $v$  (we define
$\text{mod}\{r<|z|<1\}:=-\log r$).

Let $f$ be a \pf\ rational map with non-empty Fatou set. The Julia set $\JJJ_f$ is connected and each
Fatou component $\De$ is canonically a marked disc marked by the unique
eventually periodic point $a$ in $\De$. We call $(\De, a)$ {\em a marked Fatou component of $f$}. The equipotentials of these marked Fatou components will be
called {\em equipotentials of $f$}. Notice that equipotentials
in a periodic Fatou component correspond to round circles in
B\"ottcher coordinates. We will use $\kappa$ to denote the potential function
of these marked discs.

{\noindent\bf Marked set $\P_S$}. Again
consider each $\LLL$-piece $S$  as being contained in a distinct copy $\cbar_S$ of
the Riemann sphere. Mark one point in each component of
$\cbar_S\smm S$. Set $\P_S$ to be the union of $\P\cap S$ with these marked
points. By definition a piece $S$ is an $\CCC$-piece if and only if $\#\P_S\ge 3$.

As $(F,\P)$ is \cm, there is an induced map $F_*$ on the set of $\CCC$-pieces.
Let $S_1,\cdots, S_p$ be a periodic cycle of $\CCC$-pieces for $F_*$, i.e. for $i=1,\cdots, p$
we have
 $F(E_i)=S_{i+1}$ (set $S_{p+1}=S_1$), where $E_i$ is the
unique complex $\E$-piece in $S_i$. Denote by
$\cbar_i=\cbar_{S_i}$ for simplicity.

 We will prove:

\REFTHM{Pp} Denote by $D_*$ the $(F,\P)$-transition matrix of the set of the boundary curves of $S_1\cup\cdots\cup S_p$. Let $u>0$ be any positive vector such that $D_*u<u$. Then there are pairs of \qc\ maps $(\phi_{S_i},
\psi_{S_i}):\cbar_{S_i}\to \cbar_{S_i'}$ and holomorphic maps $R_i:
\cbar_{S_i'}\to\cbar_{S_{i+1}'}$ for $1\le i\le p$
such that: \begin{enumerate}
\item[(a)] $\phi_{S_i}=\psi_{S_i}$ on $\partial S_i\cup(\P\cap {S_i})$, and $\phi_{S_i}$ is isotopic to $\psi_{S_i}$ rel $\partial
S_i\cup(\P\cap {S_i})$; \item[
(b)] $\phi_{S_{i+1}}\circ F\circ\psi_{S_i}^{-1}|_{\psi_{S_i}(E_i)}=R_i|_{\psi_{S_i}(E_i)}$; \item[
(c)] the return map $f_{S_i}:=R_{i-1}\circ\cdots\circ R_1\circ R_p\circ\cdots\circ
R_i$ is a \pf\ rational map whose conformal conjugacy class depends only on $F$ and $S_i$;\item[
(d)] for each $i\in \{1,\cdots,p\}$, for each Jordan curve $\g\subset\partial S_i$, and for $\beta_\g$
the curve in $\partial E_i$ homotopic to $\g$ within
$S_i\smm\P$, both $\phi_{S_i}(\g)$ and $\psi_{S_i}(\beta_\g)$ are
equipotentials in the same marked Fatou component of $f_{S_i}$ with
potentials \REFEQN{wp3} \kappa(\phi_{S_i}(\g))=u(\g)
\quad\text{and}\quad
\kappa(\psi_{S_i}(\beta_\g))=\frac{u(F(\beta_\g))}{\deg(F|_{\beta_\g})}. \ENDEQN
\end{enumerate}
\ENDTHM

Note that (a) and (b) together assert that  $R: \bigcup \psi_{S_i}(E_i)\to \bigcup \phi_{S_i}(S_i)$, marked by $\bigcup \phi_{S_i}(\P\cap S_i)$, with $R=R_i$ on $\psi_{S_i}(E_i)$,
is a holomorphic \rs\ c-equivalent to $(\bigcup E_i\stackrel{F}{\to} \bigcup S_i,\P\cap\bigcup S_i )$.

\subsection{Disc-marked extension and equipotentials}

For a \rs $F:\E\to \LLL$ \cm, and any $\LLL$-piece $S$, the above marking $\P_S$ makes
each complementary disc of $S$ into a marked disc. We will use $\kappa_S$ to denote
the potential function of these complementary marked discs of $S$.
Let $S_1$ and $S_2$ be $\CCC$-pieces with
$F(E_1)=S_2$ where $E_1$ is the unique complex $\E$-piece contained in $S_1$. There are many ways
to extend $F|_{E_1}$ to a branched covering. We choose the following one to rigidify the
extension.

\REFLEM{extension} Let $S_1$ and $S_2$ be $\CCC$-pieces with
$F(E_1)=S_2$ where $E_1$ is the unique complex $\E$-piece contained in $S_1$.
Let $\rho$ be a positive function defined on the set of
Jordan curves in $\partial S_1$. Then there is a quasi-regular
branched covering extension $h:\cbar_{S_1}\to
\cbar_{S_2}$ of $F|_{E_1}$ such that: \begin{enumerate}\item[
(a)] $h(\cbar_{S_1}\smm E_1)=\cbar_{S_2}\smm S_2$. \item[
(b)] $h(\P_{S_1})\subset\P_{S_2}$ and the critical values of $h$
are contained in $\P_{S_2}$. \item[
(c)] For any Jordan curve $\g\subset\partial S_1$, $h(\g)$ is an
equipotential in a complementary marked disc of $S_2$, with potential $\kappa_{S_2}(h(\g))=\rho(\g)$. \item[
(d)] $h$ is holomorphic in $\cbar_{S_1}\smm S_1$. \end{enumerate}
Such a map $h$ will be called a {\bf disc-marked extension} of $F|_{E_1}$ associated to the function $\rho$. \ENDLEM

\beginp Let $\al$ be a boundary component of $E_1$, bounding a unique
complementary component $\De_\al$ of $E_1$. Then $\eta:=F(\al)$
is a boundary curve of $S_2$, and bounds a unique complementary
marked disc $(\De_{\eta},b)$ of $S_2$. Set
$d:=\deg(F:\al\to\eta)$.

Note that $\De_\al$ may contain zero or one complementary
component of $S$. In the former case, define $h_{\al}:
\De_{\al}\to\De_{\eta}$ to be a \qc\ map if $d=1$, or a
quasi-regular map with a unique critical value $b$ if $d>1$, such
that $h_{\al}|_{\al}=F|_{\al}$.

In the latter case $\al$ is homotopic within $S_1\smm\P$ to a
unique boundary curve $\g$ of $S_1$. Let $\De_{\g}$ be the
component of $\cbar_{S_1}\smm S_1$ enclosed by $\g$. Then
$\De_{\g}\subset\subset\De_{\al}$, and $\De_{\g}$ together with the
marked point $a\in\De_{\g}$ is a complementary marked disc of
$S_1$.

Let $\eta_1$ be the equipotential in the marked disc $(\De_{\eta}, b)$ with potential
$\kappa_{S_2}(\eta_1)=\rho(\g)$. Denote
by $\De_1$ the disc enclosed by $\eta_1$ and contained in
$\De_{\eta}$. Define $h_{\g}: \De_{\g}\to\De_1$ by
$h_{\g}(z)=\varphi_1^{-1}\circ (\varphi(z))^d$, where $\varphi$
(resp. $\varphi_1$) is a conformal map from the marked disc
$(\De_{\g}, a)$ (resp. $(\De_1, b)$) onto the unit disc $\D$ with
$\varphi(a)=0$ (resp. $\varphi_1(b)=0$). Then there is a
quasi-regular covering $h_{\al\g}$ from
$\De_{\al}\smm\overline{\De}_{\g}$ onto
$\De_{\eta}\smm\overline{\De}_1$ so that
$h_{\al\g}|_{\al}=F|_{\al}$ and
$h_{\al\g}|_{\g}=h_{\g}|_{\g}$. Set $h_{\al}:=h_{\g}$ on
$\De_{\g}$ and $h_{\al}:=h_{\al\g}$ on
$\De_{\al}\smm\overline{\De}_{\g}$. Then $h_{\al}:
\De_{\al}\to\De_{\eta}$ is also quasi-regular, in particular,
holomorphic in $\De_{\g}$.

The map $F|_{E_1}$ together with $h_{\al}$ for all boundary curves of $E_1$ forms
a quasi-regular branched covering $h:\cbar_{S_1}\to \cbar_{S_2}$.
It satisfies the conditions (a)-(d). \qed

\subsection{Spherical holomorphic models}

Consider the same marking $\P_S$ (thus the function $\kappa_S$) for each $\LLL$-piece $S$ as above.

\REFLEM{pp0} Let $\rho$, resp. $\sigma$, be two positives functions defined respectively on the set of
Jordan curves in $\bigcup_{i=1}^p \partial S_i$, resp. in $\partial
S_1$. For $1\le i\le p$,  let $h_i:
\cbar_i\to\cbar_{i+1}$ be a disc-marked extension of $F: E_i\to
S_{i+1}$ associated to the function $\rho$, given by Lemma \ref{extension}.

Then there are pairs of \qc\ homeomorphisms $(\Phi_i, \Psi_i)$ of
$\cbar_i$ onto a distinct copy $\cbar'_i$ of the Riemann sphere, and holomorphic maps $R_i:
\cbar'_i\to\cbar'_{i+1}$ ($i=1,\cdots, p$), such that they satisfy the following conditions:
\begin{enumerate}
\item[(1)] $\Psi_i$ is isotopic to $\Phi_i$ rel $\P_{S_i}$, and $\Phi_i$ is holomorphic on
 $\cbar_i\smm S_i$,
    $i=1,\cdots, p$.
 \item[(2)] $R_i\equiv\Phi_{i+1}\circ h_i\circ\Phi_i^{-1}$ for $2\le i\le p$ (with $\Phi_{p+1}=\Phi_1$),
    and $R_1\equiv \Phi_2\circ h_1\circ\Psi_1^{-1}$ (see \Ref{psi}).
\item[(3)] For any $i=1,\cdots,p$, and  $f_i:=R_{i-1}\circ\cdots\circ R_1\circ R_p\circ\cdots\circ
    R_i$, we have $f_i=\Phi_i\circ h_{i-1}\circ\cdots\circ h_1\circ h_p\circ\cdots\circ
    h_i\circ\Psi_i^{-1}$ and $f_i$ is a \pf\ rational map. The conformal conjugacy class of each $f_i$ depends only on $F$ and $S_i$,
    but not on the choices of the markings, nor on the functions $\rho$ and $\sigma$, nor on $h_i,\Phi_i,\Psi_i$.
\item[(4)]
    For each Jordan curve $\g\subset\partial S_1$, the curve $\Phi_1(\g)$ is
    an equipotential of $f_1$ with potential $\kappa(\Phi_1(\g))=\sigma(\g)$.\end{enumerate}
    \ENDLEM

    Consequently we have:
    \REFCOR{Others}
(1) Fix $2\le i\le p$. The map $\Phi_i$ is holomorphic in $\cbar_i\smm S_i$.
    For each Jordan curve $\g\subset\partial S_i$,
    and for $\beta_\g$ the unique curve in $\partial E_i$ homotopic to $\g$
    within $S_i\smm\P$, both $\Phi_i(\g)$ and $\Phi_i(\beta_\g)$ are
    equipotentials in the same marked Fatou component of $f_i$. And their potentials are related as follows:
    \REFEQN{wp}\kappa(\Phi_i(\g))=\frac{\rho(\g)}{\deg(F|_{\beta_\g})}+\kappa(\Phi_i(\beta_\g)), \quad \kappa(\Phi_i(\beta_\g))=\frac{\kappa(\Phi_{i+1}\circ F(\beta_\g))}{\deg(F|_{\beta_\g})}.\ENDEQN \\
(2) For each Jordan curve $\g\subset\partial S_1$, and for
    $\beta_\g$ the curve in $\partial E_1$ homotopic to $\g$ within
    $S_1\smm\P$, the curve $\Psi_1(\beta_\g)$ is an equipotential in the
    marked Fatou component of $f_1$ that contains $\Phi_1(\g)$, with potential
    \REFEQN{wp2}\kappa(\Psi_1(\beta_\g))=\frac{\kappa(\Phi_2\circ
    F(\beta_\g))}{\deg(F|_{\beta_\g})}.\ENDEQN
\ENDCOR

See the following commutative diagram.
$$\begin{array}{ccccccccccc}
E_1 && E_2 && E_3 &&&& E_p && E_1 \\
\cap & \stackrel{F}{\searrow} & \cap &\stackrel{F}{\searrow} &
\cap &&&& \cap & \stackrel{F}{\searrow} & \cap \\
S_1 && S_2 && S_3 &&&& S_p && S_1 \\
\cap&&\cap &&\cap &&&&\cap &&\cap \\
\cbar_1 & \stackrel{h_1}{\longrightarrow}& \cbar_2 &
\stackrel{h_2}{\longrightarrow} & \cbar_3 & \longrightarrow &
\cdots & \longrightarrow & \cbar_p &
\stackrel{h_p}{\longrightarrow} & \cbar_1 \\
*\downarrow\Psi_1&& \Psi_2\downarrow\Phi_2 &&
\Psi_3\downarrow\Phi_3 &&&&
\Psi_p\downarrow\Phi_p && \Psi_1\downarrow\Phi_1 \\
\cbar'_1 & \stackrel{R_1}{\longrightarrow}& \cbar'_2 &
\stackrel{R_2}{\longrightarrow} & \cbar'_3 & \longrightarrow &
\cdots & \longrightarrow & \cbar'_p &
\stackrel{R_p}{\longrightarrow} & \cbar'_1
\end{array} $$

{\noindent\em Proof of Lemma \ref{pp0}}.\\ Denote by $H:E\to S_1$ the renormalization of $F$ relative
to $S_1$. Set
$$h:=h_p\circ\cdots\circ h_2\circ h_1: \cbar_1\to\cbar_1.$$
Then $h(\P_{S_1})\subset\P_{S_1}$ and $\P_h\subset\P_{S_1}$.
Clearly, $(h, \P_{S_1})$ is a marked extension of the renormalization
$H:E\to S_1$.

It is easy to see that the c-equivalence class of $(h, \P_{S_1})$ does not depend
on the choice of extensions.

Now the assumption that $(F,\P)$ has no renormalized obstructions implies that
$(H,\P\cap S_1)$ as a \rs, has no Thurston obstructions.
This in turn implies that $(h, \P_{S_1})$ has an orbifold distinct from $(2,2,2,2)$
and has no Thurston obstructions, as follows:
 Since the marked points in $\cbar_1\smm S_1$ map to themselves by $h$,
they are eventually $h$-periodic. Let $b$ be a periodic marked
point in $\cbar_1\smm S_1$ with period $k\ge 1$. Denote by $\De_b$
the component of $\cbar_1\smm S_1$ that contains the marked point
$b$ and $\g_b:=\partial\De_b$. Then there is unique component of
$h^{-k}(\g)$, denoted by $\beta$, such that $\beta$ is homotopic to
$\g$ rel $\P_{S_1}$. Note that $\g$ is contained in the boundary multicurve $Y$ and $\beta$ is a component of $F^{-kp}(\g)$ in $S_1$.
Thus the assumption  $\la(W_Y)<1$
implies that
$$\deg(F^{kp}: \beta\to\g)=\deg(h^k:
\beta\to\g)=\deg_b h^k>1.$$
This implies  that $h$ has a periodic critical point (in the cycle of $b$). Therefore
$(h, \P_{S_1})$  has an orbifold distinct from $(2,2,2,2)$. Now any multicurve within $\cbar\smm \P_{S_1}$
can be represented by a multicurve within $S_1\smm \P_{S_1}=S_1\smm (\P\cap S_1)$. So its $(h,\P_{S_1})$-transition
matrix is equal to its $(H,\P\cap S_1)$-transition matrix, thus has the same leading eigenvalue, which is less than one.
This implies $(h,\P_{S_1})$ has no Thurston obstructions.

{\noindent\bf Applying Thurston Theorem to get $(\phi,\psi)$ and $f_1$}. We can then apply Thurton Theorem
\ref{Thurston} to obtain a pair of quasi-conformal maps $(\phi,
\psi)$ from $\cbar_1$ onto $\cbar'_1$ and a rational map $f_1$, whose conformal conjugacy class depends only on the
c-equivalence class of
$(h,P_{S_1})$ (which depends only on  $(H,\P\cap S_1)$), such that $\psi$ is isotopic to $\phi$ rel
$\P_{S_1}$ and $f_1=\phi\circ h\circ\psi^{-1}$. In particular $f_1$ does not depend on the choice
of the functions $\rho$ and $\sigma$.
Moreover, $\P_{f_1}\subset\phi(\P_{S_1})$.

As any periodic cycle $(b)$ of marked points in $\cbar\smm S_1$ contains a critical point of $h$, the cycle $(\phi(b))$ is a
superattracting periodic cycle for $f_1$. Consequently, for every
marked point $a$ in $\cbar_1\smm S_1$, $\phi(a)$ is an eventually
superattracting periodic point of $f_1$.

{\noindent\bf From $(\phi,\psi)$ to $(\Phi_1,\Psi_1)$}.
For every marked point $a$ in $\cbar_1\smm S_1$, denote by $\De_a$
the component of $\cbar_1\smm S_1$ that contains the point $a$ and
$\g_a=\partial\De_a$. Denote by $\eta_a$ the equipotential
 of the Fatou component of $f_1$ containing $\phi(a)$ (with $\phi(a)$ as a marked point), with potential
$\kappa(\eta_a)=\sigma(\g_a)$. Then there is a
quasi-conformal map $\Phi_1$ in the isotopy (rel $\P_{S_1}$) class of
$\phi$ such that for every marked point $a$ in $\cbar_1\smm S_1$, we have $\Phi_1(\g_a)=\eta_a$
(this is because $\g_a$, resp.
$\eta_a$, is peripheral around the point $a\in\P_{S_1}$, resp.
the point $\phi(a)\in\phi(\P_{S_1})$). Moreover, $\Phi_1$ can be taken
to be holomorphic on $\bigcup_a\De_a=\cbar_1\smm S_1$.

As $\Phi_1$ is isotopic to $\phi$ rel $\P_{S_1}$, there is a \qc\ map
$\Psi_1:\cbar_1\to \cbar'_1$ such that it is isotopic to $\psi$ rel
$\P_{S_1}$ and $\Phi_1\circ h\circ\Psi_1^{-1}=f_1$.

{\noindent\bf Getting (in order) $\Phi_p,R_p, \Phi_{p-1},R_{p-1},\cdots, \Phi_2,R_2$ and then $R_1$}. This is illustrated
in the following diagrams:
\REFEQN{psi}\begin{array}{lclccclcc}
\cbar_2 & \stackrel{h_2}{\longrightarrow}& \cbar_{3} &
\longrightarrow & \cdots & \longrightarrow & \cbar_p &
\stackrel{h_p}{\longrightarrow} & \cbar_1 \\
\downarrow\Phi_2 && \downarrow\Phi_{3} &&&&
\downarrow\Phi_p && \downarrow\Phi_1 \\
\cbar'_2 & \stackrel{R_2}{\longrightarrow}& \cbar'_{3} &
\longrightarrow & \cdots & \longrightarrow & \cbar'_p &
\stackrel{R_p}{\longrightarrow} & \cbar'_1
\end{array} \quad \text{and}\quad \begin{array}{rcl}
\cbar_1 & \stackrel{h_1}{\longrightarrow} & \cbar_2\\
\Psi_1\downarrow && \downarrow \Phi_2\\
\cbar'_1 & \stackrel{R_1}{\longrightarrow} & \cbar'_2
\end{array}
\ENDEQN
More precisely pull-back the complex structure of $\cbar'_1$ to $\cbar_p$ by
$\Phi_1\circ h_p$, we have a \qc\ map $\Phi_p:\cbar_p\to\cbar'_p$
such that $R_p:=\Phi_1\circ h_p\circ\Phi_p^{-1}$ is holomorphic.

As a disc-marked extension, we know that $h_p$ is holomorphic in
$\cbar_p\smm S_p$ whose $h_p$-image is contained in $\cbar_1\smm
S_1$. Combining with the result that $\Phi_1$ is holomorphic in
$\cbar_1\smm S_1$ and the equation $R_p\circ\Phi_p=\Phi_1\circ
h_p$, we see that $\Phi_p$ is holomorphic in $\cbar_p\smm S_p$.

Inductively, for $i=p-1, \cdots, 2$, we have a \qc\ map
$\Phi_i:\cbar_i\to\cbar'_i$ such that $R_i:=\Phi_{i+1}\circ
h_i\circ\Phi_i^{-1}$ is holomorphic and $\Phi_i$ is holomorphic in
$\cbar_i\smm S_i$.

Set finally $R_1:=\Phi_2\circ h_1\circ\Psi_1^{-1}$. Then
$R_p\circ\cdots\circ R_2\circ R_1=f_1$. Therefore $R_1$ is also
holomorphic and $\Psi_1$ is holomorphic in $\cbar_1\smm S_1$.

{\noindent\bf Getting $\Psi_i$ and $f_i$}.
As a disc-marked extension, we know that the critical values of
$h_i$ is contained in $\P_{S_{i+1}}$ and
$h_i(\P_{S_i})\subset\P_{S_{i+1}}$ for $1\le i\le p$. Because
$\Psi_1$ is isotopic to $\Phi_1$ rel $\P_{S_1}$, there is a \qc\
map $\Psi_p:\cbar_p\to \cbar'_p$ such that $\Psi_p$ is isotopic to
$\Phi_p$ rel $\P_{S_p}$ and $\Psi_1\circ h_p=R_p\circ\Psi_p$.
Inductively, there is a \qc\ map $\Psi_i:\cbar_i\to \cbar'_i$ for
$i=p-1,\cdots, 2$, such that $\Psi_i$ is isotopic to $\Phi_i$ rel
$\P_{S_i}$ and $\Psi_{i+1}\circ h_i=R_i\circ\Psi_i$. Set then
$f_i:=R_{i-1}\circ\cdots\circ R_1\circ
R_p\circ\cdots\circ R_i$. Now we have
the following commutative diagrams:
$$\begin{array}{rccccrclcccl}
\cbar_i & \stackrel{h_i}{\longrightarrow}& \cdots &&
\stackrel{h_p}{\longrightarrow} & \cbar_1 &
\stackrel{h_1}\longrightarrow &\cbar_2 &\longrightarrow &\cdots &
\stackrel{h_{i-1}}{\longrightarrow} & \cbar_i \\
\Psi_i\downarrow && & \Psi_p\downarrow && \Psi_1\downarrow && \downarrow\Phi_2 &&&&
\downarrow\Phi_i \\
\cbar'_i & \stackrel{R_i}{\longrightarrow}&\cdots &&
\stackrel{R_p}{\longrightarrow} & \cbar'_1 &
\stackrel{R_1}\longrightarrow &\cbar'_2 &\longrightarrow &\cdots &
\stackrel{R_{i-1}}{\longrightarrow} & \cbar'_i
\end{array} \quad \text{and}\quad \begin{array}{rcl}
\cbar_i & \stackrel{*}{\longrightarrow} & \cbar_i\\
\Psi_i\downarrow && \downarrow \Phi_i\\
\cbar'_i & \stackrel{f_i}{\longrightarrow} & \cbar'_i
\end{array}$$
It is easy to see that $f_i=\Phi_i\circ h_{i-1}\circ\cdots\circ
h_1\circ h_p\circ\cdots\circ h_i\circ\Psi_i^{-1}$ for $i\ge 2$.
The above formula shows that $f_i$ is c-equivalent to
$h_{i-1}\circ\cdots\circ h_1\circ h_p\circ\cdots\circ h_i$, which
is \pf. So $f_i$ is also \pf\ and
$\P_{f_i}\subset\Phi_i(\P_{S_i})$. Clearly it is c-equivalent
to a marked extension of the renormalization relative to $S_i$.
Again its conformal conjugacy class depends only on $F$ and $S_i$.
 \qed

{\noindent\em Proof of Corollary \ref{Others}}.
Notice that $f_{i+1}\circ R_i=R_i\circ f_i$, i.e., $R_i$ is a
holomorphic (semi-)conjugacy from $f_i$ to $f_{i+1}$ (set
$f_{p+1}=f_1$). It is classical that their Julia sets are related
by $\JJJ(f_i)=R_i^{-1}(\JJJ(f_{i+1}))$. Note that the critical
values of $R_i$ are contained in $\Phi_{i+1}(\P_{S_{i+1}})$, which
is eventually periodic under $f_{i+1}$. We see that $R_i$ maps
equipotentials of $f_i$ to equipotentials of $f_{i+1}$.

As a disc-marked extension, for each Jordan curve
$\g\subset\partial S_p$, the curve $h_p(\g)$ lies on an equipotential in a
complementary marked disc of $S_1$. Because each Jordan curve in
$\partial\Phi_1(S_1)$ lies on an equipotential of $f_1$ and
$\Phi_1$ is holomorphic in $\cbar_1\smm S_1$, the curve $h_p(\g)$ goes to an
equipotential of $f_1$ by $\Phi_1$. This equipotential of $f_1$ is
pulled back by $R_p$ to equipotentials of $f_p$. Thus $\Phi_p(\g)$
lies on an equipotential of $f_p$. Inductively, we have that each
Jordan curve in $\partial\Phi_i(S_i)$ lies on an equipotential of
$f_i$ for $i=1,\cdots, p$.

Similarly, each curve in $\Phi_i(\partial E_i)$ lies on an
equipotential of $f_i$ for $i\ge 2$ and each curve in
$\Psi_1(\partial E_1)$ lies on an equipotential of $f_1$.

Fix $i\in \{1,\cdots,p\}$. For each Jordan curve $\g\subset\partial S_i$, and for $\beta_\g$
the curve in $\partial E_i$ homotopic to $\g$ within
$S_i\smm\P$, we have that $h_i(\beta_\g)=F(\beta_\g)$ is a curve in $\partial
S_{i+1}$. Note that $\Phi_{i+1}\circ
h_i(\beta_\g)=R_i\circ\Phi_i(\beta_\g)$ if $i\ne 1$ (with $\Phi_{p+1}=\Phi_1$) or $\Phi_2\circ
h_1(\beta_\g)=R_1\circ\Psi_1(\beta_\g)$ if $i=1$. Their potentials
are related  by:
$$\kappa(\Phi_i(\beta_\g))=\frac{\kappa(\Phi_{i+1}\circ
F(\beta_\g))}{\deg(F|_{\beta_\g})}\text{ if }i\ne 1\quad\text{or}\quad
\kappa(\Psi_1(\beta_\g))=\frac{\kappa(\Phi_2\circ
F(\beta_\g))}{\deg(F|_{\beta_\g})}\text{ if }i=1.
$$

Fix now $2\le i\le p$. By the construction of $h_i$ in Lemma \ref{extension}, the curve $h_i(\g)$ is an equipotential with potential $\rho(\g)$
in a complementary marked disc of $S_{i+1}$. We have $\text{mod}
h_i(A(\g, \beta_\g))=\rho(\g)$, where $A(\g, \beta_\g)$ is the annulus
between them. Notice that $\Phi_{i+1}$ is conformal in
$\cbar_{i+1}\smm S_{i+1}$. We also have $\text{mod}
\Phi_{i+1}\circ h_i(A(\g, \beta_\g))=\rho(\g)$. From the equation
$R_i\circ\Phi_i=\Phi_{i+1}\circ h_i$, we get
$$\kappa(\Phi_i(\g))-\kappa(\Phi_i(\beta_\g))=\text{mod}
\Phi_i(A(\g,\beta_\g))=\frac{\rho(\g)}{\deg(F|_{\beta_\g})}.$$  \qed

{\noindent\bf Remark 1}. For every $i$,  if we make a normalization by
requiring that three given distinct points in $\P_{S_i}$ (note that
$\#\P_{S_i}\ge 3$ since $S_i$ is a $\CCC$-piece) go to $(0, -1,
\infty)$ under the action of $\Phi_i$, then $f_i$ is uniquely
determined, as well as the homotopy class (rel $\P_{S_i}$) of $\Phi_i$.

{\noindent\bf Remark 2}. For a $F_*$-periodic cycle $(S_1, \cdots, S_p)$,
we have $p$ renormalizations (one for each $S_i$). Lemma \ref{pp0}
shows that none of them has Thurston obstructions if one of them
has no Thurston obstructions.

\subsection{Proof of Theorem \ref{Pp}} Fix now the positive functions $\sigma$ and $\rho$ as follows:
\REFEQN{prescribe} \forall \g\subset\partial S_1,\ \sigma(\g):=u(\g);\quad  \forall \g\subset\bigcup_{i=1}^p\partial
S_i,\
\rho(\g):=\left(u(\g)-\frac{u(F(\beta_\g))}{\deg(F|_{\beta_\g})}\right)
\deg(F|_{\beta_\g}),\ENDEQN
 where $\beta_\g$ is the curve in $\bigcup_{i=1}^p\partial E_i$
homotopic to $\g$ within $\LLL\smm\P$. Note that $\rho(\g)>0$ for every $\g$ by
the assumption $D_*u<u$.

Let  $(\Phi_i, \Psi_i,R_i,f_i)_{i=1,\cdots,p}$ be the collection of  maps derived from Lemma \ref{pp0}
with the functions $\rho$ and $\sigma$ defined above. Set $f_{S_i}:=f_i$.

Let $\g$ be a Jordan curve in $\partial S_1$. Then
$\kappa(\Phi_1(\g))=\sigma(\g)=u(\g)$ by Lemma \ref{pp0} (4).

Let $\g$ be a Jordan curve in $\partial S_p$ and $\beta_\g$ be the
curve in $\partial E_p$ homotopic to $\g$ within $S_p\smm\P$.
We have $$
\kappa(\Phi_p(\g))\stackrel{\Ref{wp}}{=}\frac{\rho(\g)}{\deg(F|_{\beta_\g})}+\frac{\kappa(\Phi_1\circ
F(\beta_\g))}{\deg(F|_{\beta_\g})}\stackrel{Lem. \ref{pp0}.(4)}{=}\frac{\rho(\g)}{\deg(F|_{\beta_\g})}+
\frac{u(F(\beta_\g))}{\deg(F|_{\beta_\g})} \stackrel{\Ref{prescribe}}{=}u(\g).$$ Inductively, for $i=p-1,\cdots, 2$, we
have $\kappa(\Phi_i(\g))=u(\g)$ for any Jordan curve
$\g\subset\partial S_i$. Therefore $\kappa(\Phi_i(\g))=u(\g)$ for any $i$ and any $\g\subset \partial S_i$.

Fix any $i\in \{1,\cdots,p\}$. Let $\beta$ be a curve in $\partial E_i$ non-null-homotopic within $S_i\smm \P$. By
(\ref{wp}) and (\ref{wp2}), we have
$$\kappa(\Phi_i(\beta))=\frac{\kappa(\Phi_{i+1}\circ
F(\beta))}{\deg(F|_{\beta})}=\frac{u(F(\beta))}{\deg(F|_{\beta})}\quad\text{
if }i\ne 1 \quad \text{and} \quad \kappa(\Psi_1(\beta))=\frac{u(F(\beta))}{\deg(F|_{\beta})}\quad\text{ if
}i=1. $$

Let $\g$ be a Jordan curve in $\partial S_1$ and $\beta_\g$ be the
Jordan curve in $\partial E_1$ homotopic to $\g$ within
$S_1\smm\P$. From the above formula and the fact that $D_*u<u$, we deduce
that $\kappa(\Psi_1(\beta_\g))<\kappa(\Phi_1(\g))$. This implies that
$\Psi_1(E_1)\subset\subset\Phi_1(S_1)$.

For $i=2,\cdots, p$, set  $\phi_{S_i}=\psi_{S_i}=\Phi_i$.
Set also $\phi_{S_1}=\Phi_1$. Obviously, (a)-(d) hold for $i\ge 2$ by the above computation.
 Now we want to define
$\psi_{S_1}$.

Notice that $\Psi_1(E_1)\subset\subset\Phi_1(S_1)$ and $\Psi_1$ is
isotopic to $\Phi_1$ rel $\P_{S_1}$. This implies that for each Jordan
curve $\g\subset\partial S_1$, both $\g$ and
$\Psi_1^{-1}\circ\Phi_1(\g)$ are contained in the same disk
component $\De$ of $\cbar_1\smm E_1$, and are homotopic within
$\De\smm\{a\}$, where $a$ is the unique point of $\P_{S_1}$ in
$\De$. Therefore there is a \qc\ map $\eta$ of $\cbar_1$  isotopic
to the identity rel $\P_{S_1}$ so that $\eta|_{E_1}=id$ and
$\eta=\Psi_1^{-1}\circ\Phi_1$ on $\cbar_1\smm  S_1$. Set
$\widehat \psi:=\Psi_1\circ\eta$. Then $\widehat\psi $ is isotopic to $\Psi_1$ therefore to
$\Phi_1$ rel $\P_{S_1}$, with $\widehat\psi |_{\cbar_1\smm
S_1}=\Phi_1|_{\cbar_1\smm S_1}$ and $\widehat\psi |_{E_1}=\Psi_1$. This gives already (b).

To get (a), i.e. $\widehat\psi $ is isotopic to
$\Phi_1$ rel $\P_{S_1}\cup\partial S_1$, we need to modify $\widehat\psi $ on the annuli of $S_1\smm E_1$.

Set $\chi=\Phi_1^{-1}\circ \widehat\psi$. Then by a purely topological argument (see Lemma \ref{topology})
there is a homeomorphism $T$ which is the identity outside $S_1\smm E_1$ such that
$\chi\circ T$ is isotopic to the identity rel $\P_{S_1}\cup(\cbar_1\smm S_1)$.
Set now $\psi_{S_1}=\widehat\psi\circ T$. We get both (a) and (b) of the Theorem. Point (d) for $i=1$
is also derived from the above computation.
\qed

\section{Proof of Theorem \ref{sb}}\label{I}

Let now $(\E\stackrel{F}{\to}\LLL,\P)$ be a marked \rs\ \cm, without boundary obstructions nor renormalized obstructions. We prove here that $(F,\P)$ is c-equivalent to a holomorphic model.  Decompose $\LLL$ into $\OOO\sqcup
\AAA\sqcup\RRR\sqcup\CCC$ as in \Ref{S}.

\subsection{Choice of the positive vector}

Let $Y$ be a boundary multicurve of $(F,\P)$. It can be chosen to be the
collection of Jordan curves in $\partial\CCC\cup
\partial_+\A$. By assumption, for the $(F,\P)$-transition matrix $W_Y$,
we have $\la(W_Y)<1$. Applying Lemma \ref{matrix}, we have a
positive vector $v\in \R^Y$ so that $W_Yv<v$, i.e. there is a
positive function  \REFEQN{v}v: Y\to\R^{+}\quad \text{such that}\quad
(W_Yv)_{\g}=\sum_{\eta\in Y}\sum_{\al\sim\g}
\frac{v(\eta)}{\deg(F:\al\to\eta)} <v(\g), \ENDEQN where the
last sum is taken over all curves $\al$ in $F^{-1}(\eta)$ that
are homotopic to $\g$ within $\LLL\smm\P$. Let $C>0$ be a constant to be determined later. Denote by ${\bf 1}$ the vector whose every entry is $1$.
Choose $M>0$ to be a large number so that $W_Y(Mv)+C{\bf 1}< Mv$, i.e. \REFEQN{v1}\forall \ \g\in Y,\quad \sum_{\eta\in Y}\sum_{\al\sim\g}
\frac{Mv(\eta)}{\deg(F:\al\to\eta)}+C<Mv(\g). \ENDEQN
 For any $\g\in Y$ the quantity
$Mv(\g)$ will be the
prescribed potential for $\phi_S(\g)$, with $S$ the $\LLL$-piece admitting $\g$ as a boundary curve.

\subsection{Definition of $(\phi_S, \psi_S)$ and $f_S$ for $\CCC$-pieces}

Assume at first that $S_1, \cdots, S_p$ are $\CCC$-pieces with
$F(E_i)=S_{i+1}$ and $S_{p+1}=S_1$, where $E_i$ is the unique
complex $\E$-piece in $S_i$. Set $u:=Mv|_{\bigcup_i\partial S_i}$. We have $D_*u<u$ for $D_*$ the $(F,\P)$-transition
matrix of the set of the boundary curves of $S_1\cup \cdots\cup S_p$. We construct $\phi_{S_i},\psi_{S_i},R_{S_i}$
and $f_{S_i}$ according to Theorem \ref{Pp} for $i=1,\cdots, p$. We do so for every periodic cycle of $F_*$.

Assume now that $S$ is a non-$F_*$-periodic $\CCC$-piece. Then there
are $\CCC$-pieces $S=:S_{-k}, S_{-k+1},$ $\cdots, S_0$ ($k> 0$) such
that $S_i$ is not $F_*$-periodic for $i< 0$ but $S_0$ is
$F_*$-periodic, and $F(E_i)=S_{i+1}$ for $i< 0$, where $E_i$ is
the unique complex $\E$-piece in $S_i$.

Denote by $\cbar_i=\cbar_{S_i}$ for  simplicity. As $S_0$ is
$F_*$-periodic, we have already constructed a \qc\ map $\phi_{S_0}:\cbar_0\to\cbar'_0$ and
a \pf\ rational map $f_0$ on $\cbar_0'$ such that they satisfy the
conditions of Theorem \ref{Pp} for $u=Mv$.

For $i=-1,-2,\cdots,-k$, let $h_i: \cbar_i\to\cbar_{i+1}$ be
a disc-marked extension of $F: E_i\to S_{i+1}$, given by Lemma \ref{extension}, associated to the
function
$$\rho(\g):=\left(Mv(\g)-\frac{Mv(F(\beta_\g))}{\deg(F|_{\beta_\g})}\right)
\deg(F|_{\beta_\g}),$$ where $\g$ is a Jordan curve in $\partial S_i$
and $\beta_\g$ is the curve in $\partial E_i$ homotopic to $\g$ within
$S_i\smm\P$. As before, there are \qc\ maps $\phi_{S_i}=\phi_i:\cbar_i\to \cbar_i'$ and holomorphic maps $R_i:\cbar'_i\to\cbar'_{i+1}$ such
that the following diagram commutes:
$$\begin{array}{lclcccl}
\cbar_{-k} & \stackrel{h_{-k}}{\longrightarrow}& \cbar_{-k+1} &
\stackrel{h_{-k+1}}{\longrightarrow} & \cdots &
\stackrel{h_{-1}}{\longrightarrow} & \cbar_0  \\
\downarrow\phi_{-k} && \downarrow\phi_{-k+1} &&&&
\downarrow\phi_0=\phi_{S_0} \\
\cbar'_{-k} & \stackrel{R_{-k}}{\longrightarrow}& \cbar'_{-k+1} &
\stackrel{R_{-k+1}}{\longrightarrow} & \cdots &
\stackrel{R_{-1}}{\longrightarrow} & \cbar'_0
\end{array} $$
Because $h_i(\P_{S_i})\subset\P_{S_{i+1}}$ and every critical
value (if exists) of $h_i$ lies on $\P_{S_{i+1}}$, we have
$R_i(\phi_i(\P_{S_i}))\subset\phi_{i+1}(\P_{S_{i+1}})$, and every
critical value of $R_{-1}\circ\cdots\circ R_{i}$ lies on
$\phi_0(\P_{S_0})$.

Set $f_{S_i}:=  R_{-1}\circ\cdots\circ R_{i}$. Let $b\in\cbar_i\smm S_i$ be a marked point.
Then $f_{S_i}\circ\phi_{i}(b)$ is the
center of a marked Fatou component $\De$ of $f_0$. The component
$\De_{\phi_{i}(b)}$ of $f_{S_i}^{-1}(\De)$ that contains $\phi_{i}(b)$ is a disc. We will
call $(\De_{\phi_{i}(b)}, \phi_{i}(b))$ a {\em canonical marked
disc}.

The name 'canonical' means that up to a M\"{o}bius transformation,
the configuration formed by these marked discs is uniquely
determined. Note that when a disc-marked extension $h_i$ is
chosen, up to a M\"{o}bius transformation, $\phi_{i}$ is uniquely
determined by $\phi_{i-1}$. As $\phi_{i-1}$ varies in its
homotopy class, $\phi_i$ varies simultaneously in its homotopy
class while $R_i$ remains unchanged. On the other hand various choices of
disc-marked extensions are related by \qc\ maps. More precisely,
if $\tilde h_i$ is another choice of the disc-marked extension,
then there is a \qc\ map $\xi$ of $\cbar_i$ isotopic to the identity rel $\P_{S_i}$, such that $\tilde
h_i=h_i\circ\xi$. Now set $\tilde \phi_{i}=\phi_{i}\circ\xi$, we
get the same holomorphic map $R_i$. This implies that the
maps $f_i$ are in dependent of the extensions (but may depend on the marking).
In particular, the canonical marked discs are independent of the
large number $M$ involved in the function $\rho$ (therefore involved in the extensions $h_i$).

\REFLEM{np} With the assumption above, for any $i=-k,\cdots, -1$, there are \qc\ maps
$\psi_{S_i}=\phi_{S_i}:\cbar_{S_i}\to \cbar_{S_i'}$ such that: \\
(1) $R_i:=\phi_{S_{i+1}}\circ h_i\circ\phi_{S_i}^{-1}$ is holomorphic and is independent of $M$. \\
(2) For any marked point $b\in\P_{S_i}\smm S_i$,
 denote by $\g_b$ the component of $\partial S_i$ that separates
 $b$ from $S_i\smm\g_b$ and by $\al_b$ the component of $\partial E_i$
 that separates $b$ from $E_i\smm\al_b$. Then both $\phi_{S_i}(\g_b)$
 and $\phi_{S_i}(\al_b)$ are equipotentials in the canonical marked
 disc $(\De_{\phi_{S_i}(b)}, \phi_{S_i}(b))$ with potentials
 \REFEQN{wp4} \kappa(\phi_{S_i}(\g_b))=Mv(\g_b)\quad\text{and}\quad
 \kappa(\phi_{S_i}(\al_b))=\kappa(\psi_{S_i}(\al_b)) =\frac{Mv(F(\al_b))}{\deg(F|_{\al_b})}.
 \ENDEQN
\ENDLEM

\beginp (1) is obvious. The proof of (2) is quite easy by following the same argument as before.
\qed

\subsection{Definition of $(\phi_S,\psi_S)$ for other $\LLL$-pieces}

Define $\phi_S=\psi_S=id$ for all $\OOO\cup\RRR$-pieces. Assume that $S$
is an $\A$-piece. Then it is a closed annulus and one of its
boundary curve, say $\g$, is contained in the boundary multicurve
$Y$. We define $\phi_S$ to be a \qc\ map from $S$ to a round
annulus $S'$ in $\cbar_{S'}$ with modulus \REFEQN{mod0}
\mod\phi_S(S)=Mv(\g).\ENDEQN

We will define a map $\psi_A$ for all annular components $A$ of $\LLL\smm \E^m$,
including the $\A$-pieces.
For this we decompose $\E$ into $\E^m\sqcup\E^2\sqcup \E^0$
as follows: \\
$\bullet$  $\E^m$ is the union of complex $\E$-pieces; \\
$\bullet$  $\E^0$ is the union of $\E$-pieces which are contained
in a disk $D\subset\LLL$ with $\#(D\cap\P)\le 1$; \\
$\bullet$  $\E^2$ is the union of $\E\smm\E^0$-pieces which are
contained in an annulus $A\subset\LLL$ with $A\cap\P=\emptyset$.

Clearly, the above three sets are mutually disjoint.
Topologically, $\E^m\subset\CCC$ and $\E^2\subset\CCC\cup\A$.
Dynamically, $F^{-1}(\OOO\cup\RRR)\subset\E^0$ and
$F^{-1}(\A)\subset \E^2\cup\E^0$.

\begin{figure}[htbp]
\begin{center}
\input{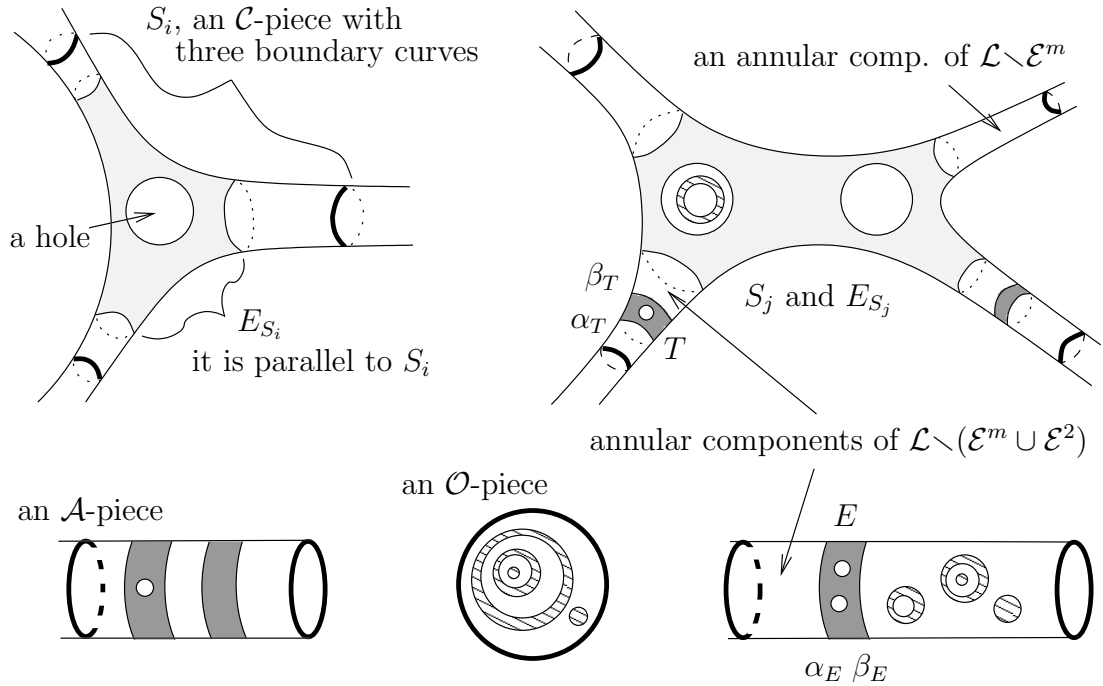}
\end{center}
\vspace{-0.5cm} \caption{\it The $\LLL$-pieces are bounded by thick
curves. Light grey ones are  $\E^m$-pieces, darker-greys, for example $E$ and $T$, are
$\E^2$-pieces. Hatched ones are $\E^0$-pieces (they may appear in
any, necessarily disc or annular, component of $\LLL\smm (\E^m\cup
\E^2)$, and may be nested in each others holes).}
\label{structure}
\end{figure}

 See Figures \ref{structure}.

{\noindent\bf 1. Definition of an auxiliary map $\varphi_E$ for $\E^2$-pieces}.

Assume that $E$ is a $\E^2$-piece and $S$ is an $\LLL$-piece with
$E\subset S$. Set ${\widehat S}:=F(E)$. Then both $S$ and ${\widehat S}$ are contained
in $\CCC\cup\AAA$. Decompose $\E^2$ into $\E^{(2,2)}\sqcup
\E^{(2,m)}$ so that $F(\E^{(2,2)})\subset\AAA$ and
$F(\E^{(2,m)})\subset\CCC$.

If ${\widehat S}$ is an $\AAA$-piece, then there is a \qc\ map $\varphi_E$
from $E$ onto a closed round annulus such that $\phi_{{\widehat S}}\circ
F\circ\varphi_E^{-1}$ is holomorphic in the interior of
$\varphi_E(E)$.

Let $\g$ be one of the two boundary curves in $\partial {\widehat S}$ with
$\g\in Y$. Then there is a Jordan curve $\beta$ in $\partial E$ so that
$F(\beta)=\g$. From (\ref{mod0}), we have: \REFEQN{mod1}
\mod\varphi_E(E)=\frac{\mod\phi_{{\widehat S}}({\widehat S})}{\deg
F|_E}=\frac{Mv(F(\beta))}{\deg (F|_{\beta})}. \ENDEQN

Now assume ${\widehat S}$ is a $\CCC$-piece. Then there is a \qr\ branched
covering $h_E: \cbar_S\to\cbar_{{\widehat S}}$ such that $h_E|_E=F|_E$,
$h_E(E^c)={\widehat S}^c$ and every critical value of $h_E$ is contained
in $\P_{{\widehat S}}$. As before, we have a \qc\ map $\varphi_E$ of $\cbar$
such that $R_E:=\phi_{{\widehat S}}\circ h_E\circ\varphi_E^{-1}$ is
holomorphic from $\cbar_{S'}$ to $\cbar_{{\widehat S}'}$.
$$\begin{array}{rcrclc} \cbar_S & \supset S & \supset  E & \stackrel{\varphi_E}{\longrightarrow} & \varphi_E(E)&\subset \cbar_{S'}\\
h_E\downarrow\ \ \  && F\downarrow && \downarrow R_E& R_E\downarrow \\
\cbar_{\widehat S} &\supset & \widehat S & \underset{\phi_{\widehat S}}{\longrightarrow} & \widehat S'&\subset \cbar_{\widehat S'}
  \end{array},\quad \begin{array}{ccc} \al,\beta & \stackrel{\varphi_E}{\longrightarrow} & \subset \De_{a'}, \De_{b'}
  \\
 F\downarrow && \downarrow R_E\\
   & \underset{\phi_{\widehat S}}{\longrightarrow} & \subset \De_a, \De_b
  \end{array}$$

Note that $\partial E$ has exactly two
boundary curves $\al$ and $\beta$ that are non-null-homotopic within $S\smm\P$. They are homotopic
to each other within $S\smm\P$. From Theorem
\ref{Pp} and Lemma \ref{np}, we know that $\phi_{{\widehat S}}\circ F(\al)$
(resp. $\phi_{{\widehat S}}\circ F(\beta)$) is an equipotential, in a marked
disc $(\De_a, a)$ (resp. $(\De_b, b)$) of the \pf\ rational map
$f_{{\widehat S}}$ when ${\widehat S}$ is $F_*$-periodic, or in a canonical marked
disc, denoted also by $(\De_a, a)$ (resp. $(\De_b, b)$) otherwise,
whose potentials are
$$\kappa(\phi_{{\widehat S}}\circ F(\al))=Mv(F(\al)),\quad
\kappa(\phi_{{\widehat S}}\circ F(\beta))=Mv(F(\beta)). $$ Let
$\De_{a'}$ (resp. $\De_{b'}$) be the
component of $R_E^{-1}(\De_a)$ (resp. $R_E^{-1}(\De_b)$) that contains
$\varphi_E(\al)$ (resp. $\varphi_E(\beta)$). Then $\De_{a'}$ and $\De_{b'}$ are disjoint discs
since neither $\De_a\smm \{a\}$ and $\De_b\smm \{b\}$ contains  critical
values of $R_E$. Set $a':=\De_{a'}\cap R_E^{-1}(a)$ and
$b':=\De_{b'}\cap R_E^{-1}(b)$. Then
$(\De_{a'}, a')$ and $(\De_{b'}, b')$ are
disjoint marked discs in $\cbar_{S'}$. Moreover they are independent of the choice of $M$,
because $(\De_a, a)$ and
$(\De_b,b)$ are independent of the choice of $M$.

Clearly, $\varphi_E(\al)$ and $\varphi_E(\beta)$ are equipotentials
with potentials
$$\kappa(\varphi_E(\al))=\frac{Mv(F(\al))}{\deg
F|_{\al}}\quad\text{and}\quad
\kappa(\varphi_E(\beta))=\frac{Mv(F(\beta))}{\deg F|_{\beta}}.$$

Let $A(E)=A(\al,\beta)$ denote the annulus bounded by $\al$ and $\beta$.
Applying Lemma \ref{Inequality}, there is a constant $C(E)>0$
which is independent of the choice of $M$, such that
\REFEQN{mod2} \frac{Mv(F(\al))}{\deg
F|_{\al}}+\frac{Mv(F(\beta))}{\deg F|_{\beta}}
\le\text{mod}\varphi_E(A(\al,\beta))\le\frac{Mv(F(\al))}{\deg F|_{\al}}+
\frac{Mv(F(\beta))}{\deg F|_{\beta}}+C(E).\ENDEQN

{\noindent\bf The constant $C$}. The set  $\E^{(2,m)}$ has only finite many
pieces $E$ with $C(E)$ independent of the choice of the
number $M$. Set $C:=\sum_E C(E)$. It is also independent of $M$.

{\noindent\bf 2. Embedding of $\varphi_E(E)$ and construction of $\psi_A$}.

Every $\E^2$-piece $E$ is contained in $\AAA$ or in an annular
component of $\CCC\smm\E^m$. We will embed $\varphi_E(E)$ into the
interior of $\AAA\cup(\CCC\smm\E^m)$ so that they are mutually
disjoint.

Assume that $S$ is an $\AAA$-piece. Let $\g$ be a boundary curve
of $S$ with $\g\in Y$. From (\ref{mod1}) and (\ref{mod2}), we have
$$\sum_{E\subset S\cap\E^{(2,2)}}\text{mod }\varphi_E(E)+
\sum_{E\subset S\cap\E^{(2,m)}} \text{mod }\varphi_E ({A(E)})
\le\sum_\beta\frac{Mv(F(\beta))}{\deg(F|_{\beta})}+C, $$ where the last
sum is taken over all the curves $\beta$ in $F^{-1}(\eta)$ for every
$\eta\in Y$ such that $\beta$ is homotopic to $\g$ within $S=S\smm \P$.

The right term is less than $Mv(\g)=\text{mod}(\phi_S(S))$ by
(\ref{v1}). Therefore, as in the non-renormalizable case, one can embed holomorphically
$\varphi_E(E)$ essentially into the interior of $\phi_S(S)$ for
every $\E^2$-piece $E\subset S$ according to the original order of
their non-null-homotopic boundary curves, so that they are
mutually disjoint. In other words, we have a
\qc\ map $\psi_S$ from $S$ onto $\phi_S(S)$, such that \\
$\bullet$ $\psi_S|_{\partial S}=\phi_S|_{\partial S}$ and
$\psi_S$ is isotopic to $\phi_S$ rel $\partial S$; \\
$\bullet$ for every $\E^2$-piece $E\subset S$,
$\varphi_E\circ\psi_S^{-1}$ is holomorphic in the interior of
$\psi_S(E)$.

Consequently, we have \\
$\bullet$ $\phi_{{\widehat S}}\circ F\circ\psi_S^{-1}$ is holomorphic in the
interior of $\psi_S(E)$ for every $\E^2$-piece $E\subset S$ with
${\widehat S}:=F(E)$.

Assume now that $S$ is an $\CCC$-piece and that $A$ is an annular component of $S\smm E_S$ where
$E_S$ is the unique complex $\E$-piece contained in $S$.  Following a similar argument as above, we
have a \qc\ map $\psi_A$ from $A$ onto $\psi_S(A)$, such that \\
$\bullet$ $\psi_A|_{\partial A}=\psi_S|_{\partial A}$ and $\psi_A$
is isotopic to $\psi_S|_A$ rel $\partial A$; \\
$\bullet$ $\phi_{{\widehat S}}\circ F\circ\psi_A^{-1}$ is holomorphic in the
interior of $\psi_A(E)$ for every $\E^2$-piece $E\subset A$ with
${\widehat S}:=F(E)$.

\subsection{Definition of $\theta_S$}

Define $\theta_S=\phi_S^{-1}\circ\psi_S$ for every $\AAA$-piece
$S$. If $S$ is a $\CCC$-piece, define
$$\theta_S=\left\{\begin{array}{l}\phi_S^{-1}\circ\psi_A\quad\text{ on
every annular component $A$ of }S\smm\E^m; \\
\phi_S^{-1}\circ\psi_S\quad\text{ otherwise.}\end{array}\right.$$
 Then
$\theta_S|_{\partial S}=id$ and $\theta_S$ is isotopic to the
identity rel $\partial S\cup(S\cap\P)$. Moreover, for every
$\E^2\cup\E^m$-piece $E$ with $E\subset S$ and $F(E)=\widehat S$, the map
$\phi_{\widehat S}\circ F\circ\theta_{S}^{-1}\phi_{S}^{-1}$ is
holomorphic in the interior of $\phi_{S}\theta_{S}(E)$.

Now if $\E^0\cup\OOO\cup\RRR=\emptyset$, the proof of Theorem
\ref{sb} is already completed. Otherwise one can follow the
argument as in the non-renormalizable case (there is no more trouble in case
$\RRR\neq\emptyset$) to modify $\theta_S$ on
$S\smm (\E^2\cup\E^m\cup\partial S)$ with the help of a suitable layer function. This ends the proof of Theorem
\ref{sb}. \qed

\section{A combination result}\label{comb}

A {\em regular puzzle} is by definition a subset of $ \cbar$ which is also a \ps.

A {\em regular open set} is by definition the complement of a
regular puzzle.

Let $U,V$ be regular open sets in $\cbar$ with $V\subset\subset
U$. Let $G:U\to V$ be a \qr\ branched covering. We say that
$(G,U,V)$ is a {\em locally holomorphic attracting system}, if
there is a finite set
$\P'\subset U$ such that:\\
$\bullet$ $G(\P')=\P'$;\\
$\bullet$ $G$ is holomorphic in a neighborhood of $\P'$ and each cycle in $\P'$ is (super)attracting;\\
$\bullet$ for any $z\in V$ the limit set of $\{G^n(z)\}$ is
contained in $\P'$.

Let $F:\E\to \LLL$ be a \rs\ \cm, in particular $F$ is quasi-regular.
We say that $F$ has no {\em analytization obstruction} is it has
no boundary obstruction, and , for each renormalization $H:E\to S$
(if any, and not necessarily \pf), either
\\
(1) $\#\P_f\cap S<\infty$ and $(H,\P_F\cap S)$ as a \rs\ has no Thurston obstructions; or
\\
(2) for the integer $p$ such that $H=F^p|_E$, each step of the
composition
$$E\overset{F}{\longrightarrow} F(E) \overset{F}{\longrightarrow} F^2(E)\overset{F}{\longrightarrow}
\cdots \overset{F}{\longrightarrow}  F^{p-1}(E)
\overset{F}{\longrightarrow} S$$ is holomorphic in the interior.

What we have proved in this paper can be reformulate in the
following stronger form:

\REFTHM{infinite} Let $G$ be a \qr\ branched covering of $\cbar$
with degree at least $2$. Assume that $\cbar=V\sqcup \LLL$ is a
splitting with
$\LLL$ a regular puzzle such that:\\
(a) $G^{-1}(V)\supset\supset V$;\\
(b) $(G, G^{-1}(V),V)$ is a locally holomorphic attracting system;\\
(c) $G:G^{-1}(\LLL)\to\LLL$ is a \rs\ \cm\ without analytization
obstructions.

Let $K$ be the union of the filled Julia set $K_H$ of each of the
holomorphic renormalizations.

\noindent Then there is a rational map $g$ and a pair of
qc-homeomorphisms $\phi,\psi$ of $\cbar$ such that
\\
$\bullet$ $\phi\circ G=g\circ \psi$;\\
$\bullet$ $\psi$ is isotopic to $\phi$ rel $\P_G\cup K$;\\
$\bullet$ the Beltrami coefficient of $\phi$ is equal to $0$
almost everywhere on $K$. \ENDTHM

\appendix

\section{Non-negative  matrices}

For a vector $v=(v_i)\in \R^n$ we write $v>0$ if every coordinate $v_i$
is strictly positive.

\REFLEM{matrix} Let $D=(a_{ij})$ be a real square matrix with
$a_{ij}\ge 0$ for each entry $a_{ij}$. Denote by $\la$ its
spectral radius, i.e. the maximal modulus of the eigenvalues. Then
$\la<1$ iff there is a vector $v>0$ such that $Dv<v$.\ENDLEM

\beginp
The following proof is provided by H.H. Rugh. Necessity: Assume
$v>0$ and $Dv<v$. Then $Dv\le av$ for some $0\le a< 1$. Define a
norm on the underlying vector space by $\|x\|=\sum_i ( v_i\cdot
|x_i|)$. Then, writing $|x|$ as the vector whose $i$-th entry is
$|x_i|$, we have  $$\|\,^t\!Dx\|=\,^tv\,^t\!D|x| =\,^t(Dv)|x|\le a
\,^tv|x|=a\|x\|\ .$$ Therefore, $\displaystyle{\la:=\max_{\la' \
\text{eigenvalue of }D}|\la'| = \max_{\la' \ \text{eigenvalue of
}\,^t\!D}|\la'| \le \|\,^t\!D\|\le a \ .}$

Sufficiency: Now assume $\la <1$. By continuity of the spectral
radius, there is $\epsilon >0$ such that the spectral radius
$\la_\epsilon$ of $D+\epsilon:=(a_{ij}+\epsilon)$ satisfies
$\la_\epsilon<1$. Now the Perron-Frobenius Theorem assures that
$\la_\epsilon$ is also an eigenvalue (called the leading
eigenvalue) and it has a strictly positive eigenvector $v>0$. So
$Dv\le(D+\epsilon)v = \la_\epsilon v < v \ .$ \qed

Note that it follows that $\la$ is also an eigenvalue of $D$
(called the leading eigenvalue). Lemma \ref{matrix} actually gives
an equivalent definition of the eigenvalues.

\REFCOR{le} Let $\la(D)$ be the leading eigenvalue of a
non-negative square matrix $D$. Then
$$\la(D)=\inf\{\la|\ \exists\ v>0\text{ such that }Dv<\la v\}.$$
\ENDCOR

\REFCOR{compare} Assume that $A$ and $B$ are non-negative
$n\times n$ matrix with $A\le B$ (i.e. each entry of $A$ is less
than or equal to the corresponding entry of $B$), then
$\la(A)\le\la(B)$. \ENDCOR

\beginp From Lemma \ref{matrix}, we see that for any
$\la_0>\la(B)$, there is a  vector $v>0$ so that
$\la_0^{-1}Bv<v$. Thus $\la_0^{-1}Av\le\la_0^{-1}Bv<v$. Again by
Lemma \ref{matrix}, we have $\la(A)<\la_0$ for any constant
$\la_0>\la(B)$. So $\la(A)\le\la(B)$. \qed

Let $A$ be an $n\times n$ matrix with a block decomposition
$$\left(\begin{array}{ccc}B_{11} & \cdots & B_{1k} \\
\vdots && \vdots \\ B_{k1} & \cdots & B_{kk}\end{array}\right)$$
where $B_{ij}$ is an $n_i\times n_j$ matrix (in particular each $B_{ii}$ is a square matrix). We say that the block
decomposition is {\em projected}\/ if for each $B_{ij}$, there is
a number $b_{ij}$ such that the summation of each column of $B_{ij}$
is equal to $b_{ij}$.

This property could be understood as the following: An $n\times n$
matrix can be considered as a linear map of $\R^n$ defined by the
left action:
$$
\left(\begin{array}{c} v_1 \\ \vdots \\ v_n\end{array}\right)
\mapsto A\left(\begin{array}{c} v_1 \\ \vdots \\
v_n\end{array}\right).$$ According to the block decomposition of
$A$, there is a corresponding decomposition of the index set
$I=\{1,\cdots, n\}$ by $ I=I_1\sqcup\cdots\sqcup I_k $ with $\#
I_i=n_i$. Define a linear projection $\pi: \R^n\to\R^k$ by
$$(\pi v)_i=\sum_{\de\in I_i}v_{\de}.$$

\REFLEM{proj} There is a $k\times k$ matrix $B$ such that
$\pi\circ A=B\circ\pi$ if and only if the block decomposition
$A=(B_{ij})$ is projected. In this case, $B=(b_{ij})$.
$$\begin{array}{rcl} \R^n & \stackrel{A}{\longrightarrow} & \R^n
\\ \pi\downarrow & & \downarrow\pi \\ \R^k & \stackrel{B}{\longrightarrow} &
\R^k \end{array} $$  \ENDLEM
\beginp
Set $A=(a_{\de\be})$. For any $v\in\R^n$,
$$(\pi\circ Av)_i=\sum_j\sum_{\be\in I_j}\left(\sum_{\de\in
I_i}a_{\de\be}\right)v_{\be},\quad \text{and}\quad
 (B\circ\pi(v))_i=\sum_j\sum_{\be\in I_j} b_{ij}v_{\be}\ .$$
If the block decomposition is projected, then for $\be\in I_j$,
$\sum_{\de\in I_i}a_{\de\be}=b_{ij}$. Therefore $\pi\circ
Av=B\circ\pi(v)$. Conversely, assume that $\pi\circ A=B\circ\pi$.
For $\be\in I_j$, let $e_{\be}\in\R^n$ be a vector whose
$\be$-entry is $1$ and $0$ elsewhere. Then $(\pi\circ
Ae_{\beta})_i=b_{ij}$, and $(B\circ\pi(e_{\be}))_i=\sum_{\de\in
I_i}a_{\de\be}$. So for $\be\in I_j$, $\sum_{\de\in
I_i}a_{\de\be}=b_{ij}$, i.e. the block decomposition is projected.
\qed

\REFTHM{mp} Assume that $A$ is a non-negative square matrix with
a projected block decomposition $A=(B_{ij})$. Set $B=(b_{ij})$.
Then $\la(A)=\la(B)$. \ENDTHM

\beginp
Let $v\neq 0$ be an eigenvector of $A$ for the leading eigenvalue
$\la(A)$, i.e. $Av=\la(A)v$. Set $u=\pi(v)$. Then $Bu=\pi\circ
Av=\pi(\la(A)v)=\la(A)\pi(v)=\la(A)u$ by the above Lemma. So
$\la(A)$ is an eigenvalue of $B$ and hence $\la(A)\le\la(B)$ since
the leading eigenvalue is the maximum of the eigenvalues.

Conversely, let $u\neq 0$ be an eigenvector of the transpose $B^t$
of $B$ for the leading eigenvalue $\la(B)$ (note that $B$ and
$B^t$ have same leading eigenvalues), i.e. $B^tu=\la(B)u$. Set
$v=(v_{\be})\in\R^n$ by $v_{\be}:=u_j$ for $\be\in I_j$. Then for
$\de\in I_i$,
$$ (A^t v)_{\de}=\sum_j\sum_{\be\in I_j} a_{\be\de} v_{\be}
=\sum_j b_{ji} v_j=(B^tu)_i=\la(B)u_i=\la(B)v_{\de}.$$ So $\la(B)$
is an eigenvalue of $A^t$. Therefore we again have
$\la(B)\le\la(A)$. \qed

\REFCOR{mp1} Let $A'$ be a non-negative square matrix with a block decomposition
$(B'_{ij})$. Assume that for each $ij$, the summation of each column of $B'_{ij}$
is at most $b_{ij}$. Set $B=(b_{ij})$. Then $\la(A')\le \la(B)$.\ENDCOR
\beginp For each $ij$, we just need to replace one entry of each column of $B'_{ij}$
by a larger number so that the summation of the column becomes exactly $b_{ij}$. Denote by $B_{ij}$
the modified matrix. Set $A=(B_{ij})$. Then $\la(A)=\la(B)$ by Theorem \ref{mp}
and $\la(A')\le \la(A)$ by Corollary \ref{compare}.\qed
\section{Reversing the Gr\"otzsch inequality}

A {\em equipotential} $\g$ in a marked disc $(\De,a)$ is a curve
mapped onto a round circle under a conformal representation
$\varphi:(\De,a)\to (\D ,0)$. The potential of $\g$ is defined to
be the modulus of the annulus between $\partial\De$ and $\g$.

\REFLEM{Inequality} Let $(D_i,z_i)$, $i=1,2$ be two disjoint
marked hyperbolic discs. Then there is a constant $C>0$
independent of $v_1>0,v_2>0$ such that, for the annulus
$A(v_1,v_2)$ between the equipotential in $D_1$ of potential $v_1
$ and the equipotential of $D_2$ of potential $v_2$, we have
$$ v_1+v_2\le\!\!\!\!\!\mod(A(v_1,v_2))\le v_1+v_2+C $$ \ENDLEM

\beginp
The left hand side is just the Gr\"otzsch inequality.

The conformal radius of a marked disc $(\De,0)$ is defined to be
the radius $r$ if there is a conformal map $\varphi: (\De,0)\to
(D(0,r),0)$ with $\varphi'(0)=1$. And the conformal radius of a
marked disc $(\widetilde{\De},\infty)$ is defined to be the
conformal radius of $(\pi(\widetilde{\De}),0)$ with $\pi(z)=1/z$.

Let $\xi$ be a M\"obius transformation of $\cbar$ with
$\xi(z_1)=0$ and $\xi(z_2)=\infty$. Any two such maps differ by a
multiplicative constant. So the product $C_1\cdot C_2$ of the
conformal radii of $(\xi(D_1),0)$ and  $(\xi(D_2),\infty)$ is
independent of the choice of $\xi$. Denote by $W_i$ the component
of $A(v_1,v_2)^c$ containing $z_i$, $i=1,2$. By Koebe
$1/4$-Theorem, $\xi(W_1)$ contains $\{|z|\le C_1 r_1/4\}$ and
$\xi(W_2)$ contains $\{|z|\ge 4/(C_2 r_2)\}$, where
$r_i=e^{-v_i}$. Therefore\\
$\mod(A(v_1,v_2))\le\log\left( \frac 4{C_2 r_2}\cdot \frac 4{C_1 r_1}\right)
=\log \left(\frac{16}{C_1C_2}\frac1{r_1r_2}\right) =\log
\frac{16}{C_1C_2}+ v_1+v_2\ .$\qed

\section{Quasi-conformal extensions}\label{Extension}
We state here several results about \qc\ maps that have been
frequently used in the paper.

\REFLEM{circle} Let $h:C_1\to C_2$ be a homeomorphism between two
quasi-circles $C_1$ and $C_2$ in $\cbar$. If $h$ can be extended
to a \qc\ map on an one-side neighborhood of $C_1$, then $h$ can
be extended to a global \qc\ homeomorphism of $\cbar$. Moreover
the extension can be chosen to be a diffeomorphism from $\cbar\smm
C_1$ onto $\cbar\smm C_2$.\ENDLEM

\REFLEM{open} Let $\Omega_i\subset \cbar$ ($i=1,2$) be two open
connected domains such that $\partial \Omega_i$ ($i=1,2$) consists
of $p\ge 0$ disjoint quasi circles (we allow the case $p=0$). Let
$\P\subset \Omega_1$ be a finite set (may or may not be empty).
Let $f:\overline \Omega_1\to \overline\Omega_2$ be an orientation
preserving homeomorphism. If, either $p=0$, or $f|_{\partial
\Omega_1}$ can be extended to a \qc\ map on an one-side
neighborhood of each curve of $\partial \Omega_1$, then there is a
\qc\ homeomorphism in the isotopy class of $f$ modulo $\partial
\Omega_1\cup \P$. \ENDLEM

\REFLEM{unit} Let $h:S^1\to S^1$ be an orientation preserving
homeomorphism of the unit circle. Assume that $h$ can be extended
as a \qc\ map $f$ on an inner neighborhood $B$ of $S^1$ (i.e.
$B\supset\{1-\ep < |z|<1\}$ for some $\ep >0$), then $h$ is
quasi-symmetric.\ENDLEM

\beginp Denote by $\mu$ the Beltrami coefficient of $f$. Denote by $\D$ the unit disc.
Let $\nu=\mu$ on B and $\nu=0$ on $\D\smm B$. By the Measurable
Riemann Mapping Theorem, there is a \qc\ homeomorphism $g$ of $\D$
whose Beltrami coefficient is $\nu$. Then $g|_{S^1}$ is
quasi-symmetric. On the other hand, $f\circ g^{-1}$ is holomorphic
on $g(B)$. Therefore $f\circ g^{-1}|_{S^1}:S^1\to S^1$ is
real-analytic, in particular quasi-symmetric. So $h=(f\circ
g^{-1})\circ g|_{S^1}$ is also quasi-symmetric.\qed

{\noindent\em Proof of Lemma \ref{circle}}. Fix $i=1,2$. By
definition of quasi-circles,  there is a \qc\ homeomorphism
$\phi_i$ of $\cbar$ such that $\phi_i(C_i)=S^1$. Furthermore
$\phi_i$ can be chosen to be diffeomorphism on $\cbar\smm C_i$ as
follows: Set $\De=\phi_i^{-1}(\D)$. Let $\psi:\De\to \D$ be a
conformal map. Then $\phi_i\circ \psi^{-1}:\D\to \D$ is a \qc\
homeomorphism. Thus its boundary map is quasi-symmetric. Let
$\eta$ be the Beurling-Ahlfors extension of this boundary map, it
is a diffeomorphism of $\D$. Now $\eta\circ \psi|_{\De}$ is again
a diffeomorphism, whose boundary map is $\phi_i|_{S^1}$.
Set $h_1=\phi_2\circ h\circ \phi_1^{-1}$. Then by Lemma \ref{unit}
this $h_1$ is quasi-symmetric, thus has a \qc\ extension to
$\cbar$. Moreover its extension can be chosen to be a
diffeomorphism outside $S^1$.
Thus $h=\phi_2^{-1}\circ h_1 \circ \phi_1$ can be extended to a
\qc\ homeomorphism of $\cbar$, and a diffeomorphism outside $C_1$.
\qed

{\noindent\em Proof of Lemma \ref{open}}. By Lemma \ref{circle} we
can assume that $\partial \Omega_i$ are smooth Jordan curves and
that $f|_{\partial \Omega_1}$ is a diffeomorphism. Then one can
find a diffeomorphism in its isotopy class rel $\partial
\Omega_1\cup \P$.\qed

\section{A lemma about isotopy}

\REFLEM{topology} Let $(D_i,a_i), i=1,\cdots k$, $k\ge 1$ be finitely many marked Jordan discs in $S^2$, with
disjoint closures. Let $P$ be a closed (or empty) set  contained in
$S^2\smm\bigsqcup \overline{D_i}$.
Assume that $h_1:S^2\to S^2$ is an orientation preserving homeomorphism, and $h: S^2\times [0,1]\to S^2$ is continuous,
such that
\\ a) $h_1|_{P\cup \bigcup D_i}=id$\\
b1) $h(\cdot,t)$ is a homeomorphism for any $t\in [0,1]$;\\
b2) $h(\cdot,0)=id$, $h(\cdot,1)=h_1$;\\
b3) $h(x,t)=x$ for any $x\in P\cup \bigcup\{a_i\}$ and any $t\in [0,1]$.

For $i=1,\cdots, k$, set $\g_i=\partial D_i$, and let $\beta_i$ be a Jordan curve disjoint from $\overline D_i$
 so that the annuli $A_i:=A(\beta_i,\g_i)$ have mutually disjoint closures and are disjoint from $P$.
Then there is a continuous map $H: S^2\times [0,1]\to S^2$ such that
\\
c1) $H(\cdot,t)$ is a homeomorphism for any $t\in [0,1]$;\\
c2) $H(\cdot,0)=id$,  $H(\cdot,1)=h_1\circ T$ with $T=id$ outside $\bigcup A_i$;\\
c3) $H(x,t)=x$ for any $x\in P\cup \bigcup D_i$ and any $t\in [0,1]$.
\ENDLEM
\beginp
 Set $h_t=h(\cdot,t)$, and $E=S^2\smm\bigcup(D_i\cup A_i)$.

 For each $i$ choose a Jordan curve $\al_i$ in $D_i$ bounding a disc $D(\al_i)$ so that  $a_i\in D(\al_i)
 \subset D_i$ and that $h_t(\beta_i)\cap \al_i=\emptyset$
 for any $t\in [0,1]$.

 Define $s:S^2\times [0,1]\to S^2$ continuous, with each $s(\cdot,t):=s_t$ a homeomorphism of $S^2$,  as follows: $s(\cdot,0)=id$.
 \\
 $s(x,t)=\left\{\begin{array}{ll} h_t^{-1}(x) & x\in \bigcup D(\al_i)\\
 \text{interpolation} & x\in \bigcup A(\al_i,\beta_i)\\
 x & x\in E
                \end{array}\right.$.
        Then $s(x,1)=\left\{\begin{array}{ll} x=h_1^{-1}(x) & x\in \bigcup D(\al_i)\\
 T_0(x) & x\in \bigcup A(\al_i,\beta_i)\\
 x & x\in E
                \end{array}\right.$,
where $T_0$ is a certain homeomorphism of $S^2$ that is identity outside $\bigcup A(\al_i,\beta_i)$. \\
Set $\xi_t=h_t\circ s_t$. Then $\xi_0=id$, \quad
 $\xi_t(x)=\left\{\begin{array}{ll} x & x\in \bigcup D(\al_i)\\
 \text{interpolation} & x\in \bigcup A(\al_i,\beta_i)\\
 h_t(x) & x\in E\\
 x=h_t(x) & x\in P\subset E
                \end{array}\right.$,\\
 and $\xi_1(x)=h_1\circ s_1(x)=\left\{\begin{array}{ll} x=h_1(x) & x\in \bigcup D(\al_i)\\
 h_1\circ T_0(x) & x\in \bigcup A(\al_i,\beta_i)\\
 h_1(x) & x\in E\\
 x=h_1(x) & x\in P\subset E
                \end{array}\right.$.

Let now $u:S^2\to S^2$ be a homeomorphism such that $u(D_i)=D(\al_i)$, $u(a_i)=a_i$ for each $i$ and $u|_{E}=id$.
Define then $v:S^2\to S^2$ be a homeomorphism such that $v|_{\bigcup D(\al_i)}=u^{-1}$ and $v|_{h_1(E)\cup E}=id$.
Set
$\zeta_t=v\circ \xi_t\circ u$. We have  $\zeta_0=\left\{\begin{array}{ll} id & x\notin \bigcup A_i
                                                         \\
v\circ u & x\in \bigcup A_i \end{array}\right.$,\\
$\zeta_t(x)=\left\{\begin{array}{ll} x=h_1(x) & x\in \bigcup D_i\\ \text{interpolation}  & x\in \bigcup A_i\\
 v\circ h_t(x) & x\in E\\
 x=h_1(x) & x\in P\subset E
                \end{array}\right.$,\quad $\zeta_1(x)=\left\{\begin{array}{ll} x=h_1(x) & x\in \bigcup D_i\\
 v\circ h_1\circ T_0\circ u(x)  & x\in \bigcup A_i\\
 h_1(x) & x\in E\\
 x=h_1(x) & x\in P\subset E
                \end{array}\right.=h_1\circ T_1(x)$\\
        for a certain homeomorphism $T_1$ of $S^2$ with $T_1=id$ outside $\bigcup A_i$.
Set finally $H_t(x)=\zeta_t\circ \zeta_0^{-1}(x)$. It has the required properties.
In particular $H(\cdot,1)= \left\{\begin{array}{ll} h_1 & x\notin \bigcup A_i
                                                         \\
v\circ h_1\circ T_0\circ v^{-1} & x\in \bigcup A_i \end{array}\right.=h_1\circ T_1\circ T_2$, with $T_2=\zeta_0^{-1}$, and $T_1\circ T_2=id$
outside $\bigcup_iA_i$.\qed
{\small CUI Guizhen, Academy of Mathematics and System Sciences, Chinese
Academy of Sciences,
Beijing 100080, People's Republic of China. Email: gzcui@math.ac.cn.

TAN Lei, Laboratoire d'Analyse, G\'eom\'etrie et Mod\'elisation,
UMR 8088, Universit\'e de Cergy-Pontoise,
2 Avenue Adolphe Chauvin, 95302 Cergy-Pontoise Cedex, France.
Email: tanlei@math.u-cergy.fr.}
\end{document}